\magnification 1200

\input amssym
\input miniltx
\input pictex

\font \bbfive = bbm5
\font \bbseven = bbm7
\font \bbten = bbm10
\font \eightbf = cmbx8
\font \eighti = cmmi8 \skewchar \eighti = '177
\font \eightit = cmti8
\font \eightrm = cmr8
\font \eightsl = cmsl8
\font \eightsy = cmsy8 \skewchar \eightsy = '60
\font \eighttt = cmtt8 \hyphenchar \eighttt = -1

\font \sixbf = cmbx6
\font \sixi = cmmi6 \skewchar \sixi = '177
\font \sixrm = cmr6
\font \sixsy = cmsy6 \skewchar \sixsy = '60
\font \tensc = cmcsc10

\scriptfont \bffam = \bbseven
\scriptscriptfont \bffam = \bbfive
\textfont \bffam = \bbten
\newskip \ttglue

\def \eightpoint {\def \rm {\fam 0 \eightrm }\relax
\textfont 0= \eightrm
\scriptfont 0 = \sixrm \scriptscriptfont 0 = \fiverm
\textfont 1 = \eighti
\scriptfont 1 = \sixi \scriptscriptfont 1 = \fivei
\textfont 2 = \eightsy
\scriptfont 2 = \sixsy \scriptscriptfont 2 = \fivesy
\textfont 3 = \tenex
\scriptfont 3 = \tenex \scriptscriptfont 3 = \tenex
\def \it {\fam \itfam \eightit }\relax
\textfont \itfam = \eightit
\def \sl {\fam \slfam \eightsl }\relax
\textfont \slfam = \eightsl
\def \bf {\fam \bffam \eightbf }\relax
\textfont \bffam = \bbseven
\scriptfont \bffam = \bbfive
\scriptscriptfont \bffam = \bbfive
\def \tt {\fam \ttfam \eighttt }\relax
\textfont \ttfam = \eighttt
\tt \ttglue = .5em plus.25em minus.15em
\normalbaselineskip = 9pt
\def \MF {{\manual opqr}\-{\manual stuq}}\relax
\let \sc = \sixrm
\let \big = \eightbig
\setbox \strutbox = \hbox {\vrule height7pt depth2pt width0pt}\relax
\normalbaselines \rm }

\def \setfont #1{\font \auxfont =#1 \auxfont } \def \withfont #1#2{{\setfont {#1}#2}}  
 \def \ifundef #1{\expandafter \ifx \csname #1\endcsname \relax } \def \undefrule {\kern 2pt \vrule width
5pt height 5pt depth 0pt \kern 2pt} \def \UndefLabels {} \def \possundef #1{\ifundef {#1}\undefrule {\eighttt
#1}\undefrule \global \edef \UndefLabels {\UndefLabels #1\par } \else \csname #1\endcsname \fi } \newcount \secno \secno
= 0 \newcount \stno \stno = 0 \newcount \eqcntr \eqcntr = 0 \ifundef {showlabel}  \fi
\def \define #1#2{\global \expandafter \edef \csname #1\endcsname {#2}} \long \def \error #1{\medskip \noindent {\bf
******* #1}} \def \fatal #1{\error {#1\par Exiting...}\end } \def \advseqnumbering {\global \advance \stno by 1 \global
\eqcntr =0} \def \current {\ifnum \secno = 0 \number \stno \else \number \secno \ifnum \stno = 0 \else .\number \stno
\fi \fi } \begingroup \catcode `\@ =0 \catcode `\\=11 @global @def @textbackslash{\} @endgroup \begingroup \catcode
`\%=11  \fi \ifundef {#1PrimarilyDefined}\define {#1}{#2}\define {#1PrimarilyDefined}{#2}\else \edef \old {\csname
#1\endcsname }\edef \new {#2}\if \old \new \else \fatal {Duplicate definition for label ``{\tt #1}'', already defined as
``{\tt \old }''.}\fi \fi } \def \newlabel #1#2{\define {#1}{#2}} \def \label #1 {\deflabel {#1}{\current }} \def
\equationmark #1 {\ifundef {InsideBlock} \advseqnumbering \eqno {(\current )} \deflabel {#1}{\current } \else \global
\advance \eqcntr by 1 \edef \subeqmarkaux {\current .\number \eqcntr } \eqno {(\subeqmarkaux )} \deflabel
{#1}{\subeqmarkaux } \fi } \def \split #1.#2.#3.#4;{\global \def \parone {#1}\global \def \partwo {#2}\global \def
\parthree {#3}\global \def \parfour {#4}} \def \NA {NA} \def \ref #1{\split #1.NA.NA.NA;(\possundef {\parone }\ifx
\partwo \NA \else .\partwo \fi )} \newcount \bibno \bibno = 0 \def \newbib
#1#2{\define {#1}{#2}} \def \Bibitem #1 #2; #3; #4 \par {\smallbreak \global \advance \bibno by 1 \item {[\possundef
{#1}]} #2, {``#3''}, #4.\par \ifundef {#1PrimarilyDefined}\else \fatal {Duplicate definition for bibliography item
``{\tt #1}'', already defined in ``{\tt [\csname #1\endcsname ]}''.}  \fi \ifundef {#1}\else \edef \prevNum {\csname
#1\endcsname } \ifnum \bibno =\prevNum \else \error {Mismatch bibliography item ``{\tt #1}'', defined earlier (in aux
file ?) as ``{\tt \prevNum }'' but should be ``{\tt \number \bibno }''.  Running again should fix this.}  \fi \fi
\define {#1PrimarilyDefined}{#2}} \def \jrn #1, #2 (#3), #4-#5;{{\sl #1}, {\bf #2} (#3), #4--#5} \def \Article #1 #2;
#3; #4 \par {\Bibitem #1 #2; #3; \jrn #4; \par } \def \references {\begingroup \bigbreak \eightpoint \centerline {\tensc
References} \nobreak \medskip \frenchspacing } \catcode `\@ =11\def \c@itrk #1{{\bf \possundef {#1}}} \def \c@ite
#1{{\rm [\c@itrk {#1}]}} \def \sc@ite [#1]#2{{\rm [\c@itrk {#2}\hskip 0.7pt:\hskip 2pt #1]}} \def \du@lcite {\if \pe@k
[\expandafter \sc@ite \else \expandafter \c@ite \fi } \def \cite {\futurelet \pe@k \du@lcite } \catcode `\@ =12 \def
\Headlines #1#2{\nopagenumbers \headline {\ifnum \pageno = 1 \hfil \else \ifodd \pageno \tensc \hfil \lcase {#1} \hfil
\folio \else \tensc \folio \hfil \lcase {#2} \hfil \fi \fi }} \def \title #1{\medskip \centerline {\withfont
{cmbx12}{\ucase {#1}}}}   \long \def \Quote #1\endQuote {\begingroup \leftskip 35pt \rightskip 35pt \parindent 17pt \eightpoint
#1\par \endgroup } \long \def \Abstract #1\endAbstract {\vskip 1cm \Quote \noindent #1\endQuote }  \def \Authors #1{\bigskip \centerline {\tensc #1}} \def \Note
#1{\footnote {}{\eightpoint #1}} \def \Date #1 {\Note {\it Date: #1.}}  \def \part #1#2{\vfill \eject \null \vskip
0.3\vsize \withfont {cmbx10 scaled 1440}{\centerline {PART #1} \vskip 1.5cm \centerline {#2}} \vfill \eject }  \def \fix {\smallskip
\noindent $\blacktriangleright $\kern 12pt}     \def \ucase #1{\edef \auxvar {\uppercase {#1}}\auxvar } \def \lcase #1{\edef
\auxvar {\lowercase {#1}}\auxvar } \def \section #1 \par {\global \advance \secno by 1 \stno = 0 \goodbreak \bigbreak
\noindent {\bf \number \secno .\enspace #1.}  \nobreak \medskip \noindent } \def \state #1 #2\par {\begingroup \def
\InsideBlock {} \medbreak \noindent \advseqnumbering {\bf \current .\enspace #1.\enspace \sl #2\par }\medbreak \endgroup
} \def \definition #1\par {\state Definition \rm #1\par } \long \def \Proof #1\endProof {\begingroup \def \InsideBlock
{} \medbreak \noindent {\it Proof.\enspace }#1 \ifmmode \eqno \endproofmarker $$ \else \hfill $\endproofmarker $
\looseness = -1 \fi \medbreak \endgroup } \def \$#1{#1 $$$$ #1} \def \explica #1#2{\mathrel {\buildrel \hbox {\sixrm #1}
\over #2}} \def \explain #1#2{\explica {\ref {#1}}{#2}}  \def \=#1{\explain
{#1}{=}} \def \pilar #1{\vrule height #1 width 0pt} \def \stake #1{\vrule depth #1 width 0pt} \newcount \fnctr \fnctr =
0 \def \fn #1{\global \advance \fnctr by 1 \edef \footnumb {$^{\number \fnctr }$}\footnote {\footnumb }{\eightpoint
#1\par \vskip -10pt}} \def \text #1{\hbox {#1}}  
 \def \Item #1{\smallskip \item {{\rm #1}}} \newcount \zitemno \zitemno = 0 \def \izitem
{\global \zitemno = 0} \def \zitemplus {\global \advance \zitemno by 1 \relax } \def \rzitem {\romannumeral \zitemno }
\def \rzitemplus {\zitemplus \rzitem } \def \zitem {\Item {{\rm (\rzitemplus )}}}  \def
\zitemmark #1 {\deflabel {#1}{\rzitem }} \newcount \nitemno \nitemno = 0  \def \nitem
{\global \advance \nitemno by 1 \Item {{\rm (\number \nitemno )}}} \newcount \aitemno \aitemno = -1 \def \boxlet
#1{\hbox to 6.5pt{\hfill #1\hfill }} \def \iaitem {\aitemno = -1} \def \aitemconv {\ifcase \aitemno a\or b\or c\or d\or
e\or f\or g\or h\or i\or j\or k\or l\or m\or n\or o\or p\or q\or r\or s\or t\or u\or v\or w\or x\or y\or z\else zzz\fi }
\def \aitem {\global \advance \aitemno by 1\Item {(\boxlet \aitemconv )}} \def \aitemmark #1 {\deflabel {#1}{\aitemconv
}}  \def \Case #1:{\medskip \noindent {\tensc Case #1:}} \def \<{\left \langle \vrule
width 0pt depth 0pt height 8pt } \def \>{\right \rangle } \def \({\big (} \def \){\big )}  \def
\and {\hbox {\quad and \quad }}   \def \imply {\mathrel {\Rightarrow }} \def \IFF {\kern 7pt\Leftrightarrow \kern 7pt} \def \IMPLY
{\kern 7pt \Rightarrow \kern 7pt} \def \for #1{\quad \forall \,#1} \def \endproofmarker {\square } \def \"#1{{\it #1}\/}
 \def \inv {^{-1}} \def \*{\otimes } \def \caldef #1{\global \expandafter \edef \csname
#1\endcsname {{\cal #1}}} \def \bfdef #1{\global \expandafter \edef \csname #1\endcsname {{\bf #1}}} \bfdef N \bfdef Z
\bfdef C \bfdef R \newlabel {SubshSect}{2} \newlabel {IntroduceShift}{2.1} \newlabel {DefineLanguage}{2.2} \newlabel
{AllSSForbWords}{2.3} \newlabel {SFTOpen}{2.5} \newlabel {OneExit}{3.3} \newlabel {DefStrongExit}{3.4} \newlabel
{EvenExit}{3.5} \newlabel {FirstPaction}{4.1} \newlabel {DescrPAction}{4.2} \newlabel {IntroPRep}{5.1} \newlabel
{CondForPrep}{5.1.1} \newlabel {SgZero}{5.2} \newlabel {FinalSpace}{5.3} \newlabel {IntroT}{5.4} \newlabel
{GenForMatsu}{6.1} \newlabel {GenForMatsuCar}{6.3} \newlabel {DefineCM}{6.4} \newlabel {Zoo}{6.5} \newlabel
{SpectrumSection}{7} \newlabel {ETTensorOne}{7.2} \newlabel {DefinePhiX}{7.3} \newlabel {XIsDense}{7.4} \newlabel
{SubsetModel}{7.5} \newlabel {XiHomeo}{7.6} \newlabel {IntroduceSpec}{7.7} \newlabel {PropForPaint}{7.8} \newlabel
{BigCompos}{7.9} \newlabel {EquivCond}{7.10} \newlabel {CharacXix}{7.11} \newlabel {BinvInXix}{7.12} \newlabel
{PropStem}{7.13} \newlabel {DefineStem}{7.14} \newlabel {StemLeftInverse}{7.15} \newlabel {StemContinuous}{7.16}
\newlabel {BooleanContinuous}{7.16.1} \newlabel {XiXMaximal}{7.17} \newlabel {SuficCritInterior}{7.18} \newlabel
{OntoContinuousSFT}{7.19} \newlabel {DefineTau}{8.1} \newlabel {PactOnAlg}{8.2} \newlabel {SpectralAction}{8.3}
\newlabel {DescrPAOnOmega}{8.4} \newlabel {TBCont}{8.4.1} \newlabel {TwoSystems}{8.5} \newlabel {Equivariance}{8.6}
\newlabel {FirstEquivar}{8.6.1} \newlabel {EquivPact}{8.7} \newlabel {Enfiabilidade}{8.8} \newlabel {CMCrossProd}{9.1}
\newlabel {FromCPtoMats}{9.2} \newlabel {CondexpOnMyElts}{9.3} \newlabel {FromMatstoCP}{9.4} \newlabel
{TensorGuy}{9.4.1} \newlabel {AlgebrasSame}{9.5} \newlabel {TwoSetsCs}{10.1.1} \newlabel {ConjugaT}{10.2.1} \newlabel
{BasisProdTop}{11.1} \newlabel {BaseTopology}{11.2} \newlabel {BigNBD}{11.2.1} \newlabel {LongAlpha}{11.2.2} \newlabel
{TrocaBeta}{11.2.3} \newlabel {HereIsU}{11.2.4} \newlabel {VerySimpleNBDs}{11.3} \newlabel {SimplerNBD}{11.3.1}
\newlabel {TopFreeSect}{12} \newlabel {IntroduceFix}{12.1} \newlabel {OldTofFree}{12.2} \newlabel {IsolatedPoint}{12.3}
\newlabel {EvenNotTopFree}{12.4} \newlabel {FixPoints}{12.5} \newlabel {MainTopFree}{12.6} \newlabel
{FoundBasicOnFP}{12.6.1} \newlabel {ClassicMinimal}{13.1} \newlabel {Reachability}{13.2} \newlabel
{PombaSimplesDiagrama}{13.3} \newlabel {BridgeForEvenShift}{13.4} \newlabel {EvenNotMinimal}{13.5} \newlabel
{CofinalityProps}{13.6} \newlabel {gbOrbit}{13.7} \newlabel {MinimalityExplained}{13.8} \newlabel
{ReachabilityMultiple}{13.9} \newlabel {ColcofOrbX}{13.12} \newlabel {XiXInV}{13.12.1} \newlabel {XisHere}{13.12.2}
\newlabel {Bejetivo}{13.12.3} \newlabel {SomeMembers}{13.12.4} \newlabel {TopFreeWhenColCof}{13.13} \newlabel
{BddCost}{13.16} \newlabel {OrbNoMeet}{13.16.1} \newlabel {CrazyInclusion}{13.17} \newlabel
{CrazyInclusionDynamical}{13.18} \newlabel {InclusionWithX}{13.18.1} \newlabel {InclusionWithXandTheta}{13.18.2}
\newlabel {MainMinimal}{13.19} \newlabel {MainHyper}{13.21} \newlabel {SuperHiperCond}{13.22} \newlabel {AppSection}{14}
\newlabel {AlgAsGpd}{14.2} \newlabel {TopFreeFinal}{14.3} \newlabel {Simplicity}{14.5} \newbib {AbadieGpd}{1} \newbib
{AdelaR}{2} \newbib {ArchSpiel}{3} \newbib {BCFS}{4} \newbib {BrownOzawa}{5} \newbib {CarlsenThesis}{6} \newbib
{CarlsenCuntzPim}{7} \newbib {CarlsenNotes}{8} \newbib {MatsuCarl}{9} \newbib {CarlsenSilvestrov}{10} \newbib
{CuntzKrieger}{11} \newbib {amena}{12} \newbib {ortho}{13} \newbib {NewLook}{14} \newbib {NHausd}{15} \newbib
{PDSFB}{16} \newbib {infinoa}{17} \newbib {ELQ}{18} \newbib {ItoTaka}{19} \newbib {KPR}{20} \newbib {KPRR}{21} \newbib
{Kurka}{22} \newbib {LindMarcus}{23} \newbib {MatsuOri}{24} \newbib {MatsumotoDimension}{25} \newbib
{MatsuAutomorph}{26} \newbib {Matsumoto}{27} \newbib {Parry}{28} \newbib {RenaultCartan}{29} \newbib {Royer}{30} \newbib
{Thomsen}{31} \def \Caixa #1{\setbox 1=\hbox {$#1$\kern 1pt}\global \edef \tamcaixa {\the \wd 1}\box 1} \def \caixa
#1{\hbox to \tamcaixa {$#1$\hfil }} \def \med #1{\mathop {\textstyle #1}\limits }  \def
\medsum {\med \sum }  \def \medcup {\med \bigcup } \def \medcap {\med \bigcap }  \caldef S \caldef O \caldef G \def \MSymbol {{\cal M}} \def \MatsAlg
{\MSymbol _X} \def \CMSymbol {{\cal O}} \def \CM {\CMSymbol _X} \def \TT {\tilde T} \def \shft {S} \def \T {{\bf T}}
\def \Lang {{\cal L}_X} \def \bigTheta {\theta } \def \Forbid {{\cal F}} \def \Folow {F} \def \Cyl {Z} \def \F {{\bf F}}
 \def \Core {{\cal D}_X} \def \SpecAlg {sp(\Core )} \def \Spec {\Omega _X} \def \XO {\Xi } \def
\SO {\Psi } \def \stem {\sigma } \def \emptyword {\varnothing } \def \pr {u} \def \prt {\tilde u} \def \et {\tilde e}
\def \qt #1{\hbox {`$#1$'}} \def \amper {\kern 1.5pt\&\kern 1.5pt} \def \onedot {\hbox {$\scriptstyle .$}} \def \nedots
{\mathinner {\mkern 1mu\raise 2.5pt\vbox {\kern 7pt\onedot }\mkern -2mu \raise 4pt\onedot \mkern -2mu\raise 5.5pt\onedot
\mkern 1mu}} \def \PPinv {\F _+\F _+\inv } \def \Fix {\hbox {\sl Fix}} \def \Orb {\hbox {\sl Orb}}  \def \cost {\hbox {\sl Cost}} \def \colcof {collectively cofinal} \def \colcofty {collective cofinality} \def
\scof {strongly cofinal} \def \scofty {strong cofinality} \def \hcof {hyper cofinal} \def \hcofty {\hcof ity} \def
\ipfimply {\global \def \BreakImply {\kern -5pt}} \def \implica (#1)(#2){(#1)$\,\imply \,$(#2)} \def \pfimply (#1)(#2){
\BreakImply \implica (#1)(#2).\global \def \BreakImply {\medskip \noindent }} \def \book #1#2{\cite [#2]{PDSFB}} \def
\explicaRef #1#2#3{\explica {\sixrm [{\sixbf \possundef {#1}}:#2]}{#3}} \def \begincenter {\begingroup \obeylines
\leftskip =0pt plus 1fill \rightskip =0pt plus 1fill}

\Headlines {Partial actions and  subshifts}{M. Dokuchaev and R. Exel}

\title {Partial actions and  subshifts}

\Authors {M. Dokuchaev and R. Exel}

\Abstract Given a finite alphabet $\Lambda $, and a not necessarily finite type subshift $X\subseteq \Lambda ^\infty $,
we introduce a partial action of the free group $\F (\Lambda )$ on a certain compactification $\Spec $ of $X$, which we
call the spectral partial action.

The space $\Spec $ has already appeared in many papers in the subject, arising as the spectrum of a commutative
C*-algebra usually denoted by $\Core $.  Since the descriptions given of $\Spec $ in the literature are often somewhat
terse and obscure, one of our main goals is to present a sensible model for it which allows for a detailed study of its
structure, as well as of the spectral partial action, from various points of view, including topological freeness and
minimality.

We then apply our results to study certain C*-algebras associated to $X$, introduced by Matsumoto and Carlsen.  Most of
the results we prove are already well known, but our proofs are hoped to be more natural and more in line with
mainstream techniques used to treat similar C*-algebras.  The clearer understanding of $\Spec $ provided by our model in
turn allows for a fine tuning of some of these results, including a necessary and sufficient condition for the
minimality of the Carlsen-Matsumoto C*-algebra $\CM $, generalizing a similar result of Thomsen.  \endAbstract

\section Introduction

The theory of C*-algebras associated to subshifts has a long and exciting history, having been initiated by Matsumoto in
\cite {MatsuOri}, later receiving invaluable contributions from many other authors, notably Carlsen, Silvestrov, and
Thomsen.  Accounts of this history may be found in \cite {MatsuCarl} and \cite {CarlsenCuntzPim}, to which the
interested reader is referred.

For now it suffices to say that, given a finite alphabet $\Lambda $, and a closed subset $X\subseteq \Lambda ^{\bf N}$
which is invariant under the left shift $$ \shft :\Lambda ^{\bf N}\to \Lambda ^{\bf N}, $$ in which case one says that
the pair $(X,\shft )$ is a (one sided) \"{subshift}, Matsumoto initiated a study of certain C*-algebras associated to
$X$, whose algebraic properties reflect certain important dynamical properties of the subshift itself, and whose
$K$-theory groups provide new invariants for subshifts.

The extensive literature in this field (see for instance the list of references in \cite {CarlsenNotes}) contains a lot
of information about the structure of these algebras, such as faithful representations, nuclearity, characterization of
simplicity, computation of $K$-theory groups and a lot more.  It is therefore a perilous task to attempt to add anything
else to the wealth of results currently available, a task we hope to be undertaking in a responsible manner.

The motivation that brought us to revisit the theory of C*-algebras associated to subshifts, and the justification for
writing the present paper, is twofold: firstly we are able to offer a sensible description of a certain topological
space, which we will denote by $\Spec $ (for the cognoscenti, we are referring to the spectrum of the ill-fated
commutative algebra $\Core $ appearing in most papers on the subject), and which has evaded all attempts at analysis,
except maybe for some somewhat obscure projective limit descriptions given in \cite {MatsumotoDimension} and \cite
[Section 2.1]{CarlsenThesis}.  See also \cite{Royer}.

Secondly we will introduce a partial action of the free group $\F =\F (\Lambda )$ on $\Spec $, called the \"{spectral
partial action}, whose associated crossed product is the Carlsen-Matsumoto C*-algebra $\CM $.  Given that description,
we may recover many known results for $\CM $, as well as give the first necessary and sufficient condition for
simplicity which applies for all subshifts, including those where the shift map is not surjective.

Our study of the space $\Spec $ is perhaps the single most important contribution we have to offer.  The method we adopt
is essentially the same one used by the second named author and M.~Laca in the analysis of Cuntz-Krieger algebras for
infinite matrices \cite {infinoa}, the crucial insight being the introduction of a partial action of the free group $\F
$.  To be more precise, for each letter $a$ in the alphabet $\Lambda $, consider the subsets $\Folow _a$ and $\Cyl _a$
of $X$, given by $$ \def \quad { } \matrix { \Folow _a & = & \{y\in X: ay\in X\}, \hbox { and} \hfill \cr \pilar {18pt}
\Cyl _a & = & \{x\in X: x=ay, \hbox { for some infinite word } y\}. } $$ It is then evident that the assignment $$
\theta _a:y\in \Folow _a \to ay\in \Cyl _a $$ is a continuous bijective map.  By iterating the $\theta _a$ and their
inverses, we obtain partially defined continuous bijective maps on $X$, thus forming what is known as a partial action
of the group $\F $ on $X$.

Partial actions on topological spaces are often required to map open sets to open sets, but except for the case of
shifts of finite type, the above partial action does not satisfy this requirement since the $\Folow _a$ may fail to be
open.  One may nevertheless use this badly behaved partial action to build a bona fide partial representation on the
(typically non-separable) Hilbert space $\ell ^2(X)$.  At this point it is evident that we have grossly violated the
topology of $X$, but alas a new commutative C*-algebra is born, generated by the range projections of all partial
isometries in our partial representation, and with it a new topological space is also born, namely the spectrum of said
algebra.

This algebra, usually denoted $\Core $, turns out to be as well known for experts in the field, as it is dreaded.
Nevertheless it plays a well known role in the theory of partial representations and, like many other commutative
algebras arising from partial representations, its spectrum $\Spec $ may be seen as a subspace of the Bernoulli space
$2^\F $, which is moreover invariant under the well known Bernoulli partial action \book {Definition}{5.10}.

Based on the algebraic relations possessed by this representation we may deduce certain special properties enjoyed by
the elements of $2^\F $ which lie in $\Spec $, and if we see $2^\F $ as the set of all subsets of the Cayley graph of
$\F $, such properties imply that these elements must have the aspect of a \"{river basin}, in the sense that there is a
main river consisting of an infinite (positive) word, together with its tributaries.  Not only must the main river\break

\null \hfill \beginpicture \setcoordinatesystem units <0.0080truecm, -0.0080truecm> point at 1000 1000 \put {$\bullet $}
at 0 0 \setdots <0.07cm> \plot 8 8 390 390 / \setsolid \put {$\circ $} at 398 398 \setdots <0.07cm> \plot 406 406 588
588 / \setsolid \put {$\circ $} at 596 596 \setdots <0.07cm> \plot 604 604 687 687 / \setsolid \put {$\circ $} at 695
695 \setdots <0.07cm> \plot 702 702 736 736 / \setsolid \put {$\circ $} at 744 744 \setdots <0.07cm> \plot 751 751 760
760 / \setsolid \put {$\circ $} at 768 768 \setdots <0.07cm> \plot 736 751 727 760 / \setsolid \put {$\circ $} at 719
768 \setdots <0.07cm> \plot 751 736 760 727 / \setsolid \put {$\circ $} at 768 719 \setdots <0.07cm> \plot 687 702 653
736 / \setsolid \put {$\circ $} at 646 744 \setdots <0.07cm> \plot 653 751 662 760 / \setsolid \put {$\circ $} at 670
768 \setdots <0.07cm> \plot 638 751 629 760 / \setsolid \put {$\circ $} at 621 768 \setdots <0.07cm> \plot 638 736 629
727 / \setsolid \put {$\circ $} at 621 719 \setdots <0.07cm> \plot 702 687 736 653 / \setsolid \put {$\circ $} at 744
646 \setdots <0.07cm> \plot 751 653 760 662 / \setsolid \put {$\circ $} at 768 670 \setdots <0.07cm> \plot 736 638 727
629 / \setsolid \put {$\circ $} at 719 621 \setdots <0.07cm> \plot 751 638 760 629 / \setsolid \put {$\circ $} at 768
621 \setdots <0.07cm> \plot 588 604 505 687 / \setsolid \put {$\circ $} at 498 695 \setdots <0.07cm> \plot 505 702 539
736 / \setsolid \put {$\circ $} at 547 744 \setdots <0.07cm> \plot 554 751 563 760 / \setsolid \put {$\circ $} at 571
768 \setdots <0.07cm> \plot 539 751 530 760 / \setsolid \put {$\circ $} at 522 768 \setdots <0.07cm> \plot 554 736 563
727 / \setsolid \put {$\circ $} at 571 719 \setdots <0.07cm> \plot 490 702 456 736 / \setsolid \put {$\circ $} at 448
744 \setdots <0.07cm> \plot 456 751 465 760 / \setsolid \put {$\circ $} at 473 768 \setdots <0.07cm> \plot 441 751 432
760 / \setsolid \put {$\circ $} at 424 768 \setdots <0.07cm> \plot 441 736 432 727 / \setsolid \put {$\circ $} at 424
719 \setdots <0.07cm> \plot 490 687 456 653 / \setsolid \put {$\circ $} at 448 646 \setdots <0.07cm> \plot 441 653 432
662 / \setsolid \put {$\circ $} at 424 670 \setdots <0.07cm> \plot 441 638 432 629 / \setsolid \put {$\circ $} at 424
621 \setdots <0.07cm> \plot 456 638 465 629 / \setsolid \put {$\circ $} at 473 621 \setdots <0.07cm> \plot 604 588 687
505 / \setsolid \put {$\circ $} at 695 498 \setdots <0.07cm> \plot 702 505 736 539 / \setsolid \put {$\circ $} at 744
547 \setdots <0.07cm> \plot 751 554 760 563 / \setsolid \put {$\circ $} at 768 571 \setdots <0.07cm> \plot 736 554 727
563 / \setsolid \put {$\circ $} at 719 571 \setdots <0.07cm> \plot 751 539 760 530 / \setsolid \put {$\circ $} at 768
522 \setdots <0.07cm> \plot 687 490 653 456 / \setsolid \put {$\circ $} at 646 448 \setdots <0.07cm> \plot 638 456 629
465 / \setsolid \put {$\circ $} at 621 473 \setdots <0.07cm> \plot 638 441 629 432 / \setsolid \put {$\circ $} at 621
424 \setdots <0.07cm> \plot 653 441 662 432 / \setsolid \put {$\circ $} at 670 424 \setdots <0.07cm> \plot 702 490 736
456 / \setsolid \put {$\circ $} at 744 448 \setdots <0.07cm> \plot 751 456 760 465 / \setsolid \put {$\circ $} at 768
473 \setdots <0.07cm> \plot 736 441 727 432 / \setsolid \put {$\circ $} at 719 424 \setdots <0.07cm> \plot 751 441 760
432 / \setsolid \put {$\circ $} at 768 424 \setdots <0.07cm> \plot 390 406 208 588 / \setsolid \put {$\circ $} at 200
596 \setdots <0.07cm> \plot 208 604 291 687 / \setsolid \put {$\circ $} at 299 695 \setdots <0.07cm> \plot 306 702 340
736 / \setsolid \put {$\circ $} at 348 744 \setdots <0.07cm> \plot 355 751 364 760 / \setsolid \put {$\circ $} at 372
768 \setdots <0.07cm> \plot 340 751 331 760 / \setsolid \put {$\circ $} at 323 768 \setdots <0.07cm> \plot 355 736 364
727 / \setsolid \put {$\circ $} at 372 719 \setdots <0.07cm> \plot 291 702 257 736 / \setsolid \put {$\circ $} at 249
744 \setdots <0.07cm> \plot 257 751 266 760 / \setsolid \put {$\circ $} at 274 768 \setdots <0.07cm> \plot 242 751 233
760 / \setsolid \put {$\circ $} at 225 768 \setdots <0.07cm> \plot 242 736 233 727 / \setsolid \put {$\circ $} at 225
719 \setdots <0.07cm> \plot 306 687 340 653 / \setsolid \put {$\circ $} at 348 646 \setdots <0.07cm> \plot 355 653 364
662 / \setsolid \put {$\circ $} at 372 670 \setdots <0.07cm> \plot 340 638 331 629 / \setsolid \put {$\circ $} at 323
621 \setdots <0.07cm> \plot 355 638 364 629 / \setsolid \put {$\circ $} at 372 621 \setdots <0.07cm> \plot 192 604 109
687 / \setsolid \put {$\circ $} at 101 695 \setdots <0.07cm> \plot 109 702 143 736 / \setsolid \put {$\circ $} at 150
744 \setdots <0.07cm> \plot 158 751 167 760 / \setsolid \put {$\circ $} at 175 768 \setdots <0.07cm> \plot 143 751 134
760 / \setsolid \put {$\circ $} at 126 768 \setdots <0.07cm> \plot 158 736 167 727 / \setsolid \put {$\circ $} at 175
719 \setdots <0.07cm> \plot 94 702 60 736 / \setsolid \put {$\circ $} at 52 744 \setdots <0.07cm> \plot 60 751 69 760 /
\setsolid \put {$\circ $} at 77 768 \setdots <0.07cm> \plot 45 751 36 760 / \setsolid \put {$\circ $} at 28 768 \setdots
<0.07cm> \plot 45 736 36 727 / \setsolid \put {$\circ $} at 28 719 \setdots <0.07cm> \plot 94 687 60 653 / \setsolid
\put {$\circ $} at 52 646 \setdots <0.07cm> \plot 45 653 36 662 / \setsolid \put {$\circ $} at 28 670 \setdots <0.07cm>
\plot 45 638 36 629 / \setsolid \put {$\circ $} at 28 621 \setdots <0.07cm> \plot 60 638 69 629 / \setsolid \put {$\circ
$} at 77 621 \setdots <0.07cm> \plot 192 588 109 505 / \setsolid \put {$\circ $} at 101 498 \setdots <0.07cm> \plot 94
505 60 539 / \setsolid \put {$\circ $} at 52 547 \setdots <0.07cm> \plot 60 554 69 563 / \setsolid \put {$\circ $} at 77
571 \setdots <0.07cm> \plot 45 554 36 563 / \setsolid \put {$\circ $} at 28 571 \setdots <0.07cm> \plot 45 539 36 530 /
\setsolid \put {$\circ $} at 28 522 \setdots <0.07cm> \plot 94 490 60 456 / \setsolid \put {$\circ $} at 52 448 \setdots
<0.07cm> \plot 45 456 36 465 / \setsolid \put {$\circ $} at 28 473 \setdots <0.07cm> \plot 45 441 36 432 / \setsolid
\put {$\circ $} at 28 424 \setdots <0.07cm> \plot 60 441 69 432 / \setsolid \put {$\circ $} at 77 424 \setdots <0.07cm>
\plot 109 490 143 456 / \setsolid \put {$\circ $} at 150 448 \setdots <0.07cm> \plot 158 456 167 465 / \setsolid \put
{$\circ $} at 175 473 \setdots <0.07cm> \plot 143 441 134 432 / \setsolid \put {$\circ $} at 126 424 \setdots <0.07cm>
\plot 158 441 167 432 / \setsolid \put {$\circ $} at 175 424 \setdots <0.07cm> \plot 406 390 588 208 / \setsolid \put
{$\circ $} at 596 200 \setdots <0.07cm> \plot 604 208 687 291 / \setsolid \put {$\circ $} at 695 299 \setdots <0.07cm>
\plot 702 306 736 340 / \setsolid \put {$\circ $} at 744 348 \setdots <0.07cm> \plot 751 355 760 364 / \setsolid \put
{$\circ $} at 768 372 \setdots <0.07cm> \plot 736 355 727 364 / \setsolid \put {$\circ $} at 719 372 \setdots <0.07cm>
\plot 751 340 760 331 / \setsolid \put {$\circ $} at 768 323 \setdots <0.07cm> \plot 687 306 653 340 / \setsolid \put
{$\circ $} at 646 348 \setdots <0.07cm> \plot 653 355 662 364 / \setsolid \put {$\circ $} at 670 372 \setdots <0.07cm>
\plot 638 355 629 364 / \setsolid \put {$\circ $} at 621 372 \setdots <0.07cm> \plot 638 340 629 331 / \setsolid \put
{$\circ $} at 621 323 \setdots <0.07cm> \plot 702 291 736 257 / \setsolid \put {$\circ $} at 744 249 \setdots <0.07cm>
\plot 751 257 760 266 / \setsolid \put {$\circ $} at 768 274 \setdots <0.07cm> \plot 736 242 727 233 / \setsolid \put
{$\circ $} at 719 225 \setdots <0.07cm> \plot 751 242 760 233 / \setsolid \put {$\circ $} at 768 225 \setdots <0.07cm>
\plot 588 192 505 109 / \setsolid \put {$\circ $} at 498 101 \setdots <0.07cm> \plot 490 109 456 143 / \setsolid \put
{$\circ $} at 448 150 \setdots <0.07cm> \plot 456 158 465 167 / \setsolid \put {$\circ $} at 473 175 \setdots <0.07cm>
\plot 441 158 432 167 / \setsolid \put {$\circ $} at 424 175 \setdots <0.07cm> \plot 441 143 432 134 / \setsolid \put
{$\circ $} at 424 126 \setdots <0.07cm> \plot 490 94 456 60 / \setsolid \put {$\circ $} at 448 52 \setdots <0.07cm>
\plot 441 60 432 69 / \setsolid \put {$\circ $} at 424 77 \setdots <0.07cm> \plot 441 45 432 36 / \setsolid \put {$\circ
$} at 424 28 \setdots <0.07cm> \plot 456 45 465 36 / \setsolid \put {$\circ $} at 473 28 \setdots <0.07cm> \plot 505 94
539 60 / \setsolid \put {$\circ $} at 547 52 \setdots <0.07cm> \plot 554 60 563 69 / \setsolid \put {$\circ $} at 571 77
\setdots <0.07cm> \plot 539 45 530 36 / \setsolid \put {$\circ $} at 522 28 \setdots <0.07cm> \plot 554 45 563 36 /
\setsolid \put {$\circ $} at 571 28 \setdots <0.07cm> \plot 604 192 687 109 / \setsolid \put {$\circ $} at 695 101
\setdots <0.07cm> \plot 702 109 736 143 / \setsolid \put {$\circ $} at 744 150 \setdots <0.07cm> \plot 751 158 760 167 /
\setsolid \put {$\circ $} at 768 175 \setdots <0.07cm> \plot 736 158 727 167 / \setsolid \put {$\circ $} at 719 175
\setdots <0.07cm> \plot 751 143 760 134 / \setsolid \put {$\circ $} at 768 126 \setdots <0.07cm> \plot 687 94 653 60 /
\setsolid \put {$\circ $} at 646 52 \setdots <0.07cm> \plot 638 60 629 69 / \setsolid \put {$\circ $} at 621 77 \setdots
<0.07cm> \plot 638 45 629 36 / \setsolid \put {$\circ $} at 621 28 \setdots <0.07cm> \plot 653 45 662 36 / \setsolid
\put {$\circ $} at 670 28 \setdots <0.07cm> \plot 702 94 736 60 / \setsolid \put {$\circ $} at 744 52 \setdots <0.07cm>
\plot 751 60 760 69 / \setsolid \put {$\circ $} at 768 77 \setdots <0.07cm> \plot 736 45 727 36 / \setsolid \put {$\circ
$} at 719 28 \setdots <0.07cm> \plot 751 45 760 36 / \setsolid \put {$\circ $} at 768 28 \arrow <0.15cm> [0.25,0.5] from
-8 8 to -239 239 \plot -239 239 -390 390 / \put {$\bullet $} at -398 398 \arrow <0.15cm> [0.25,0.5] from -390 406 to
-279 517 \plot -279 517 -208 588 / \put {$\bullet $} at -200 596 \arrow <0.15cm> [0.25,0.5] from -192 604 to -141 655
\plot -141 655 -109 687 / \put {$\bullet $} at -101 695 \setdots <0.07cm> \plot -94 702 -60 736 / \setsolid \put {$\circ
$} at -52 744 \setdots <0.07cm> \plot -45 751 -36 760 / \setsolid \put {$\circ $} at -28 768 \setdots <0.07cm> \plot -60
751 -69 760 / \setsolid \put {$\circ $} at -77 768 \setdots <0.07cm> \plot -45 736 -36 727 / \setsolid \put {$\circ $}
at -28 719 \arrow <0.10cm> [0.25,0.5] from -109 702 to -131 724 \plot -131 724 -143 736 / \put {$\bullet $} at -150 744
\setdots <0.07cm> \plot -143 751 -134 760 / \setsolid \put {$\circ $} at -126 768 \arrow <0.05cm> [0.25,0.5] from -158
751 to -165 758 \plot -165 758 -167 760 / \put {$\bullet $} at -175 768 \arrow <0.05cm> [0.25,0.5] from -167 727 to -160
734 \plot -160 734 -158 736 / \put {$\bullet $} at -175 719 \arrow <0.10cm> [0.25,0.5] from -60 653 to -82 675 \plot -82
675 -94 687 / \put {$\bullet $} at -52 646 \setdots <0.07cm> \plot -45 653 -36 662 / \setsolid \put {$\circ $} at -28
670 \arrow <0.05cm> [0.25,0.5] from -69 629 to -62 636 \plot -62 636 -60 638 / \put {$\bullet $} at -77 621 \arrow
<0.05cm> [0.25,0.5] from -36 629 to -43 636 \plot -43 636 -45 638 / \put {$\bullet $} at -28 621 \setdots <0.07cm> \plot
-208 604 -291 687 / \setsolid \put {$\circ $} at -299 695 \setdots <0.07cm> \plot -291 702 -257 736 / \setsolid \put
{$\circ $} at -249 744 \setdots <0.07cm> \plot -242 751 -233 760 / \setsolid \put {$\circ $} at -225 768 \setdots
<0.07cm> \plot -257 751 -266 760 / \setsolid \put {$\circ $} at -274 768 \setdots <0.07cm> \plot -242 736 -233 727 /
\setsolid \put {$\circ $} at -225 719 \setdots <0.07cm> \plot -306 702 -340 736 / \setsolid \put {$\circ $} at -348 744
\setdots <0.07cm> \plot -340 751 -331 760 / \setsolid \put {$\circ $} at -323 768 \setdots <0.07cm> \plot -355 751 -364
760 / \setsolid \put {$\circ $} at -372 768 \setdots <0.07cm> \plot -355 736 -364 727 / \setsolid \put {$\circ $} at
-372 719 \setdots <0.07cm> \plot -306 687 -340 653 / \setsolid \put {$\circ $} at -348 646 \setdots <0.07cm> \plot -355
653 -364 662 / \setsolid \put {$\circ $} at -372 670 \setdots <0.07cm> \plot -355 638 -364 629 / \setsolid \put {$\circ
$} at -372 621 \setdots <0.07cm> \plot -340 638 -331 629 / \setsolid \put {$\circ $} at -323 621 \arrow <0.15cm>
[0.25,0.5] from -109 505 to -161 557 \plot -161 557 -192 588 / \put {$\bullet $} at -101 498 \setdots <0.07cm> \plot -94
505 -60 539 / \setsolid \put {$\circ $} at -52 547 \setdots <0.07cm> \plot -45 554 -36 563 / \setsolid \put {$\circ $}
at -28 571 \setdots <0.07cm> \plot -60 554 -69 563 / \setsolid \put {$\circ $} at -77 571 \setdots <0.07cm> \plot -45
539 -36 530 / \setsolid \put {$\circ $} at -28 522 \arrow <0.10cm> [0.25,0.5] from -143 456 to -121 478 \plot -121 478
-109 490 / \put {$\bullet $} at -150 448 \setdots <0.07cm> \plot -158 456 -167 465 / \setsolid \put {$\circ $} at -175
473 \arrow <0.05cm> [0.25,0.5] from -167 432 to -160 439 \plot -160 439 -158 441 / \put {$\bullet $} at -175 424 \arrow
<0.05cm> [0.25,0.5] from -134 432 to -141 439 \plot -141 439 -143 441 / \put {$\bullet $} at -126 424 \arrow <0.10cm>
[0.25,0.5] from -60 456 to -82 478 \plot -82 478 -94 490 / \put {$\bullet $} at -52 448 \setdots <0.07cm> \plot -45 456
-36 465 / \setsolid \put {$\circ $} at -28 473 \arrow <0.05cm> [0.25,0.5] from -69 432 to -62 439 \plot -62 439 -60 441
/ \put {$\bullet $} at -77 424 \arrow <0.05cm> [0.25,0.5] from -36 432 to -43 439 \plot -43 439 -45 441 / \put {$\bullet
$} at -28 424 \setdots <0.07cm> \plot -406 406 -588 588 / \setsolid \put {$\circ $} at -596 596 \setdots <0.07cm> \plot
-588 604 -505 687 / \setsolid \put {$\circ $} at -498 695 \setdots <0.07cm> \plot -490 702 -456 736 / \setsolid \put
{$\circ $} at -448 744 \setdots <0.07cm> \plot -441 751 -432 760 / \setsolid \put {$\circ $} at -424 768 \setdots
<0.07cm> \plot -456 751 -465 760 / \setsolid \put {$\circ $} at -473 768 \setdots <0.07cm> \plot -441 736 -432 727 /
\setsolid \put {$\circ $} at -424 719 \setdots <0.07cm> \plot -505 702 -539 736 / \setsolid \put {$\circ $} at -547 744
\setdots <0.07cm> \plot -539 751 -530 760 / \setsolid \put {$\circ $} at -522 768 \setdots <0.07cm> \plot -554 751 -563
760 / \setsolid \put {$\circ $} at -571 768 \setdots <0.07cm> \plot -554 736 -563 727 / \setsolid \put {$\circ $} at
-571 719 \setdots <0.07cm> \plot -490 687 -456 653 / \setsolid \put {$\circ $} at -448 646 \setdots <0.07cm> \plot -441
653 -432 662 / \setsolid \put {$\circ $} at -424 670 \setdots <0.07cm> \plot -456 638 -465 629 / \setsolid \put {$\circ
$} at -473 621 \setdots <0.07cm> \plot -441 638 -432 629 / \setsolid \put {$\circ $} at -424 621 \setdots <0.07cm> \plot
-604 604 -687 687 / \setsolid \put {$\circ $} at -695 695 \setdots <0.07cm> \plot -687 702 -653 736 / \setsolid \put
{$\circ $} at -646 744 \setdots <0.07cm> \plot -638 751 -629 760 / \setsolid \put {$\circ $} at -621 768 \setdots
<0.07cm> \plot -653 751 -662 760 / \setsolid \put {$\circ $} at -670 768 \setdots <0.07cm> \plot -638 736 -629 727 /
\setsolid \put {$\circ $} at -621 719 \setdots <0.07cm> \plot -702 702 -736 736 / \setsolid \put {$\circ $} at -744 744
\setdots <0.07cm> \plot -736 751 -727 760 / \setsolid \put {$\circ $} at -719 768 \setdots <0.07cm> \plot -751 751 -760
760 / \setsolid \put {$\circ $} at -768 768 \setdots <0.07cm> \plot -751 736 -760 727 / \setsolid \put {$\circ $} at
-768 719 \setdots <0.07cm> \plot -702 687 -736 653 / \setsolid \put {$\circ $} at -744 646 \setdots <0.07cm> \plot -751
653 -760 662 / \setsolid \put {$\circ $} at -768 670 \setdots <0.07cm> \plot -751 638 -760 629 / \setsolid \put {$\circ
$} at -768 621 \setdots <0.07cm> \plot -736 638 -727 629 / \setsolid \put {$\circ $} at -719 621 \setdots <0.07cm> \plot
-604 588 -687 505 / \setsolid \put {$\circ $} at -695 498 \setdots <0.07cm> \plot -702 505 -736 539 / \setsolid \put
{$\circ $} at -744 547 \setdots <0.07cm> \plot -736 554 -727 563 / \setsolid \put {$\circ $} at -719 571 \setdots
<0.07cm> \plot -751 554 -760 563 / \setsolid \put {$\circ $} at -768 571 \setdots <0.07cm> \plot -751 539 -760 530 /
\setsolid \put {$\circ $} at -768 522 \setdots <0.07cm> \plot -702 490 -736 456 / \setsolid \put {$\circ $} at -744 448
\setdots <0.07cm> \plot -751 456 -760 465 / \setsolid \put {$\circ $} at -768 473 \setdots <0.07cm> \plot -751 441 -760
432 / \setsolid \put {$\circ $} at -768 424 \setdots <0.07cm> \plot -736 441 -727 432 / \setsolid \put {$\circ $} at
-719 424 \setdots <0.07cm> \plot -687 490 -653 456 / \setsolid \put {$\circ $} at -646 448 \setdots <0.07cm> \plot -638
456 -629 465 / \setsolid \put {$\circ $} at -621 473 \setdots <0.07cm> \plot -653 441 -662 432 / \setsolid \put {$\circ
$} at -670 424 \setdots <0.07cm> \plot -638 441 -629 432 / \setsolid \put {$\circ $} at -621 424 \arrow <0.15cm>
[0.25,0.5] from -588 208 to -477 319 \plot -477 319 -406 390 / \put {$\bullet $} at -596 200 \setdots <0.07cm> \plot
-604 208 -687 291 / \setsolid \put {$\circ $} at -695 299 \setdots <0.07cm> \plot -687 306 -653 340 / \setsolid \put
{$\circ $} at -646 348 \setdots <0.07cm> \plot -638 355 -629 364 / \setsolid \put {$\circ $} at -621 372 \setdots
<0.07cm> \plot -653 355 -662 364 / \setsolid \put {$\circ $} at -670 372 \setdots <0.07cm> \plot -638 340 -629 331 /
\setsolid \put {$\circ $} at -621 323 \setdots <0.07cm> \plot -702 306 -736 340 / \setsolid \put {$\circ $} at -744 348
\setdots <0.07cm> \plot -736 355 -727 364 / \setsolid \put {$\circ $} at -719 372 \setdots <0.07cm> \plot -751 355 -760
364 / \setsolid \put {$\circ $} at -768 372 \setdots <0.07cm> \plot -751 340 -760 331 / \setsolid \put {$\circ $} at
-768 323 \setdots <0.07cm> \plot -702 291 -736 257 / \setsolid \put {$\circ $} at -744 249 \setdots <0.07cm> \plot -751
257 -760 266 / \setsolid \put {$\circ $} at -768 274 \setdots <0.07cm> \plot -751 242 -760 233 / \setsolid \put {$\circ
$} at -768 225 \setdots <0.07cm> \plot -736 242 -727 233 / \setsolid \put {$\circ $} at -719 225 \arrow <0.15cm>
[0.25,0.5] from -687 109 to -635 161 \plot -635 161 -604 192 / \put {$\bullet $} at -695 101 \setdots <0.07cm> \plot
-702 109 -736 143 / \setsolid \put {$\circ $} at -744 150 \setdots <0.07cm> \plot -736 158 -727 167 / \setsolid \put
{$\circ $} at -719 175 \setdots <0.07cm> \plot -751 158 -760 167 / \setsolid \put {$\circ $} at -768 175 \setdots
<0.07cm> \plot -751 143 -760 134 / \setsolid \put {$\circ $} at -768 126 \arrow <0.10cm> [0.25,0.5] from -736 60 to -714
82 \plot -714 82 -702 94 / \put {$\bullet $} at -744 52 \setdots <0.07cm> \plot -751 60 -760 69 / \setsolid \put {$\circ
$} at -768 77 \arrow <0.05cm> [0.25,0.5] from -760 36 to -753 43 \plot -753 43 -751 45 / \put {$\bullet $} at -768 28
\arrow <0.05cm> [0.25,0.5] from -727 36 to -734 43 \plot -734 43 -736 45 / \put {$\bullet $} at -719 28 \arrow <0.10cm>
[0.25,0.5] from -653 60 to -675 82 \plot -675 82 -687 94 / \put {$\bullet $} at -646 52 \setdots <0.07cm> \plot -638 60
-629 69 / \setsolid \put {$\circ $} at -621 77 \arrow <0.05cm> [0.25,0.5] from -662 36 to -655 43 \plot -655 43 -653 45
/ \put {$\bullet $} at -670 28 \arrow <0.05cm> [0.25,0.5] from -629 36 to -636 43 \plot -636 43 -638 45 / \put {$\bullet
$} at -621 28 \arrow <0.15cm> [0.25,0.5] from -505 109 to -557 161 \plot -557 161 -588 192 / \put {$\bullet $} at -498
101 \setdots <0.07cm> \plot -490 109 -456 143 / \setsolid \put {$\circ $} at -448 150 \setdots <0.07cm> \plot -441 158
-432 167 / \setsolid \put {$\circ $} at -424 175 \setdots <0.07cm> \plot -456 158 -465 167 / \setsolid \put {$\circ $}
at -473 175 \setdots <0.07cm> \plot -441 143 -432 134 / \setsolid \put {$\circ $} at -424 126 \arrow <0.10cm> [0.25,0.5]
from -539 60 to -517 82 \plot -517 82 -505 94 / \put {$\bullet $} at -547 52 \setdots <0.07cm> \plot -554 60 -563 69 /
\setsolid \put {$\circ $} at -571 77 \arrow <0.05cm> [0.25,0.5] from -563 36 to -556 43 \plot -556 43 -554 45 / \put
{$\bullet $} at -571 28 \arrow <0.05cm> [0.25,0.5] from -530 36 to -537 43 \plot -537 43 -539 45 / \put {$\bullet $} at
-522 28 \arrow <0.10cm> [0.25,0.5] from -456 60 to -478 82 \plot -478 82 -490 94 / \put {$\bullet $} at -448 52 \setdots
<0.07cm> \plot -441 60 -432 69 / \setsolid \put {$\circ $} at -424 77 \arrow <0.05cm> [0.25,0.5] from -465 36 to -458 43
\plot -458 43 -456 45 / \put {$\bullet $} at -473 28 \arrow <0.05cm> [0.25,0.5] from -432 36 to -439 43 \plot -439 43
-441 45 / \put {$\bullet $} at -424 28 \arrow <0.15cm> [0.25,0.5] from -390 -390 to -159 -159 \plot -159 -159 -8 -8 /
\put {$\bullet $} at -398 -398 \setdots <0.07cm> \plot -406 -390 -588 -208 / \setsolid \put {$\circ $} at -596 -200
\setdots <0.07cm> \plot -588 -192 -505 -109 / \setsolid \put {$\circ $} at -498 -101 \setdots <0.07cm> \plot -490 -94
-456 -60 / \setsolid \put {$\circ $} at -448 -52 \setdots <0.07cm> \plot -441 -45 -432 -36 / \setsolid \put {$\circ $}
at -424 -28 \setdots <0.07cm> \plot -456 -45 -465 -36 / \setsolid \put {$\circ $} at -473 -28 \setdots <0.07cm> \plot
-441 -60 -432 -69 / \setsolid \put {$\circ $} at -424 -77 \setdots <0.07cm> \plot -505 -94 -539 -60 / \setsolid \put
{$\circ $} at -547 -52 \setdots <0.07cm> \plot -539 -45 -530 -36 / \setsolid \put {$\circ $} at -522 -28 \setdots
<0.07cm> \plot -554 -45 -563 -36 / \setsolid \put {$\circ $} at -571 -28 \setdots <0.07cm> \plot -554 -60 -563 -69 /
\setsolid \put {$\circ $} at -571 -77 \setdots <0.07cm> \plot -490 -109 -456 -143 / \setsolid \put {$\circ $} at -448
-150 \setdots <0.07cm> \plot -441 -143 -432 -134 / \setsolid \put {$\circ $} at -424 -126 \setdots <0.07cm> \plot -456
-158 -465 -167 / \setsolid \put {$\circ $} at -473 -175 \setdots <0.07cm> \plot -441 -158 -432 -167 / \setsolid \put
{$\circ $} at -424 -175 \setdots <0.07cm> \plot -604 -192 -687 -109 / \setsolid \put {$\circ $} at -695 -101 \setdots
<0.07cm> \plot -687 -94 -653 -60 / \setsolid \put {$\circ $} at -646 -52 \setdots <0.07cm> \plot -638 -45 -629 -36 /
\setsolid \put {$\circ $} at -621 -28 \setdots <0.07cm> \plot -653 -45 -662 -36 / \setsolid \put {$\circ $} at -670 -28
\setdots <0.07cm> \plot -638 -60 -629 -69 / \setsolid \put {$\circ $} at -621 -77 \setdots <0.07cm> \plot -702 -94 -736
-60 / \setsolid \put {$\circ $} at -744 -52 \setdots <0.07cm> \plot -736 -45 -727 -36 / \setsolid \put {$\circ $} at
-719 -28 \setdots <0.07cm> \plot -751 -45 -760 -36 / \setsolid \put {$\circ $} at -768 -28 \setdots <0.07cm> \plot -751
-60 -760 -69 / \setsolid \put {$\circ $} at -768 -77 \setdots <0.07cm> \plot -702 -109 -736 -143 / \setsolid \put
{$\circ $} at -744 -150 \setdots <0.07cm> \plot -751 -143 -760 -134 / \setsolid \put {$\circ $} at -768 -126 \setdots
<0.07cm> \plot -751 -158 -760 -167 / \setsolid \put {$\circ $} at -768 -175 \setdots <0.07cm> \plot -736 -158 -727 -167
/ \setsolid \put {$\circ $} at -719 -175 \setdots <0.07cm> \plot -604 -208 -687 -291 / \setsolid \put {$\circ $} at -695
-299 \setdots <0.07cm> \plot -702 -291 -736 -257 / \setsolid \put {$\circ $} at -744 -249 \setdots <0.07cm> \plot -736
-242 -727 -233 / \setsolid \put {$\circ $} at -719 -225 \setdots <0.07cm> \plot -751 -242 -760 -233 / \setsolid \put
{$\circ $} at -768 -225 \setdots <0.07cm> \plot -751 -257 -760 -266 / \setsolid \put {$\circ $} at -768 -274 \setdots
<0.07cm> \plot -702 -306 -736 -340 / \setsolid \put {$\circ $} at -744 -348 \setdots <0.07cm> \plot -751 -340 -760 -331
/ \setsolid \put {$\circ $} at -768 -323 \setdots <0.07cm> \plot -751 -355 -760 -364 / \setsolid \put {$\circ $} at -768
-372 \setdots <0.07cm> \plot -736 -355 -727 -364 / \setsolid \put {$\circ $} at -719 -372 \setdots <0.07cm> \plot -687
-306 -653 -340 / \setsolid \put {$\circ $} at -646 -348 \setdots <0.07cm> \plot -638 -340 -629 -331 / \setsolid \put
{$\circ $} at -621 -323 \setdots <0.07cm> \plot -653 -355 -662 -364 / \setsolid \put {$\circ $} at -670 -372 \setdots
<0.07cm> \plot -638 -355 -629 -364 / \setsolid \put {$\circ $} at -621 -372 \arrow <0.15cm> [0.25,0.5] from -588 -588 to
-477 -477 \plot -477 -477 -406 -406 / \put {$\bullet $} at -596 -596 \setdots <0.07cm> \plot -604 -588 -687 -505 /
\setsolid \put {$\circ $} at -695 -498 \setdots <0.07cm> \plot -687 -490 -653 -456 / \setsolid \put {$\circ $} at -646
-448 \setdots <0.07cm> \plot -638 -441 -629 -432 / \setsolid \put {$\circ $} at -621 -424 \setdots <0.07cm> \plot -653
-441 -662 -432 / \setsolid \put {$\circ $} at -670 -424 \setdots <0.07cm> \plot -638 -456 -629 -465 / \setsolid \put
{$\circ $} at -621 -473 \setdots <0.07cm> \plot -702 -490 -736 -456 / \setsolid \put {$\circ $} at -744 -448 \setdots
<0.07cm> \plot -736 -441 -727 -432 / \setsolid \put {$\circ $} at -719 -424 \setdots <0.07cm> \plot -751 -441 -760 -432
/ \setsolid \put {$\circ $} at -768 -424 \setdots <0.07cm> \plot -751 -456 -760 -465 / \setsolid \put {$\circ $} at -768
-473 \setdots <0.07cm> \plot -702 -505 -736 -539 / \setsolid \put {$\circ $} at -744 -547 \setdots <0.07cm> \plot -751
-539 -760 -530 / \setsolid \put {$\circ $} at -768 -522 \setdots <0.07cm> \plot -751 -554 -760 -563 / \setsolid \put
{$\circ $} at -768 -571 \setdots <0.07cm> \plot -736 -554 -727 -563 / \setsolid \put {$\circ $} at -719 -571 \arrow
<0.15cm> [0.25,0.5] from -687 -687 to -635 -635 \plot -635 -635 -604 -604 / \put {$\bullet $} at -695 -695 \setdots
<0.07cm> \plot -702 -687 -736 -653 / \setsolid \put {$\circ $} at -744 -646 \setdots <0.07cm> \plot -736 -638 -727 -629
/ \setsolid \put {$\circ $} at -719 -621 \setdots <0.07cm> \plot -751 -638 -760 -629 / \setsolid \put {$\circ $} at -768
-621 \setdots <0.07cm> \plot -751 -653 -760 -662 / \setsolid \put {$\circ $} at -768 -670 \arrow <0.10cm> [0.25,0.5]
from -736 -736 to -714 -714 \plot -714 -714 -702 -702 / \put {$\bullet $} at -744 -744 \setdots <0.07cm> \plot -751 -736
-760 -727 / \setsolid \put {$\circ $} at -768 -719 \arrow <0.05cm> [0.25,0.5] from -760 -760 to -753 -753 \plot -753
-753 -751 -751 / \put {$\bullet $} at -768 -768 \arrow <0.05cm> [0.25,0.5] from -727 -760 to -734 -753 \plot -734 -753
-736 -751 / \put {$\bullet $} at -719 -768 \arrow <0.10cm> [0.25,0.5] from -653 -736 to -675 -714 \plot -675 -714 -687
-702 / \put {$\bullet $} at -646 -744 \setdots <0.07cm> \plot -638 -736 -629 -727 / \setsolid \put {$\circ $} at -621
-719 \arrow <0.05cm> [0.25,0.5] from -662 -760 to -655 -753 \plot -655 -753 -653 -751 / \put {$\bullet $} at -670 -768
\arrow <0.05cm> [0.25,0.5] from -629 -760 to -636 -753 \plot -636 -753 -638 -751 / \put {$\bullet $} at -621 -768 \arrow
<0.15cm> [0.25,0.5] from -505 -687 to -557 -635 \plot -557 -635 -588 -604 / \put {$\bullet $} at -498 -695 \setdots
<0.07cm> \plot -490 -687 -456 -653 / \setsolid \put {$\circ $} at -448 -646 \setdots <0.07cm> \plot -441 -638 -432 -629
/ \setsolid \put {$\circ $} at -424 -621 \setdots <0.07cm> \plot -456 -638 -465 -629 / \setsolid \put {$\circ $} at -473
-621 \setdots <0.07cm> \plot -441 -653 -432 -662 / \setsolid \put {$\circ $} at -424 -670 \arrow <0.10cm> [0.25,0.5]
from -539 -736 to -517 -714 \plot -517 -714 -505 -702 / \put {$\bullet $} at -547 -744 \setdots <0.07cm> \plot -554 -736
-563 -727 / \setsolid \put {$\circ $} at -571 -719 \arrow <0.05cm> [0.25,0.5] from -563 -760 to -556 -753 \plot -556
-753 -554 -751 / \put {$\bullet $} at -571 -768 \arrow <0.05cm> [0.25,0.5] from -530 -760 to -537 -753 \plot -537 -753
-539 -751 / \put {$\bullet $} at -522 -768 \arrow <0.10cm> [0.25,0.5] from -456 -736 to -478 -714 \plot -478 -714 -490
-702 / \put {$\bullet $} at -448 -744 \setdots <0.07cm> \plot -441 -736 -432 -727 / \setsolid \put {$\circ $} at -424
-719 \arrow <0.05cm> [0.25,0.5] from -465 -760 to -458 -753 \plot -458 -753 -456 -751 / \put {$\bullet $} at -473 -768
\arrow <0.05cm> [0.25,0.5] from -432 -760 to -439 -753 \plot -439 -753 -441 -751 / \put {$\bullet $} at -424 -768 \arrow
<0.15cm> [0.25,0.5] from -208 -588 to -319 -477 \plot -319 -477 -390 -406 / \put {$\bullet $} at -200 -596 \setdots
<0.07cm> \plot -192 -588 -109 -505 / \setsolid \put {$\circ $} at -101 -498 \setdots <0.07cm> \plot -94 -490 -60 -456 /
\setsolid \put {$\circ $} at -52 -448 \setdots <0.07cm> \plot -45 -441 -36 -432 / \setsolid \put {$\circ $} at -28 -424
\setdots <0.07cm> \plot -60 -441 -69 -432 / \setsolid \put {$\circ $} at -77 -424 \setdots <0.07cm> \plot -45 -456 -36
-465 / \setsolid \put {$\circ $} at -28 -473 \setdots <0.07cm> \plot -109 -490 -143 -456 / \setsolid \put {$\circ $} at
-150 -448 \setdots <0.07cm> \plot -143 -441 -134 -432 / \setsolid \put {$\circ $} at -126 -424 \setdots <0.07cm> \plot
-158 -441 -167 -432 / \setsolid \put {$\circ $} at -175 -424 \setdots <0.07cm> \plot -158 -456 -167 -465 / \setsolid
\put {$\circ $} at -175 -473 \setdots <0.07cm> \plot -94 -505 -60 -539 / \setsolid \put {$\circ $} at -52 -547 \setdots
<0.07cm> \plot -45 -539 -36 -530 / \setsolid \put {$\circ $} at -28 -522 \setdots <0.07cm> \plot -60 -554 -69 -563 /
\setsolid \put {$\circ $} at -77 -571 \setdots <0.07cm> \plot -45 -554 -36 -563 / \setsolid \put {$\circ $} at -28 -571
\arrow <0.15cm> [0.25,0.5] from -291 -687 to -239 -635 \plot -239 -635 -208 -604 / \put {$\bullet $} at -299 -695
\setdots <0.07cm> \plot -306 -687 -340 -653 / \setsolid \put {$\circ $} at -348 -646 \setdots <0.07cm> \plot -340 -638
-331 -629 / \setsolid \put {$\circ $} at -323 -621 \setdots <0.07cm> \plot -355 -638 -364 -629 / \setsolid \put {$\circ
$} at -372 -621 \setdots <0.07cm> \plot -355 -653 -364 -662 / \setsolid \put {$\circ $} at -372 -670 \arrow <0.10cm>
[0.25,0.5] from -340 -736 to -318 -714 \plot -318 -714 -306 -702 / \put {$\bullet $} at -348 -744 \setdots <0.07cm>
\plot -355 -736 -364 -727 / \setsolid \put {$\circ $} at -372 -719 \arrow <0.05cm> [0.25,0.5] from -364 -760 to -357
-753 \plot -357 -753 -355 -751 / \put {$\bullet $} at -372 -768 \arrow <0.05cm> [0.25,0.5] from -331 -760 to -338 -753
\plot -338 -753 -340 -751 / \put {$\bullet $} at -323 -768 \arrow <0.10cm> [0.25,0.5] from -257 -736 to -279 -714 \plot
-279 -714 -291 -702 / \put {$\bullet $} at -249 -744 \setdots <0.07cm> \plot -242 -736 -233 -727 / \setsolid \put
{$\circ $} at -225 -719 \arrow <0.05cm> [0.25,0.5] from -266 -760 to -259 -753 \plot -259 -753 -257 -751 / \put
{$\bullet $} at -274 -768 \arrow <0.05cm> [0.25,0.5] from -233 -760 to -240 -753 \plot -240 -753 -242 -751 / \put
{$\bullet $} at -225 -768 \arrow <0.15cm> [0.25,0.5] from -109 -687 to -161 -635 \plot -161 -635 -192 -604 / \put
{$\bullet $} at -101 -695 \setdots <0.07cm> \plot -94 -687 -60 -653 / \setsolid \put {$\circ $} at -52 -646 \setdots
<0.07cm> \plot -45 -638 -36 -629 / \setsolid \put {$\circ $} at -28 -621 \setdots <0.07cm> \plot -60 -638 -69 -629 /
\setsolid \put {$\circ $} at -77 -621 \setdots <0.07cm> \plot -45 -653 -36 -662 / \setsolid \put {$\circ $} at -28 -670
\arrow <0.10cm> [0.25,0.5] from -143 -736 to -121 -714 \plot -121 -714 -109 -702 / \put {$\bullet $} at -150 -744
\setdots <0.07cm> \plot -158 -736 -167 -727 / \setsolid \put {$\circ $} at -175 -719 \arrow <0.05cm> [0.25,0.5] from
-167 -760 to -160 -753 \plot -160 -753 -158 -751 / \put {$\bullet $} at -175 -768 \arrow <0.05cm> [0.25,0.5] from -134
-760 to -141 -753 \plot -141 -753 -143 -751 / \put {$\bullet $} at -126 -768 \arrow <0.10cm> [0.25,0.5] from -60 -736 to
-82 -714 \plot -82 -714 -94 -702 / \put {$\bullet $} at -52 -744 \setdots <0.07cm> \plot -45 -736 -36 -727 / \setsolid
\put {$\circ $} at -28 -719 \arrow <0.05cm> [0.25,0.5] from -69 -760 to -62 -753 \plot -62 -753 -60 -751 / \put
{$\bullet $} at -77 -768 \arrow <0.05cm> [0.25,0.5] from -36 -760 to -43 -753 \plot -43 -753 -45 -751 / \put {$\bullet
$} at -28 -768 \arrow <0.15cm> [0.25,0.5] from 390 -390 to 159 -159 \plot 159 -159 8 -8 / \put {$\bullet $} at 398 -398
\setdots <0.07cm> \plot 406 -390 588 -208 / \setsolid \put {$\circ $} at 596 -200 \setdots <0.07cm> \plot 604 -192 687
-109 / \setsolid \put {$\circ $} at 695 -101 \setdots <0.07cm> \plot 702 -94 736 -60 / \setsolid \put {$\circ $} at 744
-52 \setdots <0.07cm> \plot 751 -45 760 -36 / \setsolid \put {$\circ $} at 768 -28 \setdots <0.07cm> \plot 736 -45 727
-36 / \setsolid \put {$\circ $} at 719 -28 \setdots <0.07cm> \plot 751 -60 760 -69 / \setsolid \put {$\circ $} at 768
-77 \setdots <0.07cm> \plot 687 -94 653 -60 / \setsolid \put {$\circ $} at 646 -52 \setdots <0.07cm> \plot 653 -45 662
-36 / \setsolid \put {$\circ $} at 670 -28 \setdots <0.07cm> \plot 638 -45 629 -36 / \setsolid \put {$\circ $} at 621
-28 \setdots <0.07cm> \plot 638 -60 629 -69 / \setsolid \put {$\circ $} at 621 -77 \setdots <0.07cm> \plot 702 -109 736
-143 / \setsolid \put {$\circ $} at 744 -150 \setdots <0.07cm> \plot 751 -143 760 -134 / \setsolid \put {$\circ $} at
768 -126 \setdots <0.07cm> \plot 736 -158 727 -167 / \setsolid \put {$\circ $} at 719 -175 \setdots <0.07cm> \plot 751
-158 760 -167 / \setsolid \put {$\circ $} at 768 -175 \setdots <0.07cm> \plot 588 -192 505 -109 / \setsolid \put {$\circ
$} at 498 -101 \setdots <0.07cm> \plot 505 -94 539 -60 / \setsolid \put {$\circ $} at 547 -52 \setdots <0.07cm> \plot
554 -45 563 -36 / \setsolid \put {$\circ $} at 571 -28 \setdots <0.07cm> \plot 539 -45 530 -36 / \setsolid \put {$\circ
$} at 522 -28 \setdots <0.07cm> \plot 554 -60 563 -69 / \setsolid \put {$\circ $} at 571 -77 \setdots <0.07cm> \plot 490
-94 456 -60 / \setsolid \put {$\circ $} at 448 -52 \setdots <0.07cm> \plot 456 -45 465 -36 / \setsolid \put {$\circ $}
at 473 -28 \setdots <0.07cm> \plot 441 -45 432 -36 / \setsolid \put {$\circ $} at 424 -28 \setdots <0.07cm> \plot 441
-60 432 -69 / \setsolid \put {$\circ $} at 424 -77 \setdots <0.07cm> \plot 490 -109 456 -143 / \setsolid \put {$\circ $}
at 448 -150 \setdots <0.07cm> \plot 441 -143 432 -134 / \setsolid \put {$\circ $} at 424 -126 \setdots <0.07cm> \plot
441 -158 432 -167 / \setsolid \put {$\circ $} at 424 -175 \setdots <0.07cm> \plot 456 -158 465 -167 / \setsolid \put
{$\circ $} at 473 -175 \setdots <0.07cm> \plot 604 -208 687 -291 / \setsolid \put {$\circ $} at 695 -299 \setdots
<0.07cm> \plot 702 -291 736 -257 / \setsolid \put {$\circ $} at 744 -249 \setdots <0.07cm> \plot 751 -242 760 -233 /
\setsolid \put {$\circ $} at 768 -225 \setdots <0.07cm> \plot 736 -242 727 -233 / \setsolid \put {$\circ $} at 719 -225
\setdots <0.07cm> \plot 751 -257 760 -266 / \setsolid \put {$\circ $} at 768 -274 \setdots <0.07cm> \plot 687 -306 653
-340 / \setsolid \put {$\circ $} at 646 -348 \setdots <0.07cm> \plot 638 -340 629 -331 / \setsolid \put {$\circ $} at
621 -323 \setdots <0.07cm> \plot 638 -355 629 -364 / \setsolid \put {$\circ $} at 621 -372 \setdots <0.07cm> \plot 653
-355 662 -364 / \setsolid \put {$\circ $} at 670 -372 \setdots <0.07cm> \plot 702 -306 736 -340 / \setsolid \put {$\circ
$} at 744 -348 \setdots <0.07cm> \plot 751 -340 760 -331 / \setsolid \put {$\circ $} at 768 -323 \setdots <0.07cm> \plot
736 -355 727 -364 / \setsolid \put {$\circ $} at 719 -372 \setdots <0.07cm> \plot 751 -355 760 -364 / \setsolid \put
{$\circ $} at 768 -372 \arrow <0.15cm> [0.25,0.5] from 208 -588 to 319 -477 \plot 319 -477 390 -406 / \put {$\bullet $}
at 200 -596 \setdots <0.07cm> \plot 192 -588 109 -505 / \setsolid \put {$\circ $} at 101 -498 \setdots <0.07cm> \plot
109 -490 143 -456 / \setsolid \put {$\circ $} at 150 -448 \setdots <0.07cm> \plot 158 -441 167 -432 / \setsolid \put
{$\circ $} at 175 -424 \setdots <0.07cm> \plot 143 -441 134 -432 / \setsolid \put {$\circ $} at 126 -424 \setdots
<0.07cm> \plot 158 -456 167 -465 / \setsolid \put {$\circ $} at 175 -473 \setdots <0.07cm> \plot 94 -490 60 -456 /
\setsolid \put {$\circ $} at 52 -448 \setdots <0.07cm> \plot 60 -441 69 -432 / \setsolid \put {$\circ $} at 77 -424
\setdots <0.07cm> \plot 45 -441 36 -432 / \setsolid \put {$\circ $} at 28 -424 \setdots <0.07cm> \plot 45 -456 36 -465 /
\setsolid \put {$\circ $} at 28 -473 \setdots <0.07cm> \plot 94 -505 60 -539 / \setsolid \put {$\circ $} at 52 -547
\setdots <0.07cm> \plot 45 -539 36 -530 / \setsolid \put {$\circ $} at 28 -522 \setdots <0.07cm> \plot 45 -554 36 -563 /
\setsolid \put {$\circ $} at 28 -571 \setdots <0.07cm> \plot 60 -554 69 -563 / \setsolid \put {$\circ $} at 77 -571
\arrow <0.15cm> [0.25,0.5] from 109 -687 to 161 -635 \plot 161 -635 192 -604 / \put {$\bullet $} at 101 -695 \setdots
<0.07cm> \plot 94 -687 60 -653 / \setsolid \put {$\circ $} at 52 -646 \setdots <0.07cm> \plot 60 -638 69 -629 /
\setsolid \put {$\circ $} at 77 -621 \setdots <0.07cm> \plot 45 -638 36 -629 / \setsolid \put {$\circ $} at 28 -621
\setdots <0.07cm> \plot 45 -653 36 -662 / \setsolid \put {$\circ $} at 28 -670 \arrow <0.10cm> [0.25,0.5] from 60 -736
to 82 -714 \plot 82 -714 94 -702 / \put {$\bullet $} at 52 -744 \setdots <0.07cm> \plot 45 -736 36 -727 / \setsolid \put
{$\circ $} at 28 -719 \arrow <0.05cm> [0.25,0.5] from 36 -760 to 43 -753 \plot 43 -753 45 -751 / \put {$\bullet $} at 28
-768 \arrow <0.05cm> [0.25,0.5] from 69 -760 to 62 -753 \plot 62 -753 60 -751 / \put {$\bullet $} at 77 -768 \arrow
<0.10cm> [0.25,0.5] from 143 -736 to 121 -714 \plot 121 -714 109 -702 / \put {$\bullet $} at 150 -744 \setdots <0.07cm>
\plot 158 -736 167 -727 / \setsolid \put {$\circ $} at 175 -719 \arrow <0.05cm> [0.25,0.5] from 134 -760 to 141 -753
\plot 141 -753 143 -751 / \put {$\bullet $} at 126 -768 \arrow <0.05cm> [0.25,0.5] from 167 -760 to 160 -753 \plot 160
-753 158 -751 / \put {$\bullet $} at 175 -768 \arrow <0.15cm> [0.25,0.5] from 291 -687 to 239 -635 \plot 239 -635 208
-604 / \put {$\bullet $} at 299 -695 \setdots <0.07cm> \plot 306 -687 340 -653 / \setsolid \put {$\circ $} at 348 -646
\setdots <0.07cm> \plot 355 -638 364 -629 / \setsolid \put {$\circ $} at 372 -621 \setdots <0.07cm> \plot 340 -638 331
-629 / \setsolid \put {$\circ $} at 323 -621 \setdots <0.07cm> \plot 355 -653 364 -662 / \setsolid \put {$\circ $} at
372 -670 \arrow <0.10cm> [0.25,0.5] from 257 -736 to 279 -714 \plot 279 -714 291 -702 / \put {$\bullet $} at 249 -744
\setdots <0.07cm> \plot 242 -736 233 -727 / \setsolid \put {$\circ $} at 225 -719 \arrow <0.05cm> [0.25,0.5] from 233
-760 to 240 -753 \plot 240 -753 242 -751 / \put {$\bullet $} at 225 -768 \arrow <0.05cm> [0.25,0.5] from 266 -760 to 259
-753 \plot 259 -753 257 -751 / \put {$\bullet $} at 274 -768 \arrow <0.10cm> [0.25,0.5] from 340 -736 to 318 -714 \plot
318 -714 306 -702 / \put {$\bullet $} at 348 -744 \setdots <0.07cm> \plot 355 -736 364 -727 / \setsolid \put {$\circ $}
at 372 -719 \arrow <0.05cm> [0.25,0.5] from 331 -760 to 338 -753 \plot 338 -753 340 -751 / \put {$\bullet $} at 323 -768
\arrow <0.05cm> [0.25,0.5] from 364 -760 to 357 -753 \plot 357 -753 355 -751 / \put {$\bullet $} at 372 -768 \arrow
<0.15cm> [0.25,0.5] from 588 -588 to 477 -477 \plot 477 -477 406 -406 / \put {$\bullet $} at 596 -596 \setdots <0.07cm>
\plot 604 -588 687 -505 / \setsolid \put {$\circ $} at 695 -498 \setdots <0.07cm> \plot 702 -490 736 -456 / \setsolid
\put {$\circ $} at 744 -448 \setdots <0.07cm> \plot 751 -441 760 -432 / \setsolid \put {$\circ $} at 768 -424 \setdots
<0.07cm> \plot 736 -441 727 -432 / \setsolid \put {$\circ $} at 719 -424 \setdots <0.07cm> \plot 751 -456 760 -465 /
\setsolid \put {$\circ $} at 768 -473 \setdots <0.07cm> \plot 687 -490 653 -456 / \setsolid \put {$\circ $} at 646 -448
\setdots <0.07cm> \plot 653 -441 662 -432 / \setsolid \put {$\circ $} at 670 -424 \setdots <0.07cm> \plot 638 -441 629
-432 / \setsolid \put {$\circ $} at 621 -424 \setdots <0.07cm> \plot 638 -456 629 -465 / \setsolid \put {$\circ $} at
621 -473 \setdots <0.07cm> \plot 702 -505 736 -539 / \setsolid \put {$\circ $} at 744 -547 \setdots <0.07cm> \plot 751
-539 760 -530 / \setsolid \put {$\circ $} at 768 -522 \setdots <0.07cm> \plot 736 -554 727 -563 / \setsolid \put {$\circ
$} at 719 -571 \setdots <0.07cm> \plot 751 -554 760 -563 / \setsolid \put {$\circ $} at 768 -571 \arrow <0.15cm>
[0.25,0.5] from 505 -687 to 557 -635 \plot 557 -635 588 -604 / \put {$\bullet $} at 498 -695 \setdots <0.07cm> \plot 490
-687 456 -653 / \setsolid \put {$\circ $} at 448 -646 \setdots <0.07cm> \plot 456 -638 465 -629 / \setsolid \put {$\circ
$} at 473 -621 \setdots <0.07cm> \plot 441 -638 432 -629 / \setsolid \put {$\circ $} at 424 -621 \setdots <0.07cm> \plot
441 -653 432 -662 / \setsolid \put {$\circ $} at 424 -670 \arrow <0.10cm> [0.25,0.5] from 456 -736 to 478 -714 \plot 478
-714 490 -702 / \put {$\bullet $} at 448 -744 \setdots <0.07cm> \plot 441 -736 432 -727 / \setsolid \put {$\circ $} at
424 -719 \arrow <0.05cm> [0.25,0.5] from 432 -760 to 439 -753 \plot 439 -753 441 -751 / \put {$\bullet $} at 424 -768
\arrow <0.05cm> [0.25,0.5] from 465 -760 to 458 -753 \plot 458 -753 456 -751 / \put {$\bullet $} at 473 -768 \arrow
<0.10cm> [0.25,0.5] from 539 -736 to 517 -714 \plot 517 -714 505 -702 / \put {$\bullet $} at 547 -744 \setdots <0.07cm>
\plot 554 -736 563 -727 / \setsolid \put {$\circ $} at 571 -719 \arrow <0.05cm> [0.25,0.5] from 530 -760 to 537 -753
\plot 537 -753 539 -751 / \put {$\bullet $} at 522 -768 \arrow <0.05cm> [0.25,0.5] from 563 -760 to 556 -753 \plot 556
-753 554 -751 / \put {$\bullet $} at 571 -768 \arrow <0.15cm> [0.25,0.5] from 687 -687 to 635 -635 \plot 635 -635 604
-604 / \put {$\bullet $} at 695 -695 \setdots <0.07cm> \plot 702 -687 736 -653 / \setsolid \put {$\circ $} at 744 -646
\setdots <0.07cm> \plot 751 -638 760 -629 / \setsolid \put {$\circ $} at 768 -621 \setdots <0.07cm> \plot 736 -638 727
-629 / \setsolid \put {$\circ $} at 719 -621 \setdots <0.07cm> \plot 751 -653 760 -662 / \setsolid \put {$\circ $} at
768 -670 \arrow <0.10cm> [0.25,0.5] from 653 -736 to 675 -714 \plot 675 -714 687 -702 / \put {$\bullet $} at 646 -744
\setdots <0.07cm> \plot 638 -736 629 -727 / \setsolid \put {$\circ $} at 621 -719 \arrow <0.05cm> [0.25,0.5] from 629
-760 to 636 -753 \plot 636 -753 638 -751 / \put {$\bullet $} at 621 -768 \arrow <0.05cm> [0.25,0.5] from 662 -760 to 655
-753 \plot 655 -753 653 -751 / \put {$\bullet $} at 670 -768 \arrow <0.10cm> [0.25,0.5] from 736 -736 to 714 -714 \plot
714 -714 702 -702 / \put {$\bullet $} at 744 -744 \setdots <0.07cm> \plot 751 -736 760 -727 / \setsolid \put {$\circ $}
at 768 -719 \arrow <0.05cm> [0.25,0.5] from 727 -760 to 734 -753 \plot 734 -753 736 -751 / \put {$\bullet $} at 719 -768
\arrow <0.05cm> [0.25,0.5] from 760 -760 to 753 -753 \plot 753 -753 751 -751 / \put {$\bullet $} at 768 -768 \endpicture
\hfill \null

\begincenter \bigskip \eightpoint

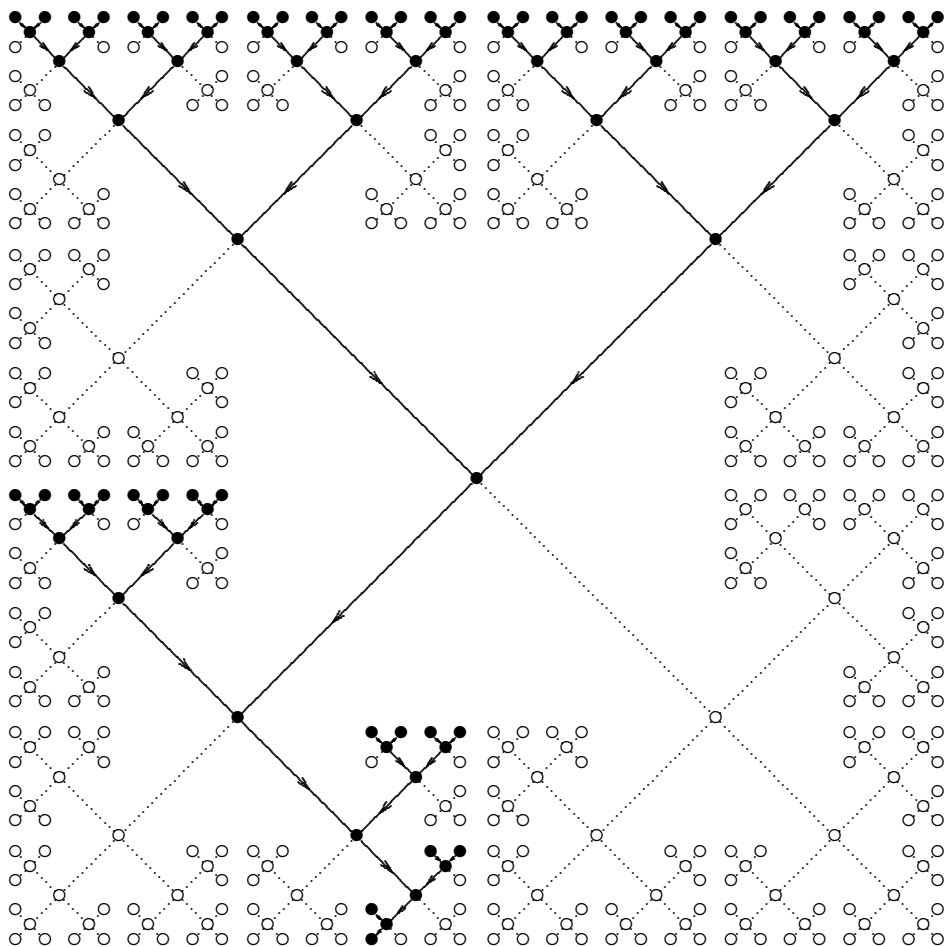
The Cayley graph of the free group picturing a river basin.
The positive generators are pointing downwards, towards
southwest and southeast, to give the idea of downstream.

\endgroup

\bigskip \bigskip \noindent form an element of the subshift $X$, but if we start anywhere in any tributary, and if we
decide to travel downstream (i.e.~following the edges of the Cayley graph corresponding to positive generators), we will
pick up another infinite word, which will merge into the main river, and which will also consist of an infinite word
belonging to $X$.

Although we do not believe it is possible to find a complete set of properties characterizing $\Spec $, we at least know
that it contains a dense copy of $X$ (not necessarily with the same topology), which may be precisely characterized and
which allows for a reasonably good handle on the other, more elusive elements of $\Spec $.  Seen from this perspective,
$\Spec $ appears slightly friendlier and we are in turn able to explore it quite efficiently.

While the methods most often used in the literature for analyzing $\CM $ are based on a rather technical study of a
certain AF-subalgebra (see e.g.~\cite [Corollary 3.9]{MatsuOri} and \cite [Theorem 4.14]{Thomsen}), our arguments are
rooted in the dynamical properties of the spectral partial action.  In particular our description of $\Spec $ is
concrete enough to allow us to find sensible necessary and sufficient conditions for this partial action to be
topologically free and minimal.  Our condition for minimality, for instance, is a lot similar to the cofinality
condition which has played an important role in characterizing simplicity for graph algebras \cite [Theorem 6.8]{KPRR},
\cite [Corollary 3.11]{KPR}.  These two crucial properties, namely topological freeness and minimality, have been
extensively used to characterize simplicity, and thus our treatment of $\CM $ is done in the same footing as for several
other better behaved C*-algebras.

We believe this new picture for the hitherto intractable spectrum of $\Core $ will allow for further advances in the
understanding of subshifts as well as of Carlsen-Matsumoto C*-algebras.

\section Subshifts

\label SubshSect We begin by fixing a nonempty finite set $\Lambda $ which will henceforth be called the \"{alphabet},
and whose elements will be called \"{letters}.

Any finite sequence of letters will be called a \"{finite word}, including the \"{empty word}, namely the word with
length zero, which will be denoted by $\emptyword $.  The set of all finite words will be denoted by $\Lambda ^*$.
Infinite sequences of letters will also be considered and we shall call them \"{infinite words}.

The best way to formalize the notion of \"{sequences}, twice referred to above, is by resorting to the Cartesian product
$\Lambda \times \Lambda \times \ldots \times \Lambda $, whose elements are therefore denoted by something like $$
(x_1,x_2,\ldots ,x_n).  $$ We shall however choose a more informal notation, denoting such a sequence by $$ x_1x_2\ldots
x_n, $$ or by $$ x_1x_2x_3\ldots $$ in the infinite case.  This is compatible with our point of view according to which
\"{sequences} are viewed as \"{words}.

If $\alpha $ is a finite word and $x$ is a finite or infinite word, then we will write $\alpha x$ for the concatenation
of $\alpha $ and $x$, namely the word obtained by juxtaposing $\alpha $ and $x$ together.

The \"{length} of a finite word $\alpha $, denoted $|\alpha |$, is the number of letters in it.

Assigning the discrete topology to $\Lambda $, the set of all infinite words, namely $\Lambda ^{\bf N}$, becomes a
compact topological space with the product topology by Tychonoff's Theorem.  The map $$ \shft : \Lambda ^{\bf N} \to
\Lambda ^{\bf N} $$ given by $$ \shft (x_1x_2x_3x_4\ldots ) = x_2x_3x_4x_5\ldots , \equationmark IntroduceShift $$ for
every $x=x_1x_2x_3x_4\ldots \in \Lambda ^{\bf N}$, is called the (left) \"{shift}.  It is easy to see that $\shft $ is
continuous.

Given a nonempty closed subset $X\subseteq \Lambda ^{\bf N}$ such that $$ \shft (X)\subseteq X, $$ we may consider the
restriction of $\shft $ to $X$, and then the pair $ (X,\shft |_X) $ is called a (one-sided) \"{subshift}.  Sometimes we
will also say that $X$ itself is a subshift, leaving the shift map to be deduced from the context.

There are many concrete situations in Mathematics where subshifts arise naturally such as in dynamical systems, Markov
chains, maps of the interval, billiards, geodesic flows, complex dynamics, information theory, automata theory and
matrix theory.  The present work is dedicated to studying subshifts from the point of view or partial dynamical systems
\cite {PDSFB}.

In what follows let us give an important example of subshifts.  A finite word $\alpha $ is said to \"{occur} in an
infinite word $x=x_1x_2x_3x_4\ldots \,$, if there are integer numbers $n\leq m$, such that $$ \alpha =x_nx_{n+1}\ldots
x_m.  $$ In other words, $\alpha $ may be found within $x$ as a contiguous block of letters.  By default we consider the
empty word $\emptyword $ as occurring in any infinite word.

Given an arbitrary subset $\Forbid \subseteq \Lambda ^*$, appropriately called the set of \"{forbidden words}, let $$ X
= X_\Forbid $$ be the set of all infinite words $x$ such that {\bf no} member of $\Forbid $ occurs in $x$.  It is easy
to see that $\shft (X_\Forbid )\subseteq X_\Forbid $, and a simple argument shows that $X_\Forbid $ is closed in
$\Lambda ^{\bf N}$, hence $X_\Forbid $ is a subshift.  If a subshift $X$ coincides with $X_\Forbid $, for some finite
set $\Forbid \subseteq \Lambda ^*$, then we say that $X$ is a \"{subshift of finite type}.

A notable example of a shift which is not of finite type, and which will often be used as a counter-example below, is as
follows: over the alphabet $\Lambda =\{0,1\}$, consider as forbidden all words of the form $$ 01^{^{2n+1}}0 =
0\underbrace {1\ldots 1}_{2n+1}0, $$ where $n\geq 0$.  Thus any odd string of 1's delimited by two 0's is forbidden.
The subshift defined by this set of forbidden words is called the \"{even shift}.  Clearly, an infinite word $x$ lies in
the space of the even shift if and only if, anytime a contiguous block of 1's occurring in $x$ is delimited by 0's,
there is an even number of said 1's.  It should be noted that an infinite word beginning with an \"{odd} (sic) number of
ones and followed by a zero is not immediately excluded.

Given a subshift $X$, the \"{language} of $X$ is defined to be the subset $$ \Lang \subseteq \Lambda ^* \equationmark
DefineLanguage $$ formed by all finite words which occur in some $x\in X$.

For future reference we cite here a well known result in Symbolic Dynamics:

\state Proposition \label AllSSForbWords \cite [Proposition 1.3.4]{LindMarcus} \ For any subshift $X$ one has that
$X=X_\Forbid $, where $\Forbid = \Lambda ^* \setminus \Lang $.

The following notions will prove to be of utmost importance in what follows:

\definition Let $\Lambda $ be a finite alphabet, and let $X\subseteq \Lambda ^{\bf N}$ be a subshift.  For each $\alpha
$ in $\Lambda ^*$, the \"{follower set} $\Folow _\alpha $, and the \"{cylinder set} $\Cyl _\alpha $ are defined by $$
\def \quad { } \matrix { \Folow _\alpha & = & \{y\in X: \alpha y\in X\}, \hbox { and}\hfill \cr \pilar {18pt} \Cyl
_\alpha & = & \{x\in X: x=\alpha y, \hbox { for some infinite word } y\}, } $$

It is well known that the $\Cyl _\alpha $ form a basis for the product topology on $X$, consisting of compact open
subsets.

Thus the follower set of $\alpha $ is the set of all infinite words which are allowed to \"{follow} $\alpha $, while the
cylinder $\Cyl _\alpha $ is the set of words which begin with the prefix $\alpha $.  If $|\alpha |=n$, then clearly
$S^n(\Cyl _\alpha ) = \Folow _\alpha $, and the restriction of $S^n$ to $\Cyl _\alpha $ gives a bijective map $$ S^n:
\Cyl _\alpha \to \Folow _\alpha .  $$

Since $\Cyl _\alpha $ is compact, the above is a homeomorphism from $\Cyl _\alpha $ to $\Folow _\alpha $.  In particular
it follows that $\Folow _\alpha $ is also compact, although it is not necessarily open as we shall now see in the
following well known result whose precise statement we have not been able to locate in the literature.

\state Proposition \label SFTOpen Given a finite alphabet $\Lambda $, and a subshift $X\subseteq \Lambda ^{\bf N}$, the
following are equivalent: \izitem \zitem $\Folow _\alpha $ is open in $X$ for every finite word $\alpha $ in $\Lambda
^*$, \zitem $\Folow _a$ is open in $X$ for every letter $a$ in $\Lambda $, \zitem $\shft $ is an open mapping on $X$,
\zitem $X$ is a subshift of finite type.

\Proof \ipfimply \pfimply (i)(ii) Obvious.

\pfimply (ii)(iii) Given an open subset $U\subseteq X$, we have for every $a$ in $\Lambda $, that $U\cap \Cyl _a$ is
open relative to $\Cyl _a$, and since $\shft $ is a homeomorphism from $\Cyl _a$ to $\Folow _a$, we conclude that $\shft
(U\cap \Cyl _a)$ is open relative to $\Folow _a$, which in turn is open relative to $X$ by hypothesis.  So $\shft (U\cap
\Cyl _a)$ is open relative to $X$, whence $$ \shft (U) = \shft \big (\medcup _{a\in A}U\cap \Cyl _a\big ) = \medcup
_{a\in A}\shft (U\cap \Cyl _a) $$ is open in $X$.

\pfimply (iii)(i) As already seen $\Cyl _\alpha $ is always open and $\Folow _\alpha =\shft ^n(\Cyl _\alpha )$, where
$n=|\alpha |$.

\BreakImply (iii)$\,\Leftrightarrow \,$(iv) This is a well known classical result in Symbolic Dynamics.  See \cite
[Theorem 1]{Parry}, \cite [Theorem 1]{ItoTaka} and \cite [Theorem 3.35]{Kurka}.  \endProof

\section Circuits

One of the most important aspects of subshifts to be discussed in this work is the question of \"{topological freeness}
(to be defined later), and which requires a careful understanding of \"{circuits}.  We will therefore set this section
aside to discuss this concept.  Besides the essential facts about circuits we shall need later, we will also present
some interesting, and perhaps unknown facts which came up in our research.

In order to motivate the concept of circuits let us first discuss the case of Markov subshifts.  So, for the time being,
we will let $X$ be a Markov subshift, which means that $X=X_\Forbid $, where the set $\Forbid $ of forbidden words
contains only words of length two.  Setting $$ a_{ij}=\left \{\matrix { 0,& \hbox {if the word \qt {ij} lies in }\Forbid
,\hfill \cr \pilar {13pt} 1,& \hbox {otherwise,} \hfill }\right .  $$ for every $i,j\in \Lambda $, the resulting 0-1
matrix $A=(a_{ij})_{i,j\in \Lambda }$ is called the \"{transition} matrix for $X$.  Excluding the uninteresting case in
which some letter of the alphabet is never used in any infinite word in $X$, the transition matrix has no zero rows.

Let $Gr(A)$ be the graph having $A$ as its adjacency matrix, so that its vertices are the elements of $\Lambda $, while
there is one edge from vertex $i$ to vertex $j$ when $a_{ij}=1$, and none otherwise.  The elements of $X$ may then be
thought of as infinite paths in $Gr(A)$, namely infinite sequences of vertices (as opposed to edges, another popular
concept) in which two successive vertices are joined by an edge.

\definition Given a Markov subshift $X$, with transition matrix $A$, we will say that a \"{circuit}\fn {Also known as a
\"{loop}, a \"{cycle}, or a \"{closed path}.} in $Gr(A)$ is a finite path $$ \gamma =x_1x_2\ldots x_n, $$ such that
$a_{x_n,x_1}=1$.

For each circuit $\gamma $, the infinite periodic path $$ \gamma ^\infty =\gamma \gamma \gamma \ldots $$ lies in $X$.
The question of topological freeness, already alluded to (but not yet defined), will be seen to be closely related the
the non-existence of circuits $\gamma $ such that $\gamma ^\infty $ is an isolated point of $X$.  If indeed $\gamma
^\infty $ is not isolated, that is, if $\{\gamma ^\infty \}$ is not an open set, we have in particular that, $$ \{\gamma
^\infty \} \neq \Cyl _\gamma .  $$ This is to say that $\Cyl _\gamma $ must contain at least one point other than
$\gamma ^\infty $.  In other words, it must be possible to prolong the finite word $\gamma $ in such a way as to obtain
an infinite path distinct from $\gamma ^\infty $.  Thus the relevance of an \"{exit} for $\gamma $, namely a letter
$x_i$ in the above expression for $\gamma $ which may be followed by at least one letter other than $x_{i+1}$ (in case
$i=n$, we take $x_{i+1}$ to mean $x_1$).

In a general subshift, not necessarily of finite type, the above considerations do not make any sense since there is no
underlying graph, but they may nevertheless be reinterpreted, as we will now see.

\definition Let $\Lambda $ be a finite alphabet, and let $X\subseteq \Lambda ^{\bf N}$ be a subshift.  \iaitem \aitem We
will say that a finite word $\gamma \in \Lambda ^*$ is a \"{circuit} (relative to $X$), provided the infinite periodic
word $\gamma ^\infty =\gamma \gamma \gamma \ldots \ $ lies in $X$.  \aitem We will say that an infinite word $y\neq
\gamma ^\infty $ is an \"{exit} for a given circuit $\gamma $, if $\gamma y$ lies in $X$.

Thus, to say that a circuit $\gamma $ has an exit is equivalent to saying that the follower set $\Folow _\gamma $ has at
least one element other than $\gamma ^\infty $.  If $\gamma $ is a circuit then $$ \gamma ^n:=\underbrace {\gamma \ldots
\gamma }_n $$ is also a circuit for any $n\in {\bf N}$.  However, even if $\gamma $ has an exit, there is no reason why
$\gamma ^n$ would also have an exit (unless $X$ is a Markov subshift).

An example of this situation is the shift (of finite type) $X_\Forbid $, on the alphabet $\{0,1\}$, where $\Forbid $
consists of a single forbidden word, namely $\qt {001}$. Evidently $\gamma =\qt {0}$ is a circuit, which admits the word
$y=\qt {111\ldots }$ as an exit.  However $\gamma ^2=\qt {00}$ clearly has no exit.

Having an exit is therefore no big deal.  A much more impressive property of a circuit $\gamma $ is for $\gamma ^n$ to
have an exit for every $n$.

\state Proposition \label OneExit Let $\Lambda $ be a finite alphabet and let $X\subseteq \Lambda ^{\bf N}$ be a
subshift.  Given a circuit $\gamma $, the following are equivalent: \izitem \zitem For each $n\in {\bf N}$, one has that
$\gamma ^n$ has an exit, \zitem There is an infinite word $z$ which is an exit for $\gamma ^n$, for all $n\in {\bf N}$.

\Proof The crucial difference between (i) and (ii), as the careful reader would have already noticed, is that in (i) it
is OK for each $\gamma ^n$ to have a different exit, while in (ii) it is required that there is one single exit which
works for all $\gamma ^n$.

It is obvious that \implica (ii)(i), so let us focus on the converse.  Assuming that $\gamma ^n$ has an exit for each
$n$, we have that the follower set $\Folow _{\gamma ^n}$ contains at least one element besides $(\gamma ^n)^\infty
=\gamma ^\infty $.  We claim that in fact $\Folow _{\gamma ^n}$ contains at least one point $z$ which does not lie in
the cylinder $\Cyl _\gamma $.

Notice that, even though $\gamma ^\infty $ lies in $\Folow _{\gamma ^n}$, it cannot be taken as the $z$ above since it
lies in the cylinder $\Cyl _\gamma $, while $z$ should not.

To prove the claim, choose $y\in \Folow _{\gamma ^n}$, with $y\neq \gamma ^\infty $, and write $$ y=\alpha _1\alpha
_2\alpha _3\ldots , $$ where each $\alpha _i$ is a finite word with the same length as $\gamma $.  If $\alpha _1\neq
\gamma $, then $y$ is not in $\Cyl _\gamma $, and there is nothing to be done.  Otherwise let $k$ be the smallest index
such that $\alpha _k\neq \gamma $, so $k>1$, and $\alpha _1=\alpha _2=\cdots =\alpha _{k-1}=\gamma $.  We may then write
$$ y = \gamma ^{k-1}z, $$ where $z=\alpha _k\alpha _{k+1}\ldots \ $ is therefore an infinite word not in $\Cyl _\gamma
$.  Since $y$ is in the follower set of $\gamma ^n$, we then have that $$ X\ni \gamma ^ny = \gamma ^n\gamma ^{k-1}z =
\gamma ^{k-1}\gamma ^nz.  $$

Recall that $X$ is invariant under the shift.  So, after applying $S$ to the above element $|\gamma |(k-1)$ times, we
conclude that $\gamma ^nz$ lies in $X$, so $z\in \Folow _{\gamma ^n}\setminus \Cyl _\gamma $, as desired.

Since $\Folow _{\gamma ^n}$ is closed and $\Cyl _\gamma $ is open, we have that $\Folow _{\gamma ^n}\setminus \Cyl
_\gamma $ is closed, and nonempty as seen above.  Observing that the $\Folow _{\gamma ^n}$ are decreasing with $n$, we
have by compactness of $X$ that $$ \medcap _{n\in {\bf N}} \Folow _{\gamma ^n}\setminus \Cyl _\gamma \neq \emptyset .
$$

Any element $z$ chosen in the above intersection is therefore not equal to $\gamma ^\infty $, because it is not in $\Cyl
_\gamma $, and it is therefore an exit for $\gamma ^n$, for all $n\in {\bf N}$.  \endProof

\definition \label DefStrongExit Let $\Lambda $ be a finite alphabet, and let $X\subseteq \Lambda ^{\bf N}$ be a
subshift.  If $\gamma $ is a circuit relative to $X$, we will say that an infinite word $y\neq \gamma ^\infty $ is a
\"{strong exit} for $\gamma $, if $\gamma ^ny \in X$, for all $n\in {\bf N}$.

It follows from \ref {OneExit} that a subshift in which all circuits have an exit, also satisfies the apparently
stronger property that all circuits have a \"{strong} exit (think about it).

An interesting example of such a subshift is as follows:

\state Proposition \label EvenExit Let $X$ be the even shift.  Then any circuit $\gamma $ relative to $X$ has a strong
exit.

\Proof Let us first assume that $$ \gamma =1^k=\underbrace {1\ldots 1}_k, $$ for some $k\geq 1$.  Then, no matter how
big is $n$, we may always exit $\gamma ^n$ via an infinite string of $0$'s.  All other circuits $\gamma $ have $\qt {0}$
somewhere, say $\gamma =\gamma '0\gamma ''$.  Observing that $$ \gamma ^n\gamma '011111\ldots $$ lies in the space of
the even shift, one may always exit $\gamma ^n$ via the infinite word $\gamma '011111\ldots \ $, which is therefore a
strong exit for $\gamma $.  \endProof

In section \ref {TopFreeSect} we will carefully study topological freeness for the partial dynamical systems we shall
encounter along the way.  However we warn the reader that the nice property of the even shift proven above (existence of
strong exits for all circuits) will be seen to be still insufficient for our purposes.

\section The standard partial action associated to a subshift

Throughout this section we will fix a finite alphabet $\Lambda $ and a subshift $X\subseteq \Lambda ^{\bf N}$.

As already seen, the shift restricts to a homeomorphism from $\Cyl _a$ to $\Folow _a$, for each $a$ in $\Lambda $.  The
inverse of this homeomorphism is clearly given by the map $$ \theta _a:\Folow _a\to \Cyl _a, $$ defined by $$ \theta
_a(y) = ay, \for y\in \Folow _a.  $$

\definition \label FirstPaction Let $\F = \F (\Lambda )$ be the free group on $\Lambda $ and let $$ \bigTheta = \big
(\{X_g\}_{g\in \F }, \{\theta _g\}_{g\in \F }\big ) $$ be the unique semi-saturated partial action of $\F $ on $X$
assigning $\theta _a$ to each $a$ in $\Lambda $, given by \book {Proposition}{4.10}.  Henceforth $\theta $ will be
referred to as the \"{standard partial action} associated to $X$.

Incidentally, to say that $\bigTheta $ is \"{semi-saturated} is to say that $$ \theta _{gh}=\theta _g\circ \theta _h, $$
whenever $|gh|=|g|+|h|$, where $|{\cdot }|$ is the usual length function on $\F $.  See \book {Definition}{4.9}.

For the case of Markov subshifts, in fact for a generalization thereof, the standard partial action was first studied in
\cite {infinoa}.

Recall from \book {Definition}{5.1} that, for a partial action on a topological space, it is usually required that the
$X_g$ be \"{open sets}.  However, since $$ X_{a\inv } = \hbox {Domain}(\theta _a) = \Folow _a, \for a\in \Lambda , $$
that requirement is not fulfilled for our $\bigTheta $ unless $X$ is a subshift of finite type by \ref {SFTOpen}.
Rather than a nuisance, this is the first indication that a non-finite type subshift conceals another partial action
which will be seen to be crucial for the analysis we will carry out later.

When considering $\bigTheta $, it is therefore best to think of it as a partial action in the category of sets (as
opposed to topological spaces).

In what follows let us give a simple description for $\bigTheta $, but first let us introduce some notation.  Denote by
$\F _+$ the subsemigroup of $\F $ generated by $\Lambda \cup \{1\}$, so that $\F _+$ may be identified with the set
$\Lambda ^*$ of all finite words in $\Lambda $.  Under this identification we shall see the empty word $\emptyword $ as
the unit of $\F $.

\state Proposition \label DescrPAction For every $g$ in $\F $ one has: \izitem \zitem If $g \not \in \PPinv $, then
$X_g$ is the empty set.  \zitem If $g=\alpha \beta \inv $, with $\alpha ,\beta \in \F _+$, and $|g|=|\alpha |+|\beta |$,
then $$ X_g=\{\alpha y\in X: y\in \Folow _\alpha \cap \Folow _\beta \}.  $$ \zitem If $g$ is as in (ii), and $x\in
X_{g\inv }$, write $x=\beta y$, for some $y\in \Folow _\beta \cap \Folow _\alpha $.  Then $$ \theta _g(x) = \theta
_{\alpha \beta \inv }(\beta y) = \alpha y.  $$

\Proof Given $g$ in $\F $, write $$ g=c_1c_2\ldots c_n, $$ with $c_i\in \Lambda \cup \Lambda \inv $ in reduced form,
meaning that $c_{i+1}\neq c_i\inv $, for every $i$.  Since $\theta $ is semi-saturated we have that $$ \theta _g=\theta
_{c_1}\circ \theta _{c_2}\circ \cdots \circ \theta _{c_n}.  $$

If $g \notin \PPinv $ then there is some $i$ such that $c_i\in \Lambda \inv $, and $c_{i+1}\in \Lambda $, say $c_i=a\inv
$ and $c_{i+1}=b$, for some $a,b\in \Lambda $.  Therefore $a\neq b$, and then $$ \theta _{c_i}\circ \theta _{c_{i+1}}=
\theta _a\inv \theta _b $$ is the empty map, because after inserting $b$ as the prefix of an infinite word $x$, we are
left with a word which does not begin with the letter $a$!  Under these conditions we then have that $\theta _g$ is the
empty map, hence $X_g$ is the empty set.  This proves (i).

Given a finite word $\alpha $, it is easy to see that $X_{\alpha \inv }=\Folow _\alpha $, $X_\alpha =\Cyl _\alpha $, and
$$ \theta _\alpha (x) = \alpha x, \for x\in X_{\alpha \inv }.  $$ If $g$ is as in (ii) we have, again by
semi-saturatedness that $\theta _g=\theta _\alpha \circ \theta _\beta \inv $, from where (ii) and (iii) follow easily.
\endProof

We should remark that the sets $X_g$, appearing in \ref {DescrPAction.ii} above, have also played an important role in
Carlsen's study of subshifts \cite [Definition 1.1.3]{CarlsenNotes}.

Due to \ref {DescrPAction.ii} there will be numerous situations below in which we will consider group elements of the
form $g=\alpha \beta \inv $, with $\alpha ,\beta \in \F _+$, and $$ |g|=|\alpha |+|\beta |.  $$ To express this fact we
will simply say that $g$ is in \"{reduced form}.  This is clearly equivalent to saying that the last letter of $\alpha $
is distinct from the last letter of $\beta $, so that no cancellation takes place.

Observe also that ``reduced form" is an attribute of the \"{presentation} of $g$ as a product of two elements, rather
than of $g$ itself.  Any element of $\PPinv $ may be written as $g=\alpha \beta \inv $, with $\alpha ,\beta \in \F _+$,
and upon canceling as many final letters of $\alpha $ and $\beta $ as necessary, one is left with a reduced form
presentation of $g$.

\section A partial representation associated to the subshift

Throughout this section we will again fix a finite alphabet $\Lambda $ and a subshift $X\subseteq \Lambda ^{\bf N}$.

Recall from \cite [9.1\amper 9.2]{PDSFB} that a \"{*-partial representation} of a group $G$ in a unital *-algebra $B$ is
a map $u:G\to B$ satisfying \begingroup \parindent 40pt \medskip \item {(PR1)} $u_1 = 1$, \medskip \item {(PR2)} $u_g
u_h u_{h\inv } = u_{gh} u_{h\inv }$.  \medskip \item {(PR3)} $u_{g\inv } = (u_g)^*$, \endgroup \medskip \noindent for
every $g$ and $h$ in $G$.

Consider the complex Hilbert space $\ell ^2(X)$, with its canonical orthonormal basis $\{\delta _x\}_{x\in X}$.  Here we
shall consider *-partial representations of $\F = \F (\Lambda )$ in the algebra of bounded linear operators on $\ell
^2(X)$, and we will refer to these simply as \"{partial representations of\/ $\F $ on $\ell ^2(X)$}.

\state Proposition \label IntroPRep For each $g$ in $\F $, denote by $\pr _g$ the unique bounded linear operator $$ \pr
_g: \ell ^2(X) \to \ell ^2(X), $$ such that for each $x$ in $X$, $$ \pr _g(\delta _x) = \left \{ \matrix { \delta
_{\theta _g(x)},& \hbox {if } x\in X_{g\inv }, \cr \pilar {12pt} 0, & \hbox {otherwise,}\hfill } \right .  $$ where
$\theta $ is the standard partial action associated to $X$.  Then the correspondence $g\mapsto \pr _g$ is a
semi-saturated partial representation of\/ $\F $ on $\ell ^2(X)$.

\Proof Given $g$ and $h$ in $\F $, and $x$ in $X$, we must prove that $$ \pr _g\pr _h\pr _{h\inv }(\delta _x) = \pr
_{gh}\pr _{h\inv }(\delta _x).  \equationmark CondForPrep $$

Suppose first that $x\in X_h\cap X_{g\inv }$.  Then the left-hand-side above equals $$ \pr _g\pr _h\pr _{h\inv }(\delta
_x) = \pr _g\pr _h\big (\delta _{\theta _{h\inv }(x)}\big ) = \pr _g(\delta _x) = \delta _{\theta _g(x)}.  $$ Observing
that $$ \theta _{h\inv }(x)\in \theta _{h\inv }\big (X_h\cap X_{g\inv }\big ) = X_{h\inv }\cap X_{(gh)\inv }, $$ the
right-hand-side of \ref {CondForPrep} equals $$ \pr _{gh}\pr _{h\inv }(\delta _x) = \pr _{gh}\big (\delta _{\theta
_{h\inv }(x)}\big ) = \delta _{\theta _{gh}(\theta _{h\inv }(x))} = \delta _{\theta _g(x)}, $$ proving \ref
{CondForPrep} in the present case.

Suppose now that $x\notin X_h\cap X_{g\inv }$.  In case $x\notin X_h$, then $\pr _{h\inv }(\delta _x)=0$, and \ref
{CondForPrep} follows trivially.  So we are left with the case that $x\in X_h\setminus X_{g\inv }$.  Under this
assumption we have $$ \pr _g\pr _h\pr _{h\inv }(\delta _x) = \pr _g\pr _h\big (\delta _{\theta _{h\inv }(x)}\big ) = \pr
_g(\delta _x) = 0.  $$ Moreover notice that $ y:= \theta _{h\inv }(x) \notin X_{(gh)\inv }, $ since otherwise $$ x =
\theta _h(y) \in \theta _h\big (X_{h\inv }\cap X_{(gh)\inv }\big ) = X_h\cap X_{g\inv }, $$ contradicting our
assumptions.  Therefore, $$ 0=\pr _{gh}(\delta _y) = \pr _{gh}(\delta _{\theta _{h\inv }(x)})= \pr _{gh}\pr _{h\inv
}(\delta _x), $$ showing that the right-hand-side of \ref {CondForPrep} also vanishes.  This proves \ref {CondForPrep},
also known as (PR2), and we leave the easy proofs of (PR1) and (PR3), as well as the fact that $S$ is semi-saturated,
for the reader.  \endProof

Either analyzing the definition of $\pr _g$ directly, or as a consequence of the above result (see \cite
[9.8.i]{PDSFB}), we have that $\pr _g$ is a partial isometry for every $g$ in $\F $, with initial space $\pilar
{11pt}\stake {7pt}\ell ^2(X_{g\inv })$ and with final space $\ell ^2(X_g)$.

More specifically, recall from \ref {DescrPAction.i} that $\theta _g$ is the empty map when $g$ is not in $\PPinv $, in
which case $$ \pr _g = 0.  \equationmark SgZero $$

On the other hand, assuming that $g$ lies in $\PPinv $, we may write $ g=\alpha \beta \inv , $ in reduced form, with
$\alpha ,\beta \in \F _+$, whence $$ \pr _g = \pr _\alpha \pr _\beta ^*, $$ by semi-saturatedness.  Based on the
description of $\theta $ given in \ref {DescrPAction.ii\amper iii}, we then have for every $x$ in $X$ that $$ \pr
_{\alpha \beta \inv }(\delta _x) = \left \{ \matrix { \delta _{\alpha y},& \hbox {if } x=\beta y, \hbox { for some }
y\in \Folow _\alpha \cap \Folow _\beta , \cr \pilar {12pt} 0, & \hbox {otherwise.}\hfill } \right .  $$

Still under the assumption that $g=\alpha \beta \inv $, in reduced form, we have by \ref {DescrPAction.ii} that the
final space of $\pr _g$ is $$ \ell ^2(X_g) = \overline {\hbox {span}} \{\delta _{\alpha y}: y\in \Folow _\alpha \cap
\Folow _\beta \}.  \equationmark FinalSpace $$

A relevant remark is that the final projections $\pr _g\pr _g^*$ are then seen to be diagonal operators relative to the
canonical orthonormal basis.  In particular these commute with each other, a well known fact from the general theory of
partial representations \cite [9.8.iv]{PDSFB}.

A special case of interest is when $g=\alpha \in \F _+$, in which case we have that $$ \pr _\alpha (\delta _x) = \left
\{ \matrix { \delta _{\alpha x},& \hbox {if } x\in \Folow _\alpha , \hfill \cr \pilar {12pt} 0, & \hbox {otherwise.}  }
\right .  $$ From this we see that the initial space of $\pr _\alpha $ is $\ell ^2(\Folow _\alpha )$, and its final
space is $\ell ^2(\Cyl _\alpha )$.

The adjoint of $\pr _\alpha $ may be described by $$ \pr _\alpha ^*(\delta _y) = \pr _{\alpha \inv }(\delta _y) = \left
\{ \matrix { \delta _{x},& \hbox {if } y\in \Cyl _\alpha , \hbox { and } y=\alpha x, \hbox { with } x\in \Folow _\alpha
, \cr \pilar {12pt} 0, & \hbox {otherwise.} \hfill } \right .  $$

\state Proposition \label IntroT Let $T$ be the operator on $\ell ^2(X)$ defined by $ T = \sum _{a\in \Lambda } \pr
_a^*.  $ Then $$ T(\delta _x) =\delta _{\shft (x)}, $$ for every $x$ in $X$, where $\shft $ is the shift map introduced
in \ref {IntroduceShift}.

\Proof Given any $y$ in $X$, let $b$ be the first letter of $y$, and write $y=bx$, for some infinite word $x$,
necessarily in $\Folow _b$.  Then $\pr _b^*(\delta _y)=\delta _x$, while $\pr _a^*(\delta _y)=0$, for all $a\neq b$.
This implies that $$ T(\delta _y) = \medsum _{a\in \Lambda } \pr _a^*(\delta _y) = \delta _x = \delta _{\shft (y)}.
\endProof

The operator $T$ may therefore be interpreted as the manifestation of the shift $\shft $ at the level of operators on
$\ell ^2(X)$.  Since each $\pr _a$ is a partial isometry, it is clear that $\Vert T\Vert \leq |\Lambda |$.

\section C*-algebras associated to subshifts

Throughout this section we will fix a finite alphabet $\Lambda $ and a subshift $X\subseteq \Lambda ^{\bf N}$.  Our goal
here is to describe two important C*-algebras that have been extensively studied in association with a subshift.

\state Proposition \label GenForMatsu Let $\pi $ be the representation of\/ $C(X)$ on $\ell ^2(X)$ defined on the
canonical orthonormal basis by $$ \pi (f)\delta _x = f(x)\delta _x, \for f\in C(X), \for x\in X.  $$ Then the following
three sets generate the same C*-algebra of operators on $\ell ^2(X)$: \izitem \zitem $\pr (\F )=\{\pr _g: g\in \F \}$,
\zitem $\pr (\Lambda ) = \{\pr _a: a\in \Lambda \}$, \zitem $\pi \big (C(X)\big ) \cup \{T\}$.

\Proof Denote the C*-algebras generated by the sets in (i), (ii) and (iii) by $B_\F $, $B_\Lambda $ and $B_T$,
respectively.

Since $\pr $ is semi-saturated, for every $g$ in $\F $, one has that $\pr _g$ may be written as a product of elements in
$\pr (\Lambda ) \cup \pr (\Lambda )^*$.  Therefore $B_\F \subseteq B_\Lambda $.

For every $\alpha $ in $\Lambda ^*$, the cylinder $\Cyl _\alpha $ is a clopen subset of $X$, so its characteristic
function, which we denote by $1_\alpha $, is a continuous function.  Moreover $\pi (1_\alpha )$ is the orthogonal
projection onto $\ell ^2(\Cyl _\alpha )$, hence it coincides with the final projections of $\pr _\alpha $, so $$ \pi
(1_\alpha )= \pr _\alpha \pr _\alpha ^*\in B_\F .  $$

It is evident that the set $\{1_\alpha : \alpha \in \Lambda ^*\}$ separates points of $X$, so by virtue of the
Stone-Weierstrass Theorem, it generates $C(X)$, as a C*-algebra.  Therefore $$ \pi \big (C(X)\big )\subseteq B_\F .  $$
Since $T\in B_\F $, by definition, we then conclude that $$ B_T\subseteq B_\F .  $$

Given $a$ in $\Lambda $, notice that $$ \pi (1_a)T^* = \medsum _{b\in \Lambda } \pr _a\pr _a^*\pr _b = \pr _a, $$
proving that $\pr _a$ lies in the C*-algebra generated by $\pi \big (C(X)\big )$ and $T$, whence $B_\Lambda \subseteq
B_T$.

We have therefore shown that $$ B_\F \subseteq B_\Lambda \subseteq B_T\subseteq B_\F , $$ whence the conclusion.
\endProof

\definition The Matsumoto C*-algebra associated to a given subshift $X$, henceforth denote by $\MatsAlg $, is the closed
*-algebra of operators on $\ell ^2(X)$ generated by any one of the sets described in the statement of \ref
{GenForMatsu}.

The algebra defined above was first introduced by Matsumoto in \cite [Lemma 4.1]{MatsuAutomorph} under the notation $\O
_X$.  However, due to the existence of several different algebras associated with subshifts usually denoted by $\O _X$
in the literature, we would rather use a new notation.  See also remark \ref {Zoo}, below.

Let $\lambda $ denote the left regular representation of the free group $\F =\F (\Lambda )$ on $\ell ^2(\F )$.  Thus,
denoting the canonical orthonormal basis of $\ell ^2(\F )$ by $\{\delta _g\}_{g\in \F }$, one has that $$ \lambda
_g(\delta _h) = \delta _{gh}, \for g,h\in \F .  $$

Regarding the partial representation $\pr $ introduced in \ref {IntroPRep}, we may define a new partial representation
$\prt $ of $\F $ on $\ell ^2(X)\otimes \ell ^2(\F )$, by tensoring $\pr $ with $\lambda $, namely $$ \prt _g = \pr
_g\otimes \lambda _g, \for g\in \F .  $$

\state Proposition \label GenForMatsuCar Let $\TT $ be the operator on $\ell ^2(X)\otimes \ell ^2(\F )$ given by $\TT =
\sum _{a\in \Lambda } \prt _a^*$.  Then the following three sets generate the same C*-algebra of operators on $\ell
^2(X)\otimes \ell ^2(\F )$: \izitem \zitem $\prt (\F )$, \zitem $\prt (\Lambda )$, \zitem $\big (\pi \big (C(X)\big
)\otimes 1\big ) \cup \{\TT \}$.

\Proof Denote the C*-algebras generated by the sets in (i), (ii) and (iii) by $\tilde B_\F $, $\tilde B_\Lambda $ and
$\tilde B_T$, respectively.  That $\tilde B_\F \subseteq \tilde B_\Lambda $ follows, as above, from the fact that $\prt
$ is semi-saturated.  For $\alpha $ in $\Lambda ^*$ one has $$ \pi (1_\alpha )\otimes 1 = \pr _\alpha \pr _\alpha
^*\otimes 1 = (\pr _\alpha \otimes \lambda _\alpha )(\pr _\alpha \otimes \lambda _\alpha )^* = \prt _\alpha \prt _\alpha
^*\in \tilde B_\F .  $$ This, plus the fact that $C(X)$ is generated by $\{1_\alpha : \alpha \in \Lambda ^*\}$, gives $
\pi \big (C(X)\big )\otimes 1\subseteq \tilde B_\F .  $ Since $\TT \in \tilde B_\F $, by definition, we then conclude
that $ \tilde B_T\subseteq \tilde B_\F .  $ Given $a$ in $\Lambda $, notice that $$ \big (\pi (1_a)\otimes 1\big )\TT ^*
= (\pr _a\pr _a^*\otimes 1) \medsum _{b\in \Lambda } \pr _b\otimes \lambda _b = \medsum _{b\in \Lambda } \pr _a\pr
_a^*\pr _b\otimes \lambda _b = \pr _a\otimes \lambda _a = \prt _a, $$ proving that $\prt _a$ lies in the C*-algebra
generated by $\pi \big (C(X)\big )\otimes 1$ and $\TT $, whence $\tilde B_\Lambda \subseteq \tilde B_T$, concluding the
proof.  \endProof

\definition \label DefineCM The Carlsen-Matsumoto C*-algebra associated to a given subshift $X$, henceforth denote by
$\CM $, is the closed *-algebra of operators on $\ell ^2(X)\otimes \ell ^2(\F )$ generated by any one of the sets
described in the statement of \ref {GenForMatsuCar}.

The algebra defined above was studied in \cite {CarlsenSilvestrov}, \cite {CarlsenCuntzPim}, \cite {CarlsenNotes} (where
it was also denoted by $\O _X$).  Its definition in the above mentioned references is not the one given in \ref
{DefineCM}, but we will prove in \ref {AlgebrasSame} that the two definitions lead to the same algebra.  We nevertheless
note that the description of $\CM $ given by \ref {GenForMatsuCar.iii} is closely related to that given in \cite
{CarlsenSilvestrov}.

\state Remark \label Zoo \rm Let us clarify the relationship between our notation and the one used in the literature.
According to \cite {CarlsenSilvestrov} there are at least three possibly non isomorphic C*-algebras associated to a
subshift in the literature, some of them defined only for two-sided subshifts.  These are: \iaitem \aitem the C*-algebra
$\O _\Lambda $ defined in \cite {MatsuOri}, \aitem the C*-algebra $\O _\Lambda $ defined in \cite {MatsuCarl}, \aitem
the C*-algebra $\O _X$ defined in \cite {CarlsenCuntzPim}.  \medskip \noindent We shall not be concerned with the first
algebra above, while the second and third ones are respectively being denoted by $\MatsAlg $ and $\CM $ in this work.

Notice that, although the notation $\CM $ emphasizes the dynamical system $(S,X)$, the construction of this algebra was
actually made based on the specific way in which $X$ is represented as a space of infinite words, as well as on the
syntactic rules of inserting a letter ahead of an infinite word.

It is then legitimate to ask what is the precise relationship between $\CMSymbol _{X_1}$ and $\CMSymbol _{X_2}$, in case
$X_1$ and $X_2$ are conjugate subshifts.  This question has been addressed for the first time in \cite {Matsumoto} where
it was proved that $\MSymbol _{X_1}$ is Morita equivalent to $\MSymbol _{X_2}$ under suitable hypothesis.  Later Carlsen
\cite [Theorems 1.3.1 \amper 1.3.5]{CarlsenNotes} found a proof of the fact that $\CMSymbol _{X_1}$ and $\CMSymbol
_{X_2}$ are Morita equivalent without any additional hypotheses.  In \cite [Section 8]{CarlsenCuntzPim}, having realized
$\CM $ as a Cuntz-Pimsner algebra, Carlsen finally proved that $\CMSymbol _{X_1}$ and $\CMSymbol _{X_2}$ are actually
isomorphic.  A further, much simpler proof of this result was found by Carlsen and Silvestrov \cite [Section
11]{CarlsenSilvestrov}, based on the description of $\CM $ as a crossed product by an endomorphism and a transfer
operator \cite {NewLook}.

The following proof, based on \cite {CarlsenSilvestrov}, is perhaps the simpler possible proof of the invariance of
Matsumoto's algebras, given that it only needs the original description of these algebras.

\state Proposition Suppose that $X_1$ and $X_2$ are conjugate subshifts.  Then $$ \MSymbol _{X_1}\simeq \MSymbol _{X_2}.
$$

\Proof Let $ \varphi :X_1 \to X_2 $ be a homeomorphism such that $\varphi \circ S_1=S_2\circ \varphi $, where $S_1$ and
$S_2$ are the shift maps on $X_1$ and $X_2$, respectively.  Define a unitary operator $U:\ell ^2(X_1)\to \ell ^2(X_2)$,
by setting $$ U(\delta _x) = \delta _{\varphi (x)}, \for x\in X_1.  $$

Decorating all of the ingredients of \ref {GenForMatsu} with subscripts to indicate whether we are speaking of $X_1$ or
$X_2$, it is easy to see that $$ U\pi _1\big (C(X_1)\big )U^*=\pi _2\big (C(X_2)\big ), \and UT_1U^*=T_2.  $$ By \ref
{GenForMatsu} we then see that $U\MSymbol _{X_1}U^* = \MSymbol _{X_2}$, so $\MSymbol _{X_1}$ and $\MSymbol _{X_2}$ are
in fact \"{spatially} isomorphic.  \endProof

Even though the above result is probably well known to the experts, we have not been able to locate it in the literature
with this exact formulation.

Given our reliance on the free group, whose rank is definitely not an invariant of the subshift, our method does not
seem appropriate to prove Carlsen's invariance Theorem \cite [Theorem 8.6]{CarlsenCuntzPim}.

\section The spectrum

\label SpectrumSection As always we fix a finite alphabet $\Lambda $ and a subshift $X\subseteq \Lambda ^{\bf N}$.
There is an important subalgebra of $\MatsAlg $ which has played a crucial role in virtually every attempt to study
subshifts from the point of view of C*-algebras (see the references given in the introduction), and which we would now
like to describe.

\definition We will denote by $\Core $ the closed *-algebra of operators on $\ell ^2(X)$ generated by the final
projections $$ e_g:= \pr _g\pr _g^*, $$ for all $g$ in $\F $.

Recall from \cite [9.8.iv]{PDSFB} that the $e_g$ commute with each other, so $\Core $ is a commutative C*-algebra, which
is moreover unital because $e_1=1$.

Since $\prt $ is also a partial representation, the final projections $$ \et _g := \prt _g\prt _g^*, $$ likewise commute
with each other, and hence generate a commutative C*-algebra.  However notice that $$ \et _g = (\pr _g\otimes \lambda
_g)(\pr _g\otimes \lambda _g)^* = \pr _g\pr _g^*\otimes 1 = e_g\otimes 1, \equationmark ETTensorOne $$ from where one
concludes that the C*-algebra generated by the $\et _g$ is nothing but $\Core \otimes 1$.

This section is dedicated to studying the spectrum of $\Core $, henceforth denoted by $$ \SpecAlg .  $$

The full description of this space requires some machinery still to be developed, but we may easily give examples of
some of its elements, as follows.  Observing that the multiplication of two diagonal operators is done by simply
multiplying the corresponding diagonal entries, we see that the assignment of a given diagonal entry to an operator
defines a multiplicative linear functional on the set of all diagonal operators.  In what follows we will refer to
diagonal entries indirectly, as eigenvalues relative to eigenvectors taken from the canonical basis.

\definition \label DefinePhiX Given any $x$ in $X$, let $\varphi _x$ be the unique linear functional on $\Core $ such
that $$ a(\delta _x) = \varphi _x(a)\delta _x, \for a\in \Core .  $$

As observed above, each $\varphi _x$ is a character on $\Core $, hence an element of $\SpecAlg $.  We will see that not
every element of $\SpecAlg $ is of the form $\varphi _x$, but the $\varphi _x$ nevertheless form a large subset of
$\SpecAlg $ in the following sense:

\state Proposition \label XIsDense The subset of\/ $\SpecAlg $ formed by all the $\varphi _x$ is a dense set.

\Proof Assume by way of contradiction that the closure of $$ \{\varphi _x:x\in X\}, $$ which we denote by $C$, is a
proper subset of $\SpecAlg $.  Picking a point $\varphi $ outside $C$ we may invoke Urysohn's Lemma to find a continuous
complex valued function $f$ on $\SpecAlg $ which vanishes on $C$, and such that $f(\varphi )=1$.  By Gelfand's Theorem
we have that $f$ is the Gelfand transform of some $a$ in $\Core $, and then for every $x$ in $X$ we have $$ \varphi
_x(a) = f(\varphi _x) =0.  $$ This implies that all diagonal entries of $a$ are zero, and since $a$ is itself a diagonal
operator, we deduce that $a=0$, and hence also that $f=0$, a contradiction.  \endProof

We thus get a map $$ \Phi : x\in X\mapsto \varphi _x\in \SpecAlg , $$ whose range is dense in $\SpecAlg $.

A crucial question in this subject is whether or not $\Phi $ is continuous.  Should this be the case, the compactness of
$X$ would imply that $\Phi $ is onto, and the mystery surrounding $\SpecAlg $ would be immediately dispelled.  However,
we will see later in \ref {OntoContinuousSFT} that for subshifts not of finite type $\SpecAlg $ is strictly bigger than
the range of $\Phi $ and this can only happen if $\Phi $ is discontinuous!

In order to describe the whole of $\SpecAlg $ it is useful to recall that $\Core $ is generated, as a C*-algebra, by the
projections $$ e_g:=\pr _g\pr _g^*, $$ so a character $\varphi $ on $\Core $ is pinned down as soon as we know the
numbers $\varphi (e_g)$, which necessarily lie in $\{0,1\}$, for every $g$ in $\F $.

To be precise, given $\varphi \in \SpecAlg $, consider the element $\xi _\varphi \in 2^\F $ given by $$ \xi _\varphi (g)
= \varphi (e_g), \for g\in \F .  $$

Identifying $2^\F $ with the set of all subsets of $\F $ as usual, we may think of $\xi _\varphi $ as the subset of $\F
$ given by $$ \xi _\varphi = \{g\in \F : \varphi (e_g)=1\}.  \equationmark SubsetModel $$

\state Proposition \label XiHomeo Considering $2^\F $ as a topological space with the product topology, the mapping $$
\SO : \varphi \in \SpecAlg \mapsto \xi _\varphi \in 2^\F , $$ is a homeomorphism from $\SpecAlg $ onto its range.

\Proof As already observed, each $\varphi $ in $\SpecAlg $ is characterized by its values on the generating idempotents,
so $\SO $ is seen to be one-to-one.  It is evident that $\SO $ is continuous, and since $\SpecAlg $ is compact, we have
that $\SO $ is a homeomorphism onto its image.  \endProof

\definition \label IntroduceSpec The range of $\SO $, which will henceforth be denoted by $$ \Spec = \SO (\SpecAlg ), $$
will be referred to as the \"{spectrum} of the subshift $X$.

As seen in \ref {XiHomeo}, we have that $\Spec $ is homeomorphic to $\SpecAlg $, and we will take the former as our main
model to study the latter.  Since $\SpecAlg $ contains a dense copy of the set\fn {We use \"{set} as opposed to
\"{topological space} to highlight the fact that it is not necessarily homeomorphic to $X$, as $\Phi $ may be
discontinuous.}  $X$ by \ref {XIsDense}, we have that $\Spec $ also contains a dense copy of $X$.  However we will see
that for subshifts not of finite type, $\Spec $ is strictly bigger than $X$.

We will generally prefer to regard a given $\xi $ in $2^\F $ as a subset of $\F $, in the spirit of \ref {SubsetModel},
rather than as a $\{0,1\}$-valued function on $\F $.

Here is the list of properties, alluded to in the introduction, that every $\xi $ in $\Spec $ must satisfy:

\state Proposition \label PropForPaint For any $\xi $ in $\Spec $ one has that: \izitem \zitem $1\in \xi $, \zitem $\xi
$ is convex, \zitem for every $g\in \xi $, there exists a unique $a$ in $\Lambda $, such that $ga\in \xi $, \zitem if
$g\in \xi $, and $\alpha $ is a finite word in $\Lambda $ such that $g\alpha \in \xi $, then $\alpha $ lies in $\Lang $.
\zitem $\xi \subseteq \PPinv $.

\def \convexity {\ref {PropForPaint.ii}}

\Proof Let $\varphi $ be the character on $\Core $ such that $\xi =\xi _\varphi $, so that $$ \xi = \{g\in \F : \varphi
(e_g)=1\}.  $$ Since $e_1=1$, one has that $\varphi (e_1)=1$, so $1\in \xi $, proving (i).

Recall from \book {Definition}{14.19} that to say that $\xi $ is \"{convex} is to say that, whenever $g,h\in \xi $, then
the \"{segment} joining $g$ and $h$, namely $$ \overline {gh}:= \{k\in \F : |g\inv h| = |g\inv k| + |k\inv h|\} $$ is
contained in $\xi $.  Assuming that $g$ and $h$ are in $\xi $, and that $k$ is in $\overline {gh}$, set $$ s=g\inv k,
\and t=k\inv h.  $$ We then have that $|st| = |s|+|t|$, so $\pr _{st}=\pr _s\pr _t$, by semi-saturatedness.  Employing
\book {Proposition}{14.5} we conclude that $ e_{st} \leq e_s, $ which is to say that $$ e_{g\inv h} \leq e_{g\inv k}.
$$ Conjugating the left-hand-side of this inequation by $\pr _g$, we get by \cite [9.8.iii]{PDSFB} that $$ \pr _g
e_{g\inv h} \pr _{g\inv } = e_h \pr _g \pr _{g\inv } = e_he_g.  $$ Doing the same relative to the right-hand-side leads
to $$ \pr _g e_{g\inv k} \pr _{g\inv } = e_ke_g, $$ so we deduce that $ e_he_g \leq e_ke_g, $ whence $$ \varphi (e_he_g)
\leq \varphi (e_ke_g).  $$

Having assumed that $g,h\in \xi $, we see that $\varphi (e_g) = \varphi (e_h) = 1$, so the left-hand-side above
evaluates to 1.  The same is therefore true for the right-hand-side, which implies that $\varphi (e_k)=1$, from where it
follows that $k\in \xi $, as desired.  This proves (ii).

In order to prove (iii), pick any $g$ in $\xi $.  Observing that $X$ is the disjoint union of the cylinders $\Cyl _a$,
as $a$ range in $\Lambda $, and that $e_a$ is the orthogonal projection onto $\ell ^2(\Cyl _a)$, we see that $$ \medsum
_{a\in \Lambda }e_a=1.  $$ Conjugating the above identity by $\pr _g$, we deduce that $$ e_g = \pr _g\pr _{g\inv } = \pr
_g\big (\medsum _{a\in \Lambda } e_a\big )\pr _{g\inv } = \medsum _{a\in \Lambda } \pr _ge_a\pr _{g\inv } = \medsum
_{a\in \Lambda } e_{ga}\pr _g\pr _{g\inv } = \medsum _{a\in \Lambda } e_{ga}e_g.  $$ Since $g\in \xi $, we have that
$\varphi (e_g)=1$, so $$ 1 = \varphi (e_g) = \medsum _{a\in \Lambda } \varphi (e_{ga}e_g) = \medsum _{a\in \Lambda }
\varphi (e_{ga}).  $$

Each $e_{ga}$ is idempotent so $\varphi (e_{ga})$ is either 0 or 1.  We then see that there exists a unique $a$ in
$\Lambda $ such that $\varphi (e_{ga})=1$, meaning that $ga\in \xi $, hence proving (iii).

Supposing, as in (iv), that $g$ and $g\alpha $ lie in $\xi $, observe that $$ e_ge_{g\alpha } = \pr _g\pr _{g\inv
}e_{g\alpha } = \pr _ge_{\alpha }\pr _{g\inv }.  $$ Since $\varphi (e_ge_{g\alpha })=1$, by hypothesis, it follows that
$e_ge_{g\alpha }\neq 0$, whence also $e_\alpha \neq 0$.  Observing that $e_\alpha $ is the orthogonal projection onto
$\ell ^2(\Cyl _\alpha )$, one deduces that $\Cyl _\alpha $ is nonempty, which implies that $\alpha $ is a word in the
language $\Lang $.

Regarding (v), let $g\in \xi $.  Then $\varphi (e_g)=1$, so $e_g$ is nonzero and hence neither is $\pr _g$.  It then
follows from \ref {SgZero} that $g$ lies in $\PPinv $.  \endProof

It should be stressed that properties \ref {PropForPaint.i-v}, which we have seen to hold for every $\xi $ in $\Spec $,
are not enough to characterize the elements in $\Spec $.  Although it would be highly desirable to find a set of
properties giving such a precise characterization, we have not been able to succeed in this task.

Nevertheless, recalling that the image of $X$ under $\Phi $ is dense in $\SpecAlg $ by \ref {XIsDense}, we already have
a somewhat satisfactory description of $\Spec $, as the closure of $X$, or rather, of its image under the following
composition of maps $$ \def \quad {\kern 8pt} \def \larw #1{\quad {\buildrel #1 \over \longrightarrow }\quad } X \larw
\Phi \SpecAlg \larw \SO \Spec \subseteq 2^\F .  $$

The map described above will acquire a special relevance in what follows, so it deserves a special notation:

\definition \label BigCompos We will denote by $$ \XO :X\to \Spec $$ the map given by the above composition, namely $\XO
=\SO \circ \Phi $.

Unraveling the appropriate definitions it is easy to see that $$ \XO (x) = \{g\in \F : e_g(\delta _x)=\delta _x\}.  $$
Notice that to say that $e_g(\delta _x)=\delta _x$ is the same as saying that that $\delta _x$ lies in the final space
of $\pr _g$, which we have seen to be $\ell ^2(X_g)$.  Thus we may alternatively describe $\XO (x)$ as $$ \XO (x) =
\{g\in \F : x\in X_g\}.  \equationmark EquivCond $$

We will now give a further, more detailed, description of $\XO (x)$.

\state Proposition \label CharacXix Given $x$ in $X$, let $$ \xi _x=\XO (x).  $$ Then $\xi _x$ consists precisely of the
elements $g$ in $\F $ such that the following conditions hold: \izitem \zitem $g$ may be written in reduced form as
$\alpha \beta \inv $, with $\alpha ,\beta \in \F _+$, \zitem $\alpha $ is a prefix of $x$ and, writing $x=\alpha y$, one
has that $y\in \Folow _\alpha \cap \Folow _\beta $.

\Proof Given $g\in \xi _x$ we may use \ref {PropForPaint.v} to write $g=\alpha \beta \inv $, and we may clearly assume
that (i) above holds.  By \ref {EquivCond} we have that $x\in X_g$, and hence (ii) follows from \ref {DescrPAction.ii}.

Conversely, assuming that $g$ satisfies (i) and (ii), we have again by \ref {DescrPAction.ii} that $x$ lies in $X_g$,
proving that $g\in \xi _x$.  \endProof

We shall often use the above result in the special case that $\alpha =\emptyword $, in which case it reads:

\state Corollary \label BinvInXix Given $x$ in $X$, and $\beta $ in $\Lambda ^*$, one has that $$ \beta \inv \in \xi _x
\iff x\in \Folow _\beta .  $$

It is instructive to view these elements within the Cayley graph of $\F $.

\bigskip \null \hskip 2cm 
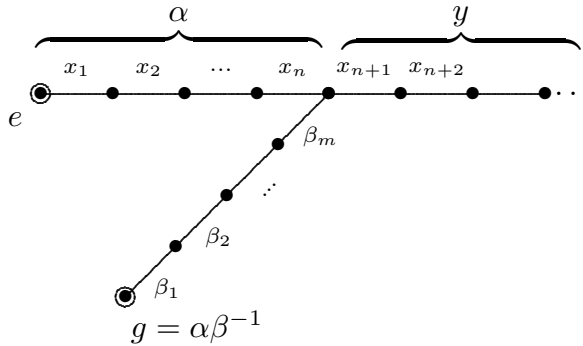
\beginpicture \setcoordinatesystem units <0.0080truecm, -0.0080truecm> point at 1000 1000 \put
{$e$} at -42 42 \circulararc 360 degrees from 0 16 center at 0 0 \put {$\bullet $} at 0 0 \plot 0 0 120 0 / \put
{$\bullet $} at 120 0 \put {$\scriptstyle x_1$} at 60 -40 \plot 120 0 240 0 / \put {$\bullet $} at 240 0 \put
{$\scriptstyle x_2$} at 180 -40 \plot 240 0 360 0 / \put {$\bullet $} at 360 0 \put {$\scriptstyle \ldots $} at 300 -40
\plot 360 0 480 0 / \put {$\bullet $} at 480 0 \put {$\scriptstyle x_n$} at 420 -40 \plot 480 0 600 0 / \put {$\bullet
$} at 600 0 \put {$\scriptstyle x_{n+1}$} at 540 -40 \plot 600 0 720 0 / \put {$\bullet $} at 720 0 \put {$\scriptstyle
x_{n+2}$} at 660 -40 \plot 720 0 840 0 / \put {$\bullet $} at 840 0 \put {$\ldots $ } at 895 0 \put {$\overbrace {\hbox
to 90pt{\hfill }}^{\textstyle \alpha }$} at 230 -110 \put {$\overbrace {\hbox to 75pt{\hfill }}^{\textstyle y}$} at 700
-110 \put {$\bullet $} at 141 339 \plot 141 339 225 255 / \put {$\bullet $} at 225 255 \put {$\scriptstyle \beta _1$} at
211 325 \plot 225 255 310 170 / \put {$\bullet $} at 310 170 \put {$\scriptstyle \beta _2$} at 296 240 \plot 310 170 395
85 / \put {$\bullet $} at 395 85 \put {$\scriptstyle \nedots $} at 381 156 \plot 395 85 480 0 / \put {$\bullet $} at 480
0 \put {$\scriptstyle \beta _m$} at 466 71 \put {$g=\alpha \beta \inv $} at 261 389 \circulararc 360 degrees from 141
355 center at 141 339 \endpicture \hfill \null

\medskip \bigskip \noindent Given any $x$ in $X$, the goal is to mark the vertices of the Cayley graph of $\F $
corresponding to the elements of $\xi _x$.  Due to \ref {PropForPaint.i}, we must always mark the unit group element.
Thereafter, beginning at the unit group element we mark all vertices according to the successive letters of
$x=x_1x_2x_3\ldots \ $, thus forming the \"{stem} of $\xi _x$.  We then choose an integer $n$ and focus on the $n^{th}$
vertex along this path, letting $$ \alpha =x_1x_2\ldots x_n, \and y=x_{n+1}x_{n+2}\ldots , $$ so that $x=\alpha y$, and
$y\in \Folow _\alpha $.  Choosing a finite word $\beta $, we then back up starting at the vertex chosen above, along the
letters of $\beta $.  However, before we do this, we must make sure that the infinite word $\beta y$ lies in $X$, which
the same as saying that $y\in \Folow _\beta $.  It is also best to choose $\beta $ such that $\beta _m\neq x_n$, since
this will avoid stepping on a vertex we have already traversed, guaranteeing \ref {CharacXix.i}.

The group element $g=\alpha \beta \inv $ is therefore an element of $\xi _x$ and, as seen in \ref {CharacXix}, all
elements of $\xi _x$ arise in this way.

One may also think of the stem of $\xi _x$ as a \"{river}, the $\beta $'s considered in the above diagram being its
\"{tributaries}, while $\xi _x$ consists of the whole \"{river basin}.

In fact, not only the $\xi _x$, but every $\xi $ in $\Spec $ has an interpretation as a river basin, but first we need
to identify the appropriate \"{rivers}.

\state Proposition \label PropStem Let $\xi $ be in $\Spec $ and let $g\in \xi $.  Then there exists a unique $x$ in $X$
such that $$ \xi \cap g\F _+ = \{g\alpha : \alpha \hbox { is a prefix of } x\}.  $$

\Proof By induction and \ref {PropForPaint.iii}, there exists an infinite sequence $x=x_1x_2x_3x_4\ldots \ $, with
$x_i\in \Lambda $, such that $$ gx_1x_2\ldots x_n\in \xi , \for n\in {\bf N}.  $$ Using \ref {PropForPaint.iv} we have
that $x_1x_2\ldots x_n\in \Lang $, for every $n$, so it follows from \ref {AllSSForbWords} that $x\in X$.  The inclusion
``$\supseteq $" relative to the sets in the statement then clearly holds.

On the other hand, given any $g\alpha \in \xi \cap g\F _+$, we claim that $\alpha $ is a prefix of $x$.  In order to
prove this, suppose otherwise, and let $\alpha '$ be the shortest prefix of $\alpha $ which is not a prefix of $x$.  By
{\convexity } we have that $g\alpha '\in \xi $, which is to say that we may assume without loss of generality that
$\alpha $ is already minimal.  Write $\alpha =y_1y_2\ldots y_n$, with $y_i$ in $\Lambda $, so that $$ \beta
:=y_1y_2\ldots y_{n-1} $$ is a prefix of $x$ by minimality.  We then have that $\beta x_n$ is a prefix of $x$, whence
$g\beta x_n$ lies in $\xi $.  But $g\alpha =g\beta y_n$ also lies in $\xi $, so $x_n=y_n$, by \ref {PropForPaint.iii},
whence $$ \alpha =\beta y_n = \beta x_n $$ is a prefix of $x$, a contradiction.  \endProof

The following concept is reminiscent of \cite [Definition 5.5]{infinoa}.

\definition \label DefineStem Given $\xi $ in $\Spec $ and $g$ in $\xi $, the unique $x$ in $X$ satisfying the
conditions of \ref {PropStem} will be called the \"{stem} of $\xi $ at $g$, and it will be denoted by $\stem _g(\xi )$.
In the special case that $g=1$, we will refer to $x$ simply as the \"{stem} of $\xi $, denoting it by $\stem (\xi )$.

Given $x$ in $X$, we have that all prefixes of $x$ lie in $\xi _x$ by \ref {CharacXix}.  From this it immediately
follows that the stem of $\xi _x$ is precisely $x$.  In symbols $$ \stem (\xi _x)=x.  $$ This shows the following:

\state Proposition \label StemLeftInverse The stem, viewed as a map $$ \stem :\Spec \mapsto X, $$ is a left inverse for
$\XO $.  Consequently $\XO $ is one-to-one and $\stem $ is onto.

We have already hinted at the fact that the map $$ \Phi :x\in X\mapsto \varphi _x\in \SpecAlg $$ may not be continuous,
in which case neither is $\XO $.  Fortunately, not all of the maps in sight are discontinuous:

\state Proposition \label StemContinuous The stem defines a continuous mapping from $\Spec $ to $X$.

\Proof Given a net $\{\xi _i\}_i$ in $\Spec $ converging to some $\xi $, let $$ x_i=\stem (\xi _i),\and x=\stem (\xi ).
$$ In order to prove the statement we need to show that $\{x_i\}_i$ converges to $x$.

By the definition of the product topology on $X\subseteq \Lambda ^{\bf N}$, given any neighborhood $U$ of $x$, we may
find a cylinder $\Cyl _\alpha $, with $$ x\in \Cyl _\alpha \subseteq U.  $$ It follows that $\alpha $ is a prefix of
$x$, whence $\alpha $ belongs to $\xi $.

By the definition of the product topology on $\Spec \subseteq 2^\F $, for every $g$ in $\F $, the function $$ \eta \in
\Spec \mapsto [g\in \eta ]\in \{0,1\}, \equationmark BooleanContinuous $$ where the brackets correspond to Boolean
value, is continuous.  Therefore $$ 1=[\alpha \in \xi ] = \lim _i\ [\alpha \in \xi _{x_i}], $$ which means that $\alpha
\in \xi _{x_i}$ for all sufficiently large $i$.  Consequently $\alpha $ is a prefix of $x_i$, so $$ x_i\in \Cyl _\alpha
\subseteq U.  $$ This concludes the proof.  \endProof

Although $\stem $ is onto, it might not be one-to-one.  This is to say that there may be many different elements in
$\Spec $ with the same stem.  An example of this situation is obtained by taking any $\xi $ together with $\xi _x$,
where $x=\stem (\xi )$.  The following result further explores the relationship between these two elements.

\state Proposition \label XiXMaximal Given $\xi $ in $\Spec $, let $x=\stem (\xi )$.  Then $\xi \subseteq \xi _x$.

\Proof Given $g$ in $\xi $, we may write $g=\alpha \beta \inv $, with $\alpha ,\beta \in \F _+$ by \ref
{PropForPaint.v}, and we may clearly assume that $g$ is in reduced form, so \ref {CharacXix.i} holds.  By {\convexity }
we have that $\alpha \in \xi \cap \F _+$, so $\alpha $ is a prefix of $x$.

Write $\alpha =x_1x_2\ldots x_n$, and $x=\alpha x_{n+1}x_{n+2}\ldots $, and notice that, for every integer $k>n$, one
has that $\alpha x_{n+1}\ldots x_k$ is a prefix of $x$, so $$ \xi \ni \alpha x_{n+1}\ldots x_k = \alpha \beta \inv \beta
x_{n+1}\ldots x_k = g\beta x_{n+1}\ldots x_k.  $$

From \ref {PropForPaint.iv} it follows that $\beta x_{n+1}\ldots x_k$ lies in $\Lang $, for every $k$, whence the
infinite word $$ z = \beta x_{n+1}x_{n+2}\ldots $$ lies in $X$ by \ref {AllSSForbWords}.  We deduce that the infinite
word $$ y = x_{n+1}x_{n+2}\ldots $$ lies in $\Folow _\beta $ and clearly also in $\Folow _\alpha $.  This proves \ref
{CharacXix.ii}, so $$ g=\alpha \beta \inv \in \xi _x.  \endProof

As a consequence we see that the \"{river basin} picture of the $\xi _x$ also applies to a general element $\xi $ in
$\Spec $.  That is, if $x=\stem (\xi )$, then $\xi $ contains the whole ``river" $x$, while it is contained in the
``river basin'' $\xi _x$ by \ref {XiXMaximal}.

In particular, conditions \ref {CharacXix.i-ii}, which are \"{necessary} and \"{sufficient} for a given group element
$g$ to lie in $\xi _x$, are seen to still be \"{necessary} for membership in $\xi $, as long as we take $x$ to be the
stem of $\xi $.  By this we mean that when $g\in \xi $, then necessarily $g\in \xi _x$, whence said conditions hold.

Incidentally, there is a special case in which these conditions are also \"{sufficient}, as we shall now see.

\state Proposition \label SuficCritInterior Given $\xi $ in $\Spec $, suppose that $\alpha $ and $\beta $ are elements
of\/ $\F _+$ satisfying conditions \ref {CharacXix.i-ii} relative to $x=\stem (\xi )$.  Suppose, in addition, that the
element $y$ referred to in condition \ref {CharacXix.ii} actually lies in the \underbar {interior} of $\Folow _\beta $.
Then $\alpha \beta \inv \in \xi $.

\Proof Since $\Phi (X)$ is dense in $\SpecAlg $, we have that $\XO (X)$ is dense in $\Spec $, so we may write $\xi =\lim
_i\xi _{x_i}$, with $x_i$ in $X$.  By hypothesis we have that $\alpha $ is a prefix of $x$, so $\alpha \in \xi $, and by
the continuity of the maps described in \ref {BooleanContinuous}, we have that $\alpha \in \xi _{x_i}$, for all
sufficiently large $i$.  It follows that $\alpha $ is a prefix of $x_i$, so $$ x_i=\alpha y_i, $$ for some infinite word
$y_i$, necessarily belonging to $\Folow _\alpha $.  Since $x_i$ converges to $x$ by \ref {StemContinuous}, we have that
$y_i$ converges to $y$, which belongs to the interior of $\Folow _\beta $, by hypothesis.  So the $y_i$ lie in $\Folow
_\beta $, again for all sufficiently large $i$.  We conclude that $\alpha \beta \inv $ satisfies \ref {CharacXix.ii}
relative to all such $x_i$, whence $\alpha \beta \inv \in \xi _{x_i}$.  By continuity of Boolean values we then have $$
[\alpha \beta \inv \in \xi ] = \lim _i\ [\alpha \beta \inv \in \xi _{x_i}] = 1, $$ so $\alpha \beta \inv \in \xi $.
\endProof

This result has the following important consequence (see also \cite [Remark 1.1.5]{CarlsenNotes}):

\state Theorem \label OntoContinuousSFT Given a subshift $X$, consider the mapping $$ \XO :x\in X\mapsto \xi _x\in \Spec
, $$ already defined in \ref {BigCompos}.  Then the following are equivalent: \izitem \zitem $X$ is a subshift of finite
type, \zitem $\XO $ is onto, \zitem $\XO $ is continuous, \zitem $\XO $ is a homeomorphism.

\Proof \ipfimply \pfimply (i)(ii) Given $\xi $ in $\Spec $, let $x$ be its stem, and we claim that fact $\xi =\xi _x$.
On the one hand we have that $\xi \subseteq \xi _x$ by \ref {XiXMaximal}.  In order to prove the reverse inclusion, pick
$g$ in $\xi _x$.  By \ref {CharacXix} we may find a decomposition $g=\alpha \beta \inv $ satisfying conditions \ref
{CharacXix.i-ii} so that, among other things, $x=\alpha y$, with $y\in \Folow _\beta $.

Since $X$ is of finite type, we have by \ref {SFTOpen} that $\Folow _\beta $ is open, so $y$ automatically lies in the
interior of $\Folow _\beta $, and we deduce from \ref {SuficCritInterior} that $g=\alpha \beta \inv \in \xi $.  This
shows that $\xi =\xi _x$, and hence that $\XO $ is onto.

\pfimply (ii)(iii) Recall from \ref {StemLeftInverse} that $\stem $ is one-to-one.  If we assume, in addition, that $\XO
$ is onto, then $\XO $ is an invertible map whose left-inverse $\stem $ is also its two-sided inverse and hence
invertible.  Since $\stem $ is a continuous invertible map defined on a compact set, it must be a homeomorphism, so its
inverse, namely $\XO $, is then seen to be continuous.

\pfimply (iii)(i) Observe that \ref {BinvInXix} allows for the following description of $\Folow _\beta $: $$ \Folow
_\beta = \{x\in X: \beta \inv \in \xi _x\} = \XO \inv \big ( \{\xi \in \Spec : \beta \inv \in \xi \}\big ), $$ which is
then the inverse image of an open set under the continuous mapping $\XO $.  We conclude that $\Folow _\beta $ is open
and then (i) follows from \ref {SFTOpen}.

We have not taken (iv) into account so far, but it is clear that it implies either (ii) or (iii), while (ii)+(iii)
easily implies (iv) since $\XO $ is one-to-one by \ref {StemLeftInverse}.  \endProof

\section The spectral partial action

Recall from \ref {FirstPaction} that $\bigTheta $ is a partial action of $\F $ on $X$ which might however not be a
\"{topological} partial action in the sense that not all of the $X_g$ are open sets.  To supersede this badly behaved
partial action we will now show that there exists a fully compliant topological partial action of $\F $ on $\Spec $
extending $\bigTheta $.  The idea will be to build an algebraic partial action on $\Core $ and then consider the
corresponding action at the level of the spectrum.

For each $g$ in $\F $, let $D_g$ be the two-sided ideal of $\Core $ generated by $e_g$, namely $$ D_g=e_g\Core , $$ and
let $\tau _g$ be the map from $D_{g\inv }$ to $D_g$ given by $$ \tau _g(a) = \pr _ga\pr _{g\inv }, \for a\in D_{g\inv }.
\equationmark DefineTau $$ As in \book {Proposition}{10.1} one may show that $$ \tau = \big (\{D_g\}_{g\in \F }, \{\tau
_g\}_{g\in \F }\big ) \equationmark PactOnAlg $$ is a C*-algebraic partial action (see \book {Definition}{11.4}) of $\F
$ on $\Core $.  In fact \book {Proposition}{10.1} is proved in a purely algebraic context, but the proof given there
carries over to the C*-algebraic setting.

By \book {Corollary}{11.6}, there exists a topological partial action $$ \vartheta = \big (\{\Omega _g\}_{g\in \F },
\{\vartheta _g\}_{g\in \F }\big ) $$ of $\F $ on $\Spec $ (here identified with the spectrum of $\Core $) linked to
$\alpha $ via the fact that $D_g$ consists of the elements of $\Core $ (here identified with the algebra $C(\Spec )$ of
all continuous complex valued functions on $\Spec $ by Gelfand's Theorem) vanishing off $\Omega _g$, plus the relation
$$ \tau _g(f)|_\xi = f\big (\vartheta _{g\inv }(\xi )\big ), $$ for every $g\in \F $, $f\in D_{g\inv }$, and $\xi \in
\Omega _g$.

\definition \label SpectralAction We will refer to the partial action $\vartheta $ introduced above as the \"{spectral
partial action} associated to $X$.

As we shall see, the spectral partial action may be seen as a the restriction of the partial Bernoulli action of $\F $
\book {Definition}{5.12}.

\state Proposition \label DescrPAOnOmega Regarding the spectral partial action associated to $X$, one has that \izitem
\zitem $\Omega _g=\{\xi \in \Spec : g\in \xi \}$, \zitem $\vartheta _g(\xi ) = g\xi = \{gh: h\in \xi \}$, \medskip
\noindent for each $g$ in $\F $, and for each $\xi \in \Omega _{g\inv }$.

\Proof Under the well known correspondence between closed two-sided ideals in $C(\Spec )$ and open sets in $\Spec $,
notice that the ideal generated by an idempotent element corresponds to the support of the latter.  In case of the
idempotent $e_g$, its support, initially viewed in $\SpecAlg $, corresponds to the set of all characters $\varphi $ such
that $\varphi (e_g)=1$.  Identifying $\SpecAlg $ and $\Spec $ by \ref {IntroduceSpec}, and noting that $$ \varphi
(e_g)=1 \explain {SubsetModel}\iff g\in \xi _\varphi , $$ for all characters $\varphi $, one sees that (i) follows.

Given $g$ in $\F $, and $\xi \in \Omega _{g\inv }$, let $\varphi $ be the character on $\Core $ corresponding to $\xi $
under $\SO $, namely such that $\xi =\xi _\varphi $.  Since $e_{g\inv }\in \xi $, by \ref {DescrPAOnOmega.i}, we have
that $\varphi (e_{g\inv })=1$.  Moreover, $\vartheta _g(\xi )$ will correspond to a character $\psi $, such that $\psi
(e_g)=1$, and $$ \psi (a) = \varphi \big (\tau _{g\inv }(a)\big ), \for a\in D_g.  $$ For any $a$ in $\Core $,
regardless of whether $a$ is in $D_g$ or not, we have that $ae_g\in D_g$, whence $$ \psi (a) = \psi (e_g)\psi (a)= \psi
(e_ga)= \varphi \big (\tau _{g\inv }(e_ga)\big ).  $$

We then have that $\vartheta _g(\xi )=\xi _\psi $, so for any $h$ in $\F $, one has that $$ h\in \vartheta _g(\xi ) \iff
h\in \xi _\psi \explain {SubsetModel}\iff \psi (e_h)=1.  \equationmark TBCont $$ Incidentally notice that $$ \psi (e_h)
= \varphi \big (\tau _{g\inv }(e_he_g)\big ) = \varphi \big (u_{g\inv }e_he_gu_g)\big ) = \varphi \big (e_{g\inv
h}e_{g\inv }\big )= \varphi (e_{g\inv h}), $$ because $\varphi (e_{g\inv })=1$.  Focusing on \ref {TBCont}, we then have
that $$ \psi (e_h)=1 \iff \varphi (e_{g\inv h}) = 1 \iff g\inv h\in \xi \iff h\in g\xi , $$ proving that $ \vartheta
_g(\xi )=g\xi , $ as desired.  \endProof

We then have partial dynamical systems $$ \bigTheta = \big (\{X_g\}_{g\in \F }, \{\theta _g\}_{g\in \F }\big ), \and
\vartheta = \big (\{\Omega _g\}_{g\in \F }, \{\vartheta _g\}_{g\in \F }\big ) \equationmark TwoSystems $$ on $X$ and
$\Spec $, respectively, and it is interesting to notice that these sets are related to each other by the maps $$ \XO
:X\to \Spec , \and \stem :\Spec \to X, $$ introduced in \ref {BigCompos} and \ref {DefineStem}, respectively.

\state Proposition \label Equivariance Both $\XO $ and $\stem $ are equivariant maps \book {Definition}{2.7} relative to
the partial dynamical systems in \ref {TwoSystems}.

\Proof Let us first prove that $$ \XO (X_g)\subseteq \Omega _g, \for g\in \F .  $$ Given $x\in X_g$, we have by \ref
{EquivCond} that $g\in \XO (x)$, whence $\XO (x)\in \Omega _g$, by \ref {DescrPAOnOmega.i}, proving the above inclusion.
We next must show that $$ \XO \big (\theta _g(x)\big ) = \vartheta _g\big (\XO (x)\big ), \equationmark FirstEquivar $$
for all $g$ in $\F $, and all $x\in X_{g\inv }$.  Given such an $x$, let $ y= \theta _g(x).  $ Then evidently $y\in
X_g$, and for any $h\in \F $, one has that $$ h\in \XO (y) \explain {EquivCond}\iff y\in X_h \iff y\in X_h\cap X_g \iff
\theta _g(x)\in X_h\cap X_g \$\iff x\in \theta _{g\inv }(X_h\cap X_g) = X_{{g\inv }h}\cap X_{g\inv } \iff x\in X_{{g\inv
}h} \explain {EquivCond}\iff {g\inv }h\in \XO (x).  $$ From this we see that $$ g\inv \XO (y)=\{{g\inv }h: h\in \XO
(y)\} =\XO (x), $$ so $$ \XO \big (\theta _g(x)\big ) = \XO (y) = g\XO (x) \={DescrPAOnOmega.ii} \vartheta _g\big (\XO
(x)\big ), $$ showing \ref {FirstEquivar}, and hence concluding the proof of the equivariance of $\XO $.

To prove that $$ \stem (\Omega _g)\subseteq X_g, \for g\in \F , $$ notice that if $g$ is not in $\PPinv $, then $g\notin
\xi $ for every $\xi $ in $\Spec $, by \ref {PropForPaint.v}.  It follows that $\Omega _g$ is empty and then there is
nothing to prove.  Otherwise write $g=\alpha \beta \inv $, with $\alpha ,\beta \in \F _+$, in reduced form.  Given $\xi
$ in $\Omega _g$, we have that $$ g\in \xi \explain {XiXMaximal}\subseteq \xi _x, $$ where $x=\stem (\xi )$.  It then
follows from \ref {CharacXix} that $x=\alpha y$, where $y\in \Folow _\alpha \cap \Folow _\beta $, which implies that
$x\in X_g$, by \ref {DescrPAction.ii}.  This proves that $\stem (\Omega _g)\subseteq X_g$, as desired.

We must finally prove that $$ \sigma \big (\vartheta _g(\xi )\big ) = \theta _g\big (\sigma (\xi )\big ), $$ for all $g$
in $\F $, and all $\xi $ in $\Omega _{g\inv }$.  Again there is nothing to do unless $g$ lies in $\PPinv $, so we may
assume that $g=\alpha \beta \inv $, with $\alpha ,\beta \in \F _+$, in reduced form.

Since $g\inv =\beta \alpha \inv \in \xi $, we have that $\beta \in \xi $, by {\convexity } so also $\beta \in \xi \cap
\F _+$.  By \ref {PropStem} the latter set consists precisely of the prefixes of the stem of $\xi $, which we will
denote by $x$ from now on.  In particular we have that $x=\beta y$, for some infinite word $y$, and we have by \ref
{DescrPAction.iii} that $$ \theta _g\big (\sigma (\xi )\big ) = \theta _{\alpha \beta \inv }\big (\beta y\big ) = \alpha
y.  $$

On the other hand, observe that if $\gamma $ is any prefix of $y$, then $\beta \gamma $ is a prefix of $x$, so $\beta
\gamma $ lies in $\xi $.  Consequently $$ \alpha \gamma = g\beta \gamma \in g\xi = \vartheta _g(\xi ).  $$ Since $\alpha
\gamma \in \F _+$ for any such $\gamma $, one sees that $\alpha \gamma $ is a prefix of the stem of $\vartheta _g(\xi
)$, from where it follows that $$ \stem \big (\vartheta _g(\xi )\big ) = \alpha y = \theta _g\big (\sigma (\xi )\big ).
\endProof

In case $X$ is a subshift of finite type the above result may be combined with \ref {OntoContinuousSFT} leading up to
the following:

\state Proposition \label EquivPact If $X$ is a subshift of finite type then the maps $$ \stem :\Spec \to X, \and \XO
:X\to \Spec $$ are mutually inverse equivariant homeomorphisms, whence the spectral partial action and the standard
partial action are equivalent.

Let us conclude this section with an important technical result, regarding membership of an element of the form $g\beta
\inv $, when we already know that $g$ is a member of a given $\xi $.  This result actually belongs in section \ref
{SpectrumSection}, but it was delayed up to now since its proof is greatly facilitated by the existence of the spectral
action.

\state Proposition \label Enfiabilidade Given $\xi $ in $\Spec $, pick $g\in \xi $, and let $\beta $ be a finite word.
Regarding the statements: \iaitem \aitem $g\beta \inv \in \xi $, and \aitem $\stem _g(\xi )\in \Folow _\beta $, \medskip
\noindent \vskip -1.3cm \hskip 4cm 
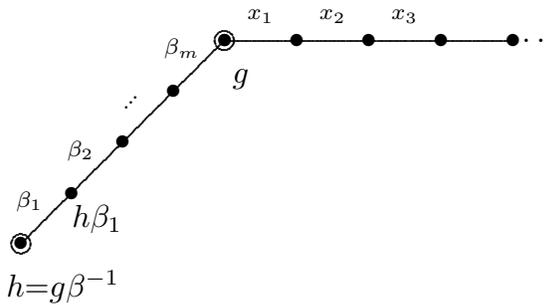
\beginpicture \setcoordinatesystem units <0.0080truecm, -0.0080truecm> point at 1000
1000 \put {$g$} at 27 65 \circulararc 360 degrees from 0 16 center at 0 0 \put {$\bullet $} at 0 0 \plot 0 0 120 0 /
\put {$\bullet $} at 120 0 \put {$\scriptstyle x_1$} at 60 -40 \plot 120 0 240 0 / \put {$\bullet $} at 240 0 \put
{$\scriptstyle x_2$} at 180 -40 \plot 240 0 360 0 / \put {$\bullet $} at 360 0 \put {$\scriptstyle x_3$} at 300 -40
\plot 360 0 480 0 / \put {$\bullet $} at 480 0 \put {$\ldots $ } at 535 0 \put {$\bullet $} at -339 339 \plot -339 339
-255 255 / \put {$\bullet $} at -255 255 \put {$\scriptstyle \beta _1$} at -325 269 \plot -255 255 -170 170 / \put
{$\bullet $} at -170 170 \put {$\scriptstyle \beta _2$} at -240 184 \plot -170 170 -85 85 / \put {$\bullet $} at -85 85
\put {$\scriptstyle \nedots $} at -156 99 \plot -85 85 0 0 / \put {$\bullet $} at 0 0 \put {$\scriptstyle \beta _m$} at
-71 14 \put {$h{=}g\beta \inv $} at -269 410 \circulararc 360 degrees from -339 355 center at -339 339 \put {$h\beta
_1$} at -212 297 \endpicture \hfill \null \medskip \noindent one has that \implica (a)(b).  Moreover, if $\xi =\xi _x$
for some $x$ in $X$, then \implica (b)(a) as well.

\Proof Let $h=g\beta \inv $, let $x=\stem _g(\xi )$, and write $$ x=x_1x_2x_3\ldots , \and \beta =\beta _1\beta _2\ldots
\beta _m.  $$

Notice that $h\beta _1$ lies in the segment joining $h$ and $g$ so, supposing (a), we have by {\convexity } that $
h\beta _1\in \xi .  $ Likewise $$ h\beta _1\beta _2\ldots \beta _i \in \xi , $$ for all $i$.  Furthermore, for every
integer $j$ we have that $$ h\beta _1\beta _2\ldots \beta _mx_1x_2\ldots x_j = gx_1x_2\ldots x_j\in \xi , $$ by the
definition of the stem of $\xi $ at $g$.  It then follows that $\stem _h(\xi )=\beta x$, so in particular $\beta x\in
X$, whence $x\in \Folow _\beta $, proving (b).

In order to prove the last sentence of the statement, let $\xi =\xi _x$ for some $x$ in $X$, and suppose that (b) holds.
In the special case that $g=1$, we have that $$ x=\stem (\xi _x) = \stem _g(\xi ) \in \Folow _\beta , $$ so $\beta \inv
$ belongs to $\xi $ by \ref {BinvInXix}, proving (a).

Dropping the assumption that $g$ is the unit group element, let us deal with the general case.  Observing that $g\in \xi
_x$, we have by \ref {EquivCond} that $x\in X_g$, so $$ y:= \theta _{g\inv }(x) $$ is well defined.  Moreover, by \ref
{Equivariance} $$ g\xi _y = \vartheta _g(\xi _y) = \xi _{\theta _g(y)} = \xi _x.  $$ In other words, $\xi _x$ is
obtained by left-translating $\xi _y$ by $g$.  It then easily follows that $$ \stem (\xi _y) = \stem _g(\xi _x) \in
\Folow _\beta , $$ so the first case treated above (i.e. $g=1$) applies for $\xi _y$, and we deduce that $\beta \inv \in
\xi _y$, whence $$ g\beta \inv \in g\xi _y=\xi _x.  \endProof

The above result plays a crucial role in understanding the elements of $\Spec $ from the point of view of their stem.
By this we mean that, once the stem of $\xi $ is marked in the Cayley graph of $\F $, and we wish to mark the remaining
group elements in $\xi $, we know from \ref {PropForPaint.v} that we need only worry about elements of the form
$g=\alpha \beta \inv $.  If we are careful to take the reduced form of $g$, then a necessary condition for it to be
marked is that $\alpha $ also be marked, in which case $\alpha $ must be a prefix of the stem.

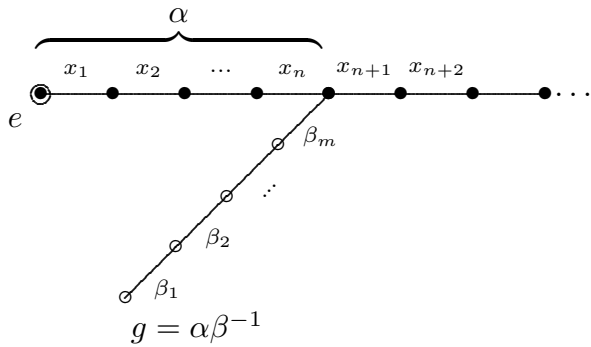
\null \hskip 2cm \beginpicture \setcoordinatesystem units <0.0080truecm, -0.0080truecm> point at 1000 1000 \put {$e$} at
-42 42 \circulararc 360 degrees from 0 16 center at 0 0 \put {$\bullet $} at 0 0 \plot 0 0 120 0 / \put {$\bullet $} at
120 0 \put {$\scriptstyle x_1$} at 60 -40 \plot 120 0 240 0 / \put {$\bullet $} at 240 0 \put {$\scriptstyle x_2$} at
180 -40 \plot 240 0 360 0 / \put {$\bullet $} at 360 0 \put {$\scriptstyle \ldots $} at 300 -40 \plot 360 0 480 0 / \put
{$\bullet $} at 480 0 \put {$\scriptstyle x_n$} at 420 -40 \plot 480 0 600 0 / \put {$\bullet $} at 600 0 \put
{$\scriptstyle x_{n+1}$} at 540 -40 \plot 600 0 720 0 / \put {$\bullet $} at 720 0 \put {$\scriptstyle x_{n+2}$} at 660
-40 \plot 720 0 840 0 / \put {$\bullet $} at 840 0 \put {$\ldots $ } at 895 0 \put {$\overbrace {\hbox to 90pt{\hfill
}}^{\textstyle \alpha }$} at 230 -110 \put {$\circ $} at 141 339 \plot 141 339 225 255 / \put {$\circ $} at 225 255 \put
{$\scriptstyle \beta _1$} at 211 325 \plot 225 255 310 170 / \put {$\circ $} at 310 170 \put {$\scriptstyle \beta _2$}
at 296 240 \plot 310 170 395 85 / \put {$\circ $} at 395 85 \put {$\scriptstyle \nedots $} at 381 156 \plot 395 85 480 0
/ \put {$\circ $} at 480 0 \put {$\scriptstyle \beta _m $} at 466 71 \put {$g=\alpha \beta \inv $} at 261 389
\endpicture \hfill \null

\bigskip We then must decide whether or not to mark $g$ itself, and this is precisely where \ref {Enfiabilidade}
intervenes: in case $\xi $ is some of the $\xi _x$, then we should mark $g$ \"{if and only if} the resulting stem at
$g$, namely $$ \beta _1\beta _2\ldots \beta _mx_{n+1}x_{n+2}\ldots $$ lies in $X$ (which is to say that
$x_{n+1}x_{n+2}\ldots \ $ lies in $\Folow _\beta $).  When $\xi $ is not necessarily a $\xi _x$, then \ref
{Enfiabilidade} does not give a definite answer, except that marking $g$ is \"{forbidden} in case the above infinite
word does not lie in $X$.

At this point it is perhaps useful to discuss an example: it is well known that the even shift is not of finite type\fn
{This will also follow from the analysis we are about to undertake.}, and hence \ref {OntoContinuousSFT} predicts the
existence of elements in $\Spec $ beyond the range of $\Xi $, meaning not of the form $\xi _x$.  In what follows we will
concretely exhibit an example of such anomalous elements.

Let $X$ be the even shift, and for each $n$, consider the infinite word $$ x_n = 1^{^{2n+1}}0^{^\infty } = \underbrace
{1\ldots 1}_{2n+1}0000\ldots $$ Since $\Spec $ is compact, there exists a subsequence, say $\{y_k\}_k = \{x_{n_k}\}_k$,
such that $\{\xi _{y_k}\}_k$, converges to some $\xi \in \Spec $.  Our next goal will be to prove that $\xi $ is not of
the form $\xi _x$, for any $x$ in $X$.

The fact that $\{x_n\}_n$, and hence also $\{y_k\}_k$ converges to the infinite word $$ 1^\infty =1111\ldots $$ relative
to the topology of $X$, does not imply that $\xi _{y_k}$ converges to $\xi _{1^\infty }$, as the correspondence $x\to
\xi _x$ is not known to be continuous (it will soon be evident that it is discontinuous at $1^\infty $).  Nevertheless,
the continuity of the stem \ref {StemContinuous} implies that $$ \stem (\xi ) = \lim _k\stem (\xi _{y_k}) = \lim _k y_k
= 1^\infty .  $$

Therefore, if $\xi =\xi _x$, for some $x$, then $$ 1^\infty = \stem (\xi ) = \stem (\xi _x) = x.  $$ So, to prove that
$\xi $ is not equal to any $\xi _x$ we therefore only need to verify that $\xi \neq \xi _{1^\infty }$.

Observing that $y_k$ is not in the follower set of the finite word $$ \beta = \qt 0, $$ we have by \ref {BinvInXix} that
$\beta \inv $ is not in $\xi _{y_k}$.  By \ref {BooleanContinuous} it follows that $$ [\beta \inv \in \xi ] = \lim _k \
[\beta \inv \in \xi _{y_k}] = 0, $$ so $\beta \inv \notin \xi $.  Nevertheless, $1^\infty $ does belong to the follower
set of $\beta $, hence $\beta \inv \in \xi _{1^\infty }$, again by \ref {BinvInXix}.  This proves that $\xi \neq \xi
_{1^\infty }$, so $\xi $ is not in the range of $\XO $, whence $\XO $ is not onto, and we then deduce from \ref
{OntoContinuousSFT} that $\XO $ is not continuous.  We also recover the well known fact that the even shift is not of
finite type.

This example also illustrates that the implication ``\implica (b)(a)" in \ref {Enfiabilidade} may indeed fail, since the
stem of $\xi $ lies in $\Folow _\beta $, and yet $\beta \inv $ is not in $\xi $.

\section Partial crossed product description of $\CM $

As always we fix a finite alphabet $\Lambda $ and a subshift $X\subseteq \Lambda ^{\bf N}$.  In this section we plan to
prove that $\CM $ is isomorphic to the crossed product of $C(\Spec )$ (also known as $\Core $) by the spectral partial
action of the free group $\F $. In symbols $$ \CM \simeq C(\Spec )\rtimes _\vartheta \F .  \equationmark CMCrossProd $$

We will also show that the associated semi-direct product Fell bundle is amenable \cite {amena}, whence the full and
reduced crossed products coincide.

Regarding the partial action $\tau $ of $\F $ on $\Core $ introduced in \ref {PactOnAlg}, recall that $\vartheta $ was
defined as the partial action on the spectrum of $\Core $ induced by $\tau $.  Thus, moving in the opposite direction,
the partial dynamical system induced by $\vartheta $ on $C(\Spec )$ is equivalent to $\tau $.  In order to prove \ref
{CMCrossProd}, it therefore suffices to prove that $$ \CM \simeq \Core \rtimes _\tau \F .  $$

\state Proposition \label FromCPtoMats There exists a surjective *-homomorphism $$ \varphi :\Core \rtimes _\tau \F \to
\CM , $$ such that $\varphi (a\delta _g)=(a\otimes 1)\prt _g=a\pr _g\otimes \lambda _g$, for all $g$ in $\F $, and every
$a\in D_g$.

\Proof We claim that the pair $(j,\prt )$ is a covariant representation (see \book {Definition}{9.10}) of $\tau $ in
$\CM $, where $$ j:a\in \Core \mapsto a\otimes 1\in \CM .  $$ To see this we pick any $g$ in $\F $, and $a$ in $D_{g\inv
}$, and compute $$ \prt _gj(a)\prt _{g\inv } = (\pr _g\otimes \lambda _g) (a\otimes 1)(\pr _{g\inv }\otimes \lambda
_{g\inv }) = \pr _ga\pr _{g\inv }\otimes 1 \={DefineTau} \tau _g(a)\otimes 1 = j\big (\tau _g(a)\big ).  $$

This shows that indeed $(j,\prt )$ is a covariant representation, so the existence of $\varphi $ follows from \book
{Proposition}{13.1}.  Given $\alpha $ in $\Lambda ^*$, notice that $e_\alpha \in D_\alpha $, so $$ \varphi (e_\alpha
\delta _\alpha ) = e_\alpha \pr _\alpha \otimes \lambda _\alpha = \pr _\alpha \otimes \lambda _\alpha = \prt _\alpha .
$$ This shows that every $\prt _\alpha $ lies in the range of $\varphi $, whence $\varphi $ is onto.  \endProof

In order to find a map in the opposite direction, let us consider the standard conditional expectation $E$ from the
algebra of all bounded operators on $\ell ^2(X)\otimes \ell ^2(\F )$ onto the subalgebra of all diagonal operators
relative to the standard orthonormal basis.  Thus, if $t$ is any bounded operator on $\ell ^2(X)\otimes \ell ^2(\F )$,
then $E(t)$ is the operator whose off diagonal entries are zero, and whose diagonal entries are the same as those of
$t$.

\state Lemma \label CondexpOnMyElts Let $t=\et _{h_1}\et _{h_2}\ldots \et _{h_n}\prt _g$, where $h_1,h_2,\ldots
,h_n,\,g\in \F $.  Then $$ E(t) = \left \{\matrix { t, & \hbox {if } g=1,\hfill \cr \pilar {12pt} 0, & \hbox
{otherwise.}  }\right .  $$

\Proof We have $$ E(t) = E(\et _{h_1}\et _{h_2}\ldots \et _{h_n}\prt _g) = \et _{h_1}\et _{h_2}\ldots \et _{h_n}E(\prt
_g), $$ because the $\et _h$'s are diagonal operators by \ref {ETTensorOne}. When $g=1$ we have that $\prt _g=1$, so it
is easy to see that $E(t)=t$.  When $g\neq 1$, notice that for all $x\in X$, and $h\in \F $, $$ \langle \prt _g(\delta
_x\otimes \delta _h), \delta _x\otimes \delta _h\rangle = \langle \pr _g(\delta _x),\delta _x\rangle \langle \delta
_{gh},\delta _h\rangle =0, $$ since $gh\neq h$.  Thus all diagonal entries for $\prt _g$ vanish, whence $E(\prt _g)=0$,
concluding the proof.  \endProof

Observe that the elements of the form $\et _{h_1}\et _{h_2}\ldots \et _{h_n}\prt _g$ span a *-subalgebra of $\CM $ by
\book {Proposition}{9.8}.  It is therefore a dense subalgebra, and \ref {CondexpOnMyElts} then implies that $\CM $ is
invariant under $E$.  Again by \ref {CondexpOnMyElts} one sees that the $$ E(\CM )\subseteq \Core \otimes 1, $$ so we
may see $E$ as a conditional expectation from $\CM $ to $\Core \otimes 1$.

In our next result we will refer to the \"{reduced partial crossed product} of $\Core $ by $\F $ under $\tau $ (see
\book {Definition}{17.10}), which we will denote by $\Core \rtimes ^{\hbox {\sixrm red}}_\tau \F $.

\state Proposition \label FromMatstoCP There exists a *-homomorphism $$ \psi :\CM \to \Core \rtimes ^{\hbox {\sixrm
red}}_\tau \F , $$ such that $ \psi (\prt _\alpha ) = e_\alpha \delta _\alpha , $ for all $\alpha $ in $\Lambda ^*$.

\Proof For each $g$ in $\F $, let $B_g$ be the closed linear subspace of\/ $\CM $ given by $$ B_g = (D_g\otimes 1)\prt
_g = (D_g\pr _g \otimes \lambda _g).  $$

We claim that $B_gB_h\subseteq B_{gh}$, for all $g,h\in \F $.  In fact, given $a\in D_g$, and $b\in D_h$, we have $$
(a\pr _g \otimes \lambda _g) (b\pr _h \otimes \lambda _h) = a\pr _g b\pr _h \otimes \lambda _{gh}, \equationmark
TensorGuy $$ while $$ a\pr _g b\pr _h= e_ga\pr _g be_{g\inv } \pr _h= \pr _g\pr _{g\inv }a\pr _gb\pr _{g\inv } \pr _g\pr
_h\explicaRef {PDSFB}{9.1.iii}{=} \pr _g\tau _{g\inv }(a)b\pr _{g\inv } \pr _{gh} \$= \tau _g(\tau _{g\inv }(a)b\big )
\pr _{gh} \in \tau _g(D_{g\inv }\cap D_h) \pr _{gh} \subseteq D_{gh}\pr _{gh}.  $$ Therefore the element described in
\ref {TensorGuy} lies in $B_{gh}$, proving the claim.  We leave it for the reader to prove that $(B_g)^*=B_{g\inv }$, as
well as that $\sum _{g\in \F }B_g$ is a dense subspace of $\CM $, which, when combined with the conditional expectation
$E$, provided above, verifies all of the assumptions of \book {Theorem}{19.1}, which in turn provides for the desired
map $\psi $.  \endProof

Combining the above with \ref {FromCPtoMats} we arrive at the main result of this section:

\state Theorem \label AlgebrasSame Let $\Lambda $ be a finite alphabet, and let $X\subseteq \Lambda ^{\bf N}$ be a
subshift.  Then: \izitem \zitem The semi-direct product bundle \book {Proposition}{16.6} corresponding to the spectral
partial action satisfies the approximation property \book {Definition}{20.4}, and hence is amenable \book
{Definition}{20.1}.  \zitem The Carlsen-Matsumoto C*-algebra $\CM $ is naturally isomorphic to both the full and the
reduced crossed product of $C(\Spec )$ by the free group $\F (\Lambda )$ under the spectral partial action.  In symbols
$$ \def \quad { \ } \matrix {\CM & \simeq & C(\Spec )\rtimes _\vartheta \F \hfill \cr \pilar {18pt} & \simeq & C(\Spec
)\rtimes ^{\hbox {\sixrm red}}_\vartheta \F .}  $$

\Proof The first point follows from \cite [Theorems 4.1 \& 6.3]{ortho} (see also \book {Theorem}{20.13}).

As for (ii), regarding the maps $\varphi $ and $\psi $ provided by \ref {FromCPtoMats} and \ref {FromMatstoCP}, it is
easy to see that the composition $\psi \circ \varphi $ is the regular representation \book {Definition}{17.6} of
$C(\Spec )\rtimes _\vartheta \F $, which is an isomorphism by (i).  In particular $\varphi $ is one-to-one, but since we
already saw that it is onto in \ref {FromCPtoMats}, we deduce that it is an isomorphism.  If both $\varphi $ and the
composition $\psi \circ \varphi $ are isomorphisms, then so is $\psi $, whence $\CM $ is isomorphic to the reduced
crossed product as well.  \endProof

In the introduction of \cite {CarlsenCuntzPim} Carlsen argues that, contrary to $\CM $, Matsumoto's algebra $\MatsAlg $
does not have good universal properties and hence they should be considered as the class of \"{reduced C*-algebras
associated with subshifts}.  Interpreting the term \"{reduced} as one usually does when speaking of crossed products,
the above result seems to indicate that $\CM $ should be seen as both full and reduced C*-algebras and that $\MatsAlg $
is just an epimorphic image of $\CM $.

\section Comparison with Carlsen's description of $\CM $

As always we fix a finite alphabet $\Lambda $ and a subshift $X\subseteq \Lambda ^{\bf N}$.

In this section we will prove the fact, already hinted at after \ref {DefineCM}, that $\CM $ is isomorphic to the
algebra introduced by Carlsen in \cite [Definition 5.1]{CarlsenCuntzPim}.

We have already observed that the $e_g$ are diagonal operators relative to the canonical orthonormal basis of $\ell
^2(X)$.  On the other hand, the algebra of all diagonal operators is clearly isomorphic to $\ell ^\infty (X)$, so we may
see $\Core $ as a subalgebra of $\ell ^\infty (X)$.

\state Proposition For all finite words $\alpha $ and $\beta $ in $\Lambda ^*$, let $$ C(\beta ,\alpha ) = \{\alpha y\in
X: y\in \Folow _\alpha \cap \Folow _\beta \}.  $$ Then $\Core $ coincides with the closed *-subalgebra of $\ell ^\infty
(X)$ generated by the characteristic functions $1_{C(\beta ,\alpha )}$, as $\alpha $ and $\beta $ range in $\Lambda ^*$.

\Proof In case $\alpha $ and $\beta $ are such that $|\alpha \beta \inv | = |\alpha | +|\beta \inv |$, that is, when
$g:=\alpha \beta \inv $ is in reduced form, we have seen in \ref {DescrPAction.ii} that $X_g=C(\beta ,\alpha )$.  Since
$e_g$ is the final projection of $\pr _g$, whose range is $\ell ^2(X_g)$ by \ref {FinalSpace}, we have that $e_g$
coincides with $1_{C(\beta ,\alpha )}$ (up to the above identification of diagonal operators and bounded functions).
For general $\alpha $ and $\beta $, let $\gamma $ be the longest common suffix of $\alpha $ and $\beta $, so we may find
$\alpha '$ and $\beta '$ such that $$ \alpha =\alpha '\gamma ,\quad \beta =\beta '\gamma , $$ and $\alpha '{\beta '}\inv
$ is the reduced form for the group element $g:=\alpha \beta \inv $.  We leave it for the reader to check that $$
C(\beta ,\alpha ) = C(\beta ',\alpha ') \cap C(\emptyword ,\alpha ), $$ from where it follows that $$ 1_{C(\beta ,\alpha
)} = 1_{C(\beta ',\alpha ')} 1_{C(\emptyword ,\alpha )} = e_ge_\alpha = e_{\alpha \beta \inv }e_\alpha .  \equationmark
TwoSetsCs $$ We therefore conclude that the algebra generated by all of the $1_{C(\beta ,\alpha )}$ is the same as the
algebra generated by all of the $e_g$.  \endProof

As a consequence we see that $\Core $ is the same as the algebra denoted $\tilde \Core $ studied in \cite [Definition
4.1]{CarlsenCuntzPim}, or the algebra denoted $\Core $ studied in \cite [Lemma 7]{CarlsenSilvestrov}, and which also
appears in many other papers dealing with C*-algebras associated with subshifts (see the references given in the
introduction).

\state Theorem For every subshift $X$, denote by $\CM '$ the C*-algebra introduced by Carlsen in \cite [Definition
5.1]{CarlsenCuntzPim}.  Then there is an isomorphism $$ \varphi :\CM '\to \CM , $$ such that $\varphi (S_\alpha )=\prt
_\alpha $, for all $\alpha $ in $\Lambda ^*$, where $S_\alpha $ is the partial isometry given in \cite [Definition
5.3]{CarlsenCuntzPim}.

\Proof For each $\alpha \in \Lambda ^*$, let $$ T_\alpha =\prt _\alpha =\pr _\alpha \otimes \lambda _\alpha , $$ and
notice that for any given $\alpha ,\beta \in \Lambda ^*$, we have $$ T_\alpha T_\beta ^*T_\beta T_\alpha ^* = \prt
_\alpha \et _{\beta \inv }\prt _{\alpha \inv } \explicaRef {PDSFB}{9.8.iii}{=} \et _{\alpha \beta \inv }\prt _\alpha
\prt _{\alpha \inv } \$= \et _{\alpha \beta \inv }\et _\alpha = (e_{\alpha \beta \inv }e_\alpha )\otimes 1 \={TwoSetsCs}
1_{C(\beta ,\alpha )}\otimes 1.  $$

Using the universal property of $\CM '$ \cite [Remark 7.3]{CarlsenCuntzPim}, we see that there exists a *-homomorphism
$\varphi :\CM '\to \CM $ sending each $S_\alpha $ to $T_\alpha $, and which is therefore necessarily onto.

In order to complete the proof it now suffices to prove that $\varphi $ is injective, and for this we will employ \cite
[Theorem 13]{CarlsenSilvestrov}, which demands that we build a suitable action of the circle group $\T $ on $\CM $.
Consider the unique group homomorphism $$ \varepsilon :\F \to {\bf Z}, $$ such that $\varepsilon (a)=1$, for every $a\in
\Lambda $, so that in particular $\varepsilon (\alpha )=|\alpha |$, if $\alpha $ is a finite word.  Moreover, for each
$z$ in $\T $, let $V_z$ be the unitary operator defined on the canonical orthonormal basis of $\ell ^2(X)\otimes \ell
^2(\F )$ by $$ V_z(\delta _x\otimes \delta _g) = z^{\varepsilon (g)}\delta _x\otimes \delta _g, \for x\in X, \for g\in
\F .  $$ We then claim that $$ V_zT_\alpha V_{z\inv } = z^{|\alpha |}T_\alpha , \for z\in \T , \for \alpha \in \Lambda
^*.  \equationmark ConjugaT $$ To prove it we compute on a general element $\delta _x\otimes \delta _g$ of the canonical
orthonormal basis: $$ V_zT_\alpha V_{z\inv }(\delta _x\otimes \delta _g) = z^{-\varepsilon (g)}V_z(\pr _\alpha \otimes
\lambda _\alpha )(\delta _x\otimes \delta _g) = z^{-\varepsilon (g)}V_z\big (\pr _\alpha (\delta _x)\otimes \delta
_{\alpha g}\big ) \$= z^{-\varepsilon (g)+\varepsilon (\alpha g)} \big (\pr _\alpha (\delta _x)\otimes \delta _{\alpha
g}\big ) = z^{\varepsilon (\alpha )} T_\alpha (\delta _x\otimes \delta _g).  $$

This proves \ref {ConjugaT}, so it follows that the formula $$ \gamma _z(a) = V_z aV_{z\inv }, \for x\in \T , \for a\in
\CM , $$ defines a strongly continuous action of $\T $ on $\CM $ satisfying point {\it (2)} of \cite [Theorem
13]{CarlsenSilvestrov}, whence $\varphi $ is injective.  \endProof

\section The topology on the spectrum

The topology on $\Spec $, being induced from the product topology of $2^\F $, admits a basis formed by the open sets $$
\def \quad {\kern 5pt} U_{g_1,g_2,\ldots ,g_n;\,h_1,h_2,\ldots ,h_m} = \left \{\matrix {\eta \in \Spec :& g_1\in \eta ,&
g_2\in \eta ,& \ldots \ , & g_n\in \eta \cr \pilar {11pt} &h_1\notin \eta ,& h_2\notin \eta ,& \ldots \ ,& h_m\notin
\eta }\right \}, \equationmark BasisProdTop $$ where $g_1,g_2,\ldots ,g_n;\ h_1,h_2,\ldots ,h_m$ range in $\F $.  Given
the special nature of elements of $\Spec $, we may restrict to sets of a somewhat special nature, as follows.

\state Proposition \label BaseTopology Given $\alpha ,\ \beta _1,\beta _2,\ldots ,\beta _n,\ \gamma _1,\gamma _2,\ldots
,\gamma _m\in \F _+$, consider the subset of\/ $\Spec $ given by $$ \def \quad {\kern 5pt} V_{\alpha ;\,\beta _1,\beta
_2,\ldots ,\beta _n;\,\gamma _1,\gamma _2,\ldots ,\gamma _m} = \left \{\matrix {\eta \in \Spec :& \alpha &\in &\eta
,\hfill \cr \pilar {11pt} & \alpha \beta _i\inv &\in &\eta , & \hbox { for } i=1,\ldots ,n, \cr \pilar {11pt} & \alpha
\gamma _j\inv &\notin &\eta , &\hbox { for } j=1,\ldots ,m}\right \}, $$ Then the collection consisting of all sets of
the above form is a basis for the topology of $\Spec $.

\Proof Let $A$ be an open subset of $\Spec $, and let $\xi \in A$.  It suffices to prove that there is an open set $V$
of the form described in the statement with $\xi \in V\subseteq A$.

By the definition of the product topology on $2^\F $ there are $g_1,g_2,\ldots ,g_n;\,h_1,h_2,\ldots ,$ $h_m\in \F $, as
in \ref {BasisProdTop}, such that $$ \xi \in U_{g_1,g_2,\ldots ,g_n;\,h_1,h_2,\ldots ,h_m}\subseteq A.  \equationmark
BigNBD $$

Since the elements of $\Spec $ only contain group elements of the form $\alpha \beta \inv $, with $\alpha ,\beta \in \F
_+$ by \ref {PropForPaint.v}, we may ignore the $g_i$ and the $h_j$ which are not of this form without affecting \ref
{BigNBD}.  We may therefore assume that $$ g_i=\alpha _i\beta _i\inv , \and h_j=\mu _j\gamma _j\inv , $$ in reduced
form, with $\alpha _i,\beta _i,\mu _j,\gamma _j\in \F _+$.  Let $\alpha $ be any prefix of the stem of $\xi $, long
enough so that $$ |\alpha |\geq \max \big \{|\alpha _1|, |\alpha _2|, \ldots , |\alpha _n|, |\mu _1|, |\mu _2|, \ldots ,
|\mu _m|\big \}, \equationmark LongAlpha $$ and observe that since $\alpha \in \xi $, then $$ \def \quad {\kern 5pt}
\matrix {\xi &\in & U_{\alpha ,g_1,g_2,\ldots ,g_n;\,h_1,h_2,\ldots ,h_m} \subseteq \hfill \cr \pilar {12pt} &\subseteq
& U_{\phantom {\alpha ,}g_1,g_2,\ldots ,g_n;\,h_1,h_2,\ldots ,h_m} \subseteq A.}  $$

The reader is asked to compare the display above with \ref {BigNBD}, paying special attention to the important detail
that $\alpha $ was inserted ahead of the $g_i$'s in the subscripts of the first occurrence of $U$ above.  The proof will
consist in showing that this occurrence of $U$ coincides with the set displayed in the statement for suitable choices of
$\beta _i$ and $\gamma _j$.

Notice that, by assumption the $g_i\in \xi $, so {\convexity } implies that $\alpha _i\in \xi $, and then the $\alpha
_i$ are necessarily prefixes of the stem of $\xi $ by \ref {PropStem}.  Consequently the $\alpha _i$ are also prefixes
of the long $\alpha $ chosen above, and we may find suitable elements $\delta _i$ in $\F _+$, such that $$ \alpha
=\alpha _i\delta _i, $$ for all $i$.  Letting $\beta _i'=\beta _i\delta _i$, observe that $$ \alpha {\beta _i'}\inv =
\alpha \delta _i\inv \beta _i\inv = \alpha _i\beta _i\inv = g_i, \equationmark TrocaBeta $$ so, upon replacing each
$\beta _i$ by $\beta _i'$, we may suppose that $g_i=\alpha \beta _i\inv $.  It should however be noticed that this
presentation of $g_i$ is no longer in reduced form.

Next we should treat the $h_j$.  Firstly, let us consider those $h_j$ whose corresponding $\mu _j$ is \"{not} a prefix
of the stem of $\xi $.  We then claim that, for any $\eta $ in $\Spec $, one has that $$ \alpha \in \eta \imply \mu
_j\notin \eta .  $$ Otherwise, if both $\alpha $ and $\mu _j$ lie in $\eta $, then both would be prefixes of the stem of
$\eta $, in which case $\mu _j$ would be a prefix of $\alpha $, since the former it is shorter than the latter by \ref
{LongAlpha}.  This would entail that $\mu _j$ is a prefix of the stem of $\xi $, contradicting our assumptions, and
hence proving our claim.  By {\convexity } we have that $$ \mu _j\notin \eta \imply \mu _j\gamma _j\inv \notin \eta , $$
so a combination of these two implications yields $$ \alpha \in \eta \imply h_j\inv \notin \eta .  $$

We then see that $h_j$ may be deleted from the list of subscripts of our $$ U_{\alpha ,g_1,g_2,\ldots
,g_n;\,h_1,h_2,\ldots ,h_m}, \equationmark HereIsU $$ since the condition ``$\alpha \in \eta $" in the first line of
\ref {BasisProdTop} already gives ``$h_j\inv \notin \eta $", in the second.

After deleting such $h_j$'s, we may assume that $\mu _j$ is a prefix of the stem of $\xi $, for every $j$, and, again by
\ref {LongAlpha}, $\mu _j$ is necessarily a prefix of $\alpha $.  Arguing as in \ref {TrocaBeta} we may then \"{stretch}
each $\mu _j$ and $\gamma _j$ by the same amount, and hence assume that the $\mu _j$ all coincide with $\alpha $, so
that $h_j=\alpha \gamma _j\inv $.

The description of the set in \ref {HereIsU} is then identical to the description of the set displayed in the statement,
so we have concluded the proof.  \endProof

Up to the statement of the above result, whenever we considered an expression of the form ``$\alpha \beta \inv $", this
was supposed to be in reduced form.  However the reader should be warned that this is no longer the case, especially
after \ref {TrocaBeta}, were we deliberately gave up on reduced forms in exchange for working with a single $\alpha $.

There is a further simplification which may be bestowed upon the general form of the open sets described in \ref
{BaseTopology}, provided we are concerned with neighborhoods of points in the range of $\XO $, namely we may do away
with the $\gamma _j$.  Making this idea precise is our next goal.

\state Proposition \label VerySimpleNBDs Given $x$ in $X$, let $\xi _x=\XO (x)$.  Then the collection of all sets of the
form $$ \def \quad {\kern 5pt} V_{\alpha ;\,\beta _1,\beta _2,\ldots ,\beta _n} = \left \{\matrix {\eta \in \Spec :&
\alpha &\in &\eta ,\hfill \cr \pilar {11pt} & \alpha \beta _i\inv &\in &\eta , & \hbox { for } i=1,\ldots ,n}\right \},
$$ for $\alpha ,\ \beta _1,\beta _2,\ldots ,\beta _n\in \F _+$, which moreover contains $\xi _x$, forms a neighborhood
base for $\xi _x$.

\Proof Needless to say, the above sets are special cases of the sets in \ref {BaseTopology}, corresponding to taking
$m=0$, meaning that the conditions ``$\alpha \gamma _j\inv \notin \eta $" are now absent.

In order to prove the statement, we must show that for every open set $U$ containing $\xi _x$, there are $\alpha ,\
\beta _1,\beta _2,\ldots ,\beta _n\in \F _+$ such that $$ \xi _x\in V_{\alpha ;\,\beta _1,\beta _2,\ldots ,\beta _n}
\subseteq U.  $$ Using \ref {BaseTopology} we may clearly suppose that $$ U=V_{\alpha ;\,\beta _1,\beta _2,\ldots ,\beta
_n;\,\gamma _1,\gamma _2,\ldots ,\gamma _m}, $$ for suitable $\alpha ,\ \beta _1,\beta _2,\ldots ,\beta _n,\ \gamma
_1,\gamma _2,\ldots ,\gamma _m\in \F _+$.  Observing that $\xi _x$ lies in $U$, we have that $\alpha \in \xi _x$, while
$\alpha \gamma _j\inv \notin \xi _x$, for every $j=1,\ldots ,m$.

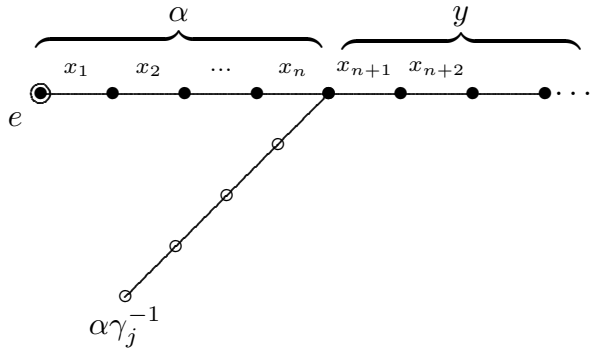
\null \hskip 2cm \beginpicture \setcoordinatesystem units <0.0080truecm, -0.0080truecm> point at 1000 1000 \put {$e$} at
-42 42 \circulararc 360 degrees from 0 16 center at 0 0 \put {$\bullet $} at 0 0 \plot 0 0 120 0 / \put {$\bullet $} at
120 0 \put {$\scriptstyle x_1$} at 60 -40 \plot 120 0 240 0 / \put {$\bullet $} at 240 0 \put {$\scriptstyle x_2$} at
180 -40 \plot 240 0 360 0 / \put {$\bullet $} at 360 0 \put {$\scriptstyle \ldots $} at 300 -40 \plot 360 0 480 0 / \put
{$\bullet $} at 480 0 \put {$\scriptstyle x_n$} at 420 -40 \plot 480 0 600 0 / \put {$\bullet $} at 600 0 \put
{$\scriptstyle x_{n+1}$} at 540 -40 \plot 600 0 720 0 / \put {$\bullet $} at 720 0 \put {$\scriptstyle x_{n+2}$} at 660
-40 \plot 720 0 840 0 / \put {$\bullet $} at 840 0 \put {$\ldots $ } at 895 0 \put {$\overbrace {\hbox to 90pt{\hfill
}}^{\textstyle \alpha }$} at 230 -110 \put {$\overbrace {\hbox to 75pt{\hfill }}^{\textstyle y}$} at 700 -110 \put
{$\circ $} at 141 339 \plot 141 339 225 255 / \put {$\circ $} at 225 255 \put {$\scriptstyle \null $} at 211 325 \plot
225 255 310 170 / \put {$\circ $} at 310 170 \put {$\scriptstyle \null $} at 296 240 \plot 310 170 395 85 / \put {$\circ
$} at 395 85 \put {$\scriptstyle \null $} at 381 156 \plot 395 85 480 0 / \put {$\circ $} at 480 0 \put {$\scriptstyle
\null $} at 466 71 \put {$\alpha \gamma _j\inv $} at 141 389 \endpicture \hfill \null

\bigskip

It follows that $\alpha $ is a prefix of the stem of $\xi _x$, also known as $x$, so we may write $x=\alpha y$ for some
infinite word $y$.  The stem of $\xi _x$ at $\alpha $ is therefore $y$ and the fact, noted above, that $\alpha \gamma
_j\inv \notin \xi _x$, together with \ref {Enfiabilidade} leads to the conclusion that $y$ is not in the follower set of
$\gamma _j$, which is to say that $$ \gamma _jy \notin X.  $$ By \ref {AllSSForbWords} we then deduce that $\gamma _jy$
has some prefix which is a forbidden word, namely a finite word not in the language $\Lang $.  By increasing its length
we may suppose that this forbidden prefix is of the form $$ \gamma _j\delta _j $$ where $\delta _j$ is a prefix of $y$.
Denote by $\delta $ the longer among the $\delta _j$, and set $$ \alpha '=\alpha \delta , \and \beta _i'=\beta _i\delta
, $$ for every $i=1,\ldots ,n.$ It follows that $$ \alpha '{\beta _i'}\inv = \alpha \beta _i\inv , $$ and the proof will
be concluded once we show that $$ \xi _x\in V_{\alpha ';\,\beta _1',\beta _2',\ldots ,\beta _n'} \subseteq V_{\alpha
;\,\beta _1,\beta _2,\ldots ,\beta _n;\,\gamma _1,\gamma _2,\ldots ,\gamma _m}.  \equationmark SimplerNBD $$

We leave the easy ``$\in $" for the reader to check and concentrate on the ``$\subseteq $".  We thus pick $\eta $ in the
set appearing in the left-hand-side above, and we note that \iaitem \aitem $\alpha '\in \eta $, \aitem $\alpha '{\beta
_j'}\inv \in \eta $, \medskip \noindent for every $j=1,\ldots ,m$.  By (a) and {\convexity } we have that $\alpha \in
\eta $, and clearly $$ \alpha \beta _i\inv =\alpha '{\beta _j'}\inv \in \eta .  $$

Thus, in order to prove that $\eta $ lies in the set appearing in the right-hand-side of \ref {SimplerNBD}, we must only
check that $\alpha \gamma _j\inv \notin \eta $.  For this, observe that $$ \alpha \gamma _j\inv = \alpha \delta (\gamma
_j\delta )\inv = \alpha '(\gamma _j\delta )\inv .  $$ Assuming by contradiction that this element belongs to $\eta $, we
have by (a) and \ref {Enfiabilidade} that $$ \stem _{\alpha '}(\eta )\in \Folow _{\gamma _j\delta }.  $$

Letting $z= \stem _{\alpha '}(\eta )$, this means that $\gamma _j\delta z\in X$, but $\gamma _j\delta $ admits the
forbidden word $\gamma _j\delta _j$ as a prefix, a contradiction.  This proves our claim that $\alpha \gamma _j\inv
\notin \eta $, and hence that $$ \eta \in V_{\alpha ;\,\beta _1,\beta _2,\ldots ,\beta _n;\,\gamma _1,\gamma _2,\ldots
,\gamma _m}, $$ showing \ref {SimplerNBD}.  \endProof

Speaking of a neighborhood of the form $V_{\alpha ;\,\beta _1,\beta _2,\ldots ,\beta _n}$, as above, notice that for
every $\eta $ in this set, one automatically has that $$ \alpha \beta \inv \in \eta , $$ for $\beta =\emptyword $, as
well as for $\beta =\alpha $, regardless of whether or not the words $\emptyword $ and $\alpha $ are among the $\beta
_i$'s.  Therefore, should one so wish, these two words may be added to the $\beta _i$'s without altering the resulting
neighborhood.  That is $$ V_{\alpha ;\,\beta _1,\beta _2,\ldots ,\beta _n} = V_{\alpha ;\,\emptyword ,\alpha ,\beta
_1,\beta _2,\ldots ,\beta _n}.  $$

\section Topological freeness

\label TopFreeSect Recall from \cite [Definition 2.1]{ELQ} (see also \cite {ArchSpiel} and \cite [Section 29]{PDSFB})
that a topological partial action $$ \rho = \big (\{Y_g\}_{g\in G}, \{\rho _g\}_{g\in G}\big ) $$ of a group $G$ on a
space $Y$ is said to be \"{topologically free} if, for every $g\neq 1$, the set of fixed points for $\rho _g$, namely $$
\Fix _g= \{x\in Y_{g\inv }: \rho _g(x)=x\}, \equationmark IntroduceFix $$ has empty interior.

In this section we will give necessary and sufficient conditions for a general subshift to have a topologically free
spectral partial action.  However, let us begin by reviewing the well known Markov case.

Observing that Markov subshifts are of finite type, we see that the spectral partial action is equivalent to the
standard partial action by \ref {EquivPact}.

The following result, inspired by condition (I) of \cite {CuntzKrieger}, characterizes topological freeness in terms of
circuits and exits.  It was first proved for row-finite matrices in \cite [Lemma 3.4]{KPR} in the language of groupoids.
A generalization for infinite matrices was given in \cite [Proposition 12.2]{infinoa}.

\state Proposition \label OldTofFree Let $X$ be a Markov subshift.  Then the standard partial action $\bigTheta $ of\/
$\F $ on $X$ (which, by \ref {EquivPact} is equivalent to the spectral partial action) is topologically free if and only
if every circuit has an exit.

For subshifts not of finite type, the situation is a lot more delicate.  Even if all circuits have strong exits (see
\ref {DefStrongExit}), the spectral partial action may fail to be topologically free.

To see this, consider the even shift.  In section \ref {SubshSect} we have proven that all of its circuits have strong
exits, and yet its associated spectral partial action is not topologically free for the following reason: consider the
infinite word $$ 1^\infty =11111\ldots , $$ which is fixed by $\theta _g$, where $g$ is the element\fn {Warning: this is
one of the generators of $\F $ and not the unit group element!}  of $\F $ corresponding to the word $\qt {1}$.  By \ref
{Equivariance} we have that $\XO $ is equivariant, so $\xi _{1^\infty }$ is fixed by $\vartheta _g$.

We will show that $\vartheta $ is not topologically free by showing that $\xi _{1^\infty }$ is an isolated point in
$\Spec $ (even though $1^\infty $ is not isolated in $X$).  In fact, let $$ \beta _1=01,\and \beta _2=011, $$ and
consider the open set $ V_{\emptyword ; \beta _1, \beta _2} $ described in \ref {VerySimpleNBDs}.  Notice that for any
$x$ in $X$ one has that $$ \xi _x \in V_{\emptyword ; \beta _1, \beta _2} \iff \beta _1\inv \kern -3pt,\kern 3pt \beta
_2\inv \in \xi _x \explain {Enfiabilidade}\iff x\in \Folow _{\beta _1}\cap \Folow _{\beta _2}.  $$

We leave it for the reader to check that the rules of the even shift imply that $ \Folow _{\beta _1}\cap \Folow _{\beta
_2} = \{1^\infty \}, $ so the only $\xi _x$ in $V_{\emptyword ; \beta _1, \beta _2}$ is $\xi _{1^\infty }$.  Since the
set of all $\xi _x$ is dense in $\Spec $ by \ref {XIsDense}, and since $V_{\emptyword ; \beta _1, \beta _2}$ is open, a
simple exercise in Topology gives that $$ V_{\emptyword ; \beta _1, \beta _2}=\{\xi _{1^\infty }\}.  \equationmark
IsolatedPoint $$ The above is then an open set of fixed points, whence $\vartheta $ is not topologically free.

We have therefore proven:

\state Proposition \label EvenNotTopFree The spectral partial action associated to the even shift is not topologically
free.

One might wonder if the standard (as opposed to spectral) partial action $\theta $ for the even shift is topologically
free.  Although we believe this is not a well posed question, since the standard partial is not topological (the $X_g$
are not all open), one might decide to ignore this and insist in checking whether the interior of any set of fixed
points is empty.  In this case the answer is easily seen to be positive, so in this sense the standard partial of the
even shift is topologically free.

Before we give the appropriate characterization of topological freeness for the spectral partial action of general
subshifts, let us understand their fixed points a little better.  The following result is entirely similar to known
results for Markov subshifts, but we give a full proof, which takes no more than a few lines, for the convenience of the
reader:

\state Proposition \label FixPoints Let $\Lambda $ be a finite alphabet, let $X\subseteq \Lambda ^{\bf N}$ be a
subshift, and denote by $\theta $ the standard partial action of\/ $\F $ on $X$.  Given $g\in \F \setminus \{1\}$, let
$x\in X$ be a fixed point for $\theta _g$.  Then \izitem \zitem $g$ admits a decomposition in reduced form as $\nu
\alpha ^{\pm 1}\nu \inv $, with $\alpha ,\nu \in \F _+$, \zitem $x=\nu \alpha ^\infty $, whence $\alpha $ is a circuit,
and $x$ is the unique fixed point for $g$.

\Proof Since $x$ lies in the domain of $\theta _g$, that domain is nonempty, whence $g\in \PPinv $ by \ref
{DescrPAction.i}.  We may therefore write $g$ in reduced form as $\mu \nu \inv $, with $\mu ,\nu \in \F _+$.  Since
$x\in X_{g\inv }$, we have by \ref {DescrPAction.ii} that $x=\nu y$, for some $y\in \Folow _\nu \cap \Folow _\mu $, and
then $$ \nu y = x = \theta _{\mu \nu \inv }(x) = \mu y.  $$

Therefore either $\nu $ is a prefix of $\mu $, or vice-versa.  We assume without loss of generality that $|\mu |>|\nu |$
(one cannot have $|\mu |=|\nu |$ because $g\neq 1$), and so we may write $\mu =\nu \alpha $, for some $\alpha \in
\Lambda ^*$.  It then follows from the above that $\nu y = \nu \alpha y$, so $y = \alpha y$, whence $y=\alpha ^\infty $.
Summarizing we have $$ g=\mu \nu \inv = \nu \alpha \nu \inv , \and x=\nu y=\nu \alpha ^\infty .  $$ completing the
proof.  \endProof

We may now present the main result of this section.  Contrary to one might expect, this result does not explicitly
mention the existence of exits for circuits, but please see the remarks after the proof below for an interpretation in
terms of exits for circuits.

\state Theorem \label MainTopFree Let $\Lambda $ be a finite alphabet and let $X\subseteq \Lambda ^{\bf N}$ be a
subshift.  Then the following are equivalent: \izitem \zitem the spectral partial action $\vartheta $ associated to $X$
is topologically free, \zitem for every $\beta _1,\beta _2,\ldots ,\beta _n$ in $\Lambda ^*$, and for every circuit
$\gamma $ such that $$ \gamma ^\infty \in \medcap _{i=1}^n\Folow _{\beta _i}, $$ one has that $\bigcap _{i=1}^n\Folow
_{\beta _i}$ contains some element other than $\gamma ^\infty $.

\Proof \ipfimply \pfimply (i)(ii) Pick $\beta _1,\beta _2,\ldots ,\beta _n$ in $\Lambda ^*$, and let $\gamma $ be a
circuit such that $\gamma ^\infty \in \bigcap _{i=1}^n\Folow _{\beta _i}$.

Notice that $\gamma ^\infty $ is a fixed point for $\theta _\gamma $, so $\xi _{\gamma ^\infty }$ is a fixed point for
$\vartheta _\gamma $ by \ref {Equivariance}.  Also observe that $\beta _i\inv \in \xi _{\gamma ^\infty }$, for all $i$,
by \ref {BinvInXix}, so $$ \xi _{\gamma ^\infty } \in V_{\emptyword ;\beta _1,\beta _2,\ldots ,\beta _n}.  $$ Should
$\xi _{\gamma ^\infty }$ be the only element of this set, the singleton $\{\xi _{\gamma ^\infty }\}$ would be an open
set consisting of fixed points for $\vartheta _\gamma $, contradicting (i).  Thus $$ V_{\emptyword ;\beta _1,\beta
_2,\ldots ,\beta _n} \setminus \{\xi _{\gamma ^\infty }\} $$ is not empty, hence it contains some $\xi _x$, for $x$ in
$X$, as $\XO (X)$ is dense in $\Spec $.  Thus necessarily $x\neq \gamma ^\infty $, and since each $\beta _i\inv $
belongs to $\xi _x$, we have that $x\in \Folow _{\beta _i}$, again by \ref {BinvInXix}.  This proves (ii).

\pfimply (ii)(i) Arguing by contradiction, let $g\in \F \setminus \{1\}$, and let $U$ be a nonempty open subset of
$\Omega _{g\inv }$ consisting of fixed points for $\vartheta _g$.  Since $\XO (X)$ is dense in $\Spec $ by \ref
{XIsDense}, there exists some $x$ in $X$ such that $\xi _x\in U$, whence $\xi _x$ is fixed by $\vartheta _g$.

By \ref {Equivariance} we conclude that $x = \stem (\xi _x)$ is fixed by $\theta _g$, so \ref {FixPoints.i} provides
$\alpha ,\nu \in \F _+$, such that $g=\nu \alpha ^{\pm 1}\nu \inv $, in reduced form.  Upon replacing $g$ by $g\inv $ if
necessary (and keeping $U$ unaltered), we may assume without loss of generality that $g=\nu \alpha \nu \inv $.
Observing that $\vartheta $ is semi-saturated, we have that $$ \vartheta _g=\vartheta _\nu \circ \vartheta _\alpha \circ
\vartheta _{\nu \inv }, $$ so in particular $\Omega _{g\inv }\subseteq \Omega _\nu $.  Setting $W=\vartheta _{\nu \inv
}(U)$ we then have that $W\subseteq \Omega _{\alpha \inv }$, and it is clear that $\vartheta _\alpha $ is the identity
on $W$.

The upshot of the above argument is that, if $\vartheta $ is not topologically free, then there exists a finite word
$\alpha $, and a nonempty open set $W$ consisting of fixed points for $\vartheta _\alpha $.  Employing \ref {XIsDense}
once more, there exists some $x$ in $X$ such that $\xi _x\in W$.  Therefore $\xi _x$ is fixed by $\vartheta _\alpha $,
whence $x=\stem (\xi _x)$ is fixed by $\theta _\alpha $, and we deduce from \ref {FixPoints.ii} that $x=\alpha ^\infty
$.

Using the special neighborhood base of $\xi _x$ provided by \ref {VerySimpleNBDs}, we may then find finite words $\mu
,\beta _1,\beta _2,\ldots ,\beta _n\in \Lambda ^*$, such that $$ \xi _x\in V_{\mu ;\,\beta _1,\beta _2,\ldots ,\beta _n}
\subseteq W.  \equationmark FoundBasicOnFP $$ In particular we have that $\mu $ lies in $\xi _x$, which in turn implies
that $\mu $ is a prefix of $x$ by \ref {CharacXix}.  Being a prefix of $x=\alpha ^\infty $ therefore implies that $\mu $
is a prefix of $\alpha ^n$, for some $n$, whence there exists a finite word $\gamma $ such that $\alpha ^n=\mu \gamma $.
Setting $\mu '=\mu \gamma $, and $\beta '_i=\beta _i\gamma $, we have that $$ \xi _x\in V_{\mu ';\,\beta '_1,\beta
'_2,\ldots ,\beta '_n} \subseteq V_{\mu ;\,\beta _1,\beta _2,\ldots ,\beta _n} \subseteq W.  $$ Replacing $\mu $ by $\mu
'$, and each $\beta _i$ by the corresponding $\beta '_i$, we then may assume without loss of generality that \ref
{FoundBasicOnFP} reads $$ \xi _x\in V_{\alpha ^n;\,\beta _1,\beta _2,\ldots ,\beta _n} \subseteq W.  $$

Among other things we then have that both $\alpha ^n$ and $\alpha ^n\beta _i\inv $ lie in $\xi _x$, for all $i$, so \ref
{Enfiabilidade} gives $$ \stem _{\alpha ^n}(\xi _x)\in \Folow _{\beta _i}.  $$

Recalling that $x=\alpha ^\infty $, and staring at the definition of the stem, will make it clear that $\stem _{\alpha
^n}(\xi _x)=\alpha ^\infty $.  So $$ \alpha ^\infty \in \medcap _{i=1}^n\Folow _{\beta _i}.\phantom {\cap \Folow
_{\alpha ^n}.}  $$ Noticing that $\alpha ^\infty \in \Folow _{\alpha ^n}$, we may soup up the above conclusion by
writing $$ \alpha ^\infty \in \medcap _{i=1}^n\Folow _{\beta _i} \cap \Folow _{\alpha ^n}, $$ and then we may use
hypothesis (ii) to produce an infinite word $y\neq \alpha ^\infty $ belonging to $\bigcap _{i=1}^n\Folow _{\beta _i}$,
as well as to $\Folow _{\alpha ^n}$.  The infinite word $$ z=\alpha ^ny $$ thus lies in $X$ and the stem of $\xi _z$ at
$\alpha ^n$ is clearly $y$.  For that reason, and using \ref {Enfiabilidade}, we have that $\alpha ^n\beta _i\inv \in
\xi _z$, for all $i$, whence $$ \xi _z\in V_{\alpha ^n;\,\beta _1,\beta _2,\ldots ,\beta _n}\subseteq W, $$ so $\xi _z$
is a fixed point for $\vartheta _\alpha $, whence $z$ is a fixed point for $\theta _\alpha $.  However the only such
fixed point is $x$, which is manifestly different from $z$.  We have thus reached a contradiction, hence proving that
$\vartheta $ is topologically free.  \endProof

When one compares our last result to the well known characterization of topological freeness for Markov subshifts given
in \ref {OldTofFree}, one might wonder what happened to the role of exits for circuits.  Although property \ref
{MainTopFree.ii} does not explicitly mention exits, it is closely related to that concept.  To see this let $\gamma $ be
a circuit.  Then evidently $\gamma ^\infty \in \Folow _{\gamma ^n}$, for every $n$, so the existence of an element $x$
in $\Folow _{\gamma ^n}$, other than $\gamma ^\infty $, as provided by \ref {MainTopFree.ii}, gives an exit for $\gamma
^n$.  If this indeed holds for every $n$, then $\gamma $ has a strong exit by \ref {OneExit}.

We therefore see that \ref {MainTopFree.ii} implies that all circuits have strong exits. However we should not forget
that the existence of strong exits by itself is not enough to guarantee topological freeness for the spectral partial
action, as the example of the even shift in \ref {EvenExit} above shows.

\section Minimality

Recall from \cite [Definition 2.8]{ELQ} that a topological partial action $$ \rho = \big (\{Y_g\}_{g\in G}, \{\rho
_g\}_{g\in G}\big ) $$ of a group $G$ on a space $Y$ is said to be \"{minimal} if there are no nontrivial $\rho
$-invariant closed subsets.  This is clearly equivalent to saying that for every $y$ in $Y$, the \"{orbit of $y$},
namely the set $$ \Orb (y)= \{\rho _g(y): g\in G,\ Y_{g\inv }\ni y\}, $$ is dense in $Y$.

Minimality is well understood for Markov subshifts as well as for many other partial dynamical systems related to it.
See, for example \cite {KPR} and \cite {KPRR}.

For the purpose of comparison, lets us state the following well known result.

\state Proposition \label ClassicMinimal Let $X$ be a Markov subshift with transition matrix $A$.  Then the standard
partial action $\bigTheta $ of\/ $\F $ on $X$ (which, by \ref {EquivPact} is equivalent to the spectral partial action)
is minimal if and only if, for any vertex $a$, and any infinite path $x$ in $Gr(A)$, there exists a finite path starting
in $a$ and ending in some vertex of $x$.

The property described in the above result is sometimes referred to as \"{cofinality}.  Motivated by this concept it is
natural to consider the following property applicable for general subshifts.

\definition \label Reachability Let $\Lambda $ be a finite alphabet and let $X\subseteq \Lambda ^{\bf N}$ be a subshift.
We shall say that an infinite word $x$ in $X$ may be \"{reached} from a finite word $\beta $ in $\Lang $, when there
exists a finite word $\gamma $, and a prefix $\alpha $ of $x$, such that, upon writing $x=\alpha y$, one has that $\beta
\gamma y\in X$.

\null \hskip 2cm \beginpicture \setcoordinatesystem units <0.0080truecm, -0.0080truecm> point at 1000 1000 \put {$e$} at
-42 42 \circulararc 360 degrees from 0 16 center at 0 0 \put {$\bullet $} at 0 0 \plot 0 0 120 0 / \put {$\bullet $} at
120 0 \put {$\scriptstyle x_1$} at 60 -40 \plot 120 0 240 0 / \put {$\bullet $} at 240 0 \put {$\scriptstyle x_2$} at
180 -40 \plot 240 0 360 0 / \put {$\bullet $} at 360 0 \put {$\scriptstyle \ldots $} at 300 -40 \plot 360 0 480 0 / \put
{$\bullet $} at 480 0 \put {$\scriptstyle x_n$} at 420 -40 \plot 480 0 600 0 / \put {$\bullet $} at 600 0 \put
{$\scriptstyle x_{n+1}$} at 540 -40 \plot 600 0 720 0 / \put {$\bullet $} at 720 0 \put {$\scriptstyle x_{n+2}$} at 660
-40 \plot 720 0 840 0 / \put {$\bullet $} at 840 0 \put {$\ldots $ } at 895 0 \put {$\overbrace {\hbox to 90pt{\hfill
}}^{\textstyle \alpha }$} at 230 -110 \put {$\overbrace {\hbox to 75pt{\hfill }}^{\textstyle y}$} at 700 -110
\circulararc 360 degrees from 480 16 center at 480 0 \put {$\bullet $} at 141 339 \plot 141 339 225 255 / \put {$\bullet
$} at 225 255 \put {$\scriptstyle \gamma _1$} at 211 325 \plot 225 255 310 170 / \put {$\bullet $} at 310 170 \put
{$\scriptstyle \gamma _2$} at 296 240 \plot 310 170 395 85 / \put {$\bullet $} at 395 85 \put {$\scriptstyle \nedots $}
at 381 156 \plot 395 85 480 0 / \put {$\bullet $} at 480 0 \put {$\scriptstyle \gamma _m$} at 466 71 \circulararc 360
degrees from 141 355 center at 141 339 \put {$\bullet $} at -212 410 \plot -212 410 -95 386 / \put {$\bullet $} at -95
386 \put {$\scriptstyle \beta _1$} at -148 427 \plot -95 386 23 363 / \put {$\bullet $} at 23 363 \put {$\scriptstyle
\beta _2$} at -30 404 \plot 23 363 141 339 / \put {$\bullet $} at 141 339 \put {$\scriptstyle \ldots $} at 88 381
\endpicture \hfill \null

\centerline {\eightrm Diagram \advseqnumbering (\current )} \label PombaSimplesDiagrama

\bigskip \noindent We shall moreover say that $X$ is \"{cofinal} if, for every $\beta \in \Lang $, and every $x$ in $X$,
one has that $x$ may be reached from $\beta $.

One could think of $\gamma $ as being the \"{bridge} allowing one to travel from $\beta $ in order to reach $x$ at
$\alpha $.  A quicker way to express \ref {Reachability} is to say that $$ \exists \gamma \in \Lambda ^*,\ \exists n\in
{\bf N}: \ S^n(x)\in \Folow _{\beta \gamma }.  $$

Substituting the above notion of cofinality for the classical notion, one might be tempted to conjecture a
generalization of Proposition \ref {ClassicMinimal} above to arbitrary subshifts, but once more the even shift stands as
a counter-example.  It is easy to see that the even shift satisfies the property just mentioned, because, given any $x$
and $\beta $, as in the above diagram, one could simply take $\alpha =\emptyword $, and then either $$ \gamma
=\emptyword ,\quad \hbox {or}\quad \gamma =\qt {1} \equationmark BridgeForEvenShift $$ will always work as a \"{bridge},
as one of them will correctly adjust the parity of the amount of $1$'s between the last $\qt 0$ of $\beta $ and the
first $\qt 0$ of $x$.

However the spectral partial action for the even shift is not minimal.  To see this let us return to the situation
presented in \ref {IsolatedPoint}, when we have verified that the open set $V_{\emptyword ; \beta _1, \beta _2}$
consists of a single point, namely $\xi _{1^\infty }$.  Should the spectral partial action for the even shift be
minimal, then $\Orb (\xi _x)$ should intersect $V_{\emptyword ; \beta _1, \beta _2}$, for every $x$ in $X$, but this is
not the case, e.g.~for the element $x=\qt {00000\ldots }$.  In fact, if $$ \vartheta _g(\xi _{0^\infty })\in
V_{\emptyword ; \beta _1, \beta _2}, $$ for some $g$ in $\F $, then $\vartheta _g(\xi _{0^\infty })=\xi _{1^\infty }$,
whence by \ref {Equivariance} one has that $\theta _g(0^\infty )=1^\infty $, and this is manifestly impossible.

For future reference let us highlight the conclusion reached above:

\state Proposition \label EvenNotMinimal The spectral partial action associated to the even shift is not minimal.

As in the previous section, one might decide to ignore that the standard partial action for the even shift is not a
topological partial action and ask whether or not it is minimal in the sense that all orbits are dense.  The answer is
then easily seen to be positive, so in this sense the standard partial for the even shift is minimal.

Even though cofinality for a subshift does not imply minimality for the spectral partial action, as in the Markov case,
the former concept is quite relevant in our study of minimality, having the following dynamical interpretation.

\state Proposition \label CofinalityProps Let $X$ be a subshift.  Then the following are equivalent: \izitem \zitem $X$
is cofinal, \zitem for every $x$ in $X$, and for every $\beta $ in $\Lang $, there is some $g$ in $\F $ such that both
$g$ and $g\beta $ lie in $\xi _x$, \zitem for every $x$ in $X$, and for every $\beta $ in $\Lang $, one has that $\Orb
(\xi _x)\cap \Omega _\beta \neq \emptyset $.

\Proof \ipfimply \pfimply (i)(ii) Given $x$ and $\beta $, as in (ii), choose $\alpha $, $\gamma $ and $y$, as in (i), so
that $x=\alpha y$, and $\beta \gamma y\in X$.  We then have that $\stem _\alpha (\xi _x)=y$, so the fact that $y\in
\Folow _\gamma $ yields $\alpha \gamma \inv \in \xi _x$, by \ref {Enfiabilidade}.  Similarly, the fact that $y\in \Folow
_{\beta \gamma }$ yields $\alpha (\beta \gamma )\inv \in \xi _x$.  It is then clear that $g=\alpha (\beta \gamma )\inv $
satisfies the conditions of (ii).

\pfimply (ii)(i) Given $\beta $ in $\Lang $, and $x$ in $X$, choose $g$ in $\F $ such that $g,g\beta \in \xi _x$.  By
\ref {CharacXix} we may write $g\beta =\alpha \gamma \inv $ in reduced form, with $\alpha ,\gamma \in \F _+$, such that
$x=\alpha y$, for some infinite word $y\in \Folow _\alpha \cap \Folow _\gamma $.  Observing that $\xi _x$ contains both
$\alpha $ and $$ g = g\beta \beta \inv = \alpha \gamma \inv \beta \inv = \alpha (\beta \gamma )\inv , $$ we deduce from
\ref {Enfiabilidade} that $y=\stem _\alpha (\xi _x)\in \Folow _{\beta \gamma }$, whence $\beta \gamma y\in X$.  This
shows that $x$ may be reached from $\beta $.

The equivalence of (ii) and (iii) follows from the next simple result which we will also use later in a slightly more
general situation. \endProof

\state Proposition \label gbOrbit Let $X$ be a subshift and let $\xi \in \Spec $.  Given $g\in \F $, and $\beta \in
\Lang $, the following are equivalent: \izitem \zitem both $g$ and $g\beta $ lie in $\xi $, \zitem $\xi \in \Omega _g$,
and $\vartheta _{g\inv }(\xi )\in \Omega _\beta $.

\Proof Notice that $g\in \xi $, if and only if $\xi \in \Omega _g$, by \ref {DescrPAOnOmega.i}, and in this case $$
g\beta \in \xi \iff \beta \in g\inv \xi = \vartheta _{g\inv }(\xi ) \iff \vartheta _{g\inv }(\xi )\in \Omega _\beta .
\endProof

Recall that $\vartheta $ is minimal if and only if, for every $\xi $ in $\Spec $, and for every nonempty open set
$U\subseteq \Spec $, one has that $$ \Orb (\xi )\cap U\neq \emptyset .  \equationmark MinimalityExplained $$

Expressing cofinality in terms of \ref {CofinalityProps.iii} may thus be interpreted as a weak form of minimality in the
sense that the above holds when $\xi $ has the form $\xi _x$, and $U$ has the form $\Omega _\beta $.  One therefore has
a clear perspective that cofinality is a weaker property than minimality, and in fact strictly weaker, as the example of
the even shift shows.

If we are to find a set of properties related to cofinality, characterizing minimality for the spectral partial action,
we must therefore significantly strengthen the notion of cofinality.  We will in fact do this in two directions,
introducing the notions of \"{\colcofty } and \"{\scofty }.  The first will lead to Theorem \ref {ColcofOrbX}, ensuring
that \ref {MinimalityExplained} holds for an arbitrary $U$, but still under the restriction that $\xi =\xi _x$, while
the second will show up in Theorem \ref {BddCost}, yielding \ref {MinimalityExplained} for an arbitrary $\xi $, but
under the restriction that $U=\Omega _\beta $.  Together these two will eventually prove $\vartheta $ to be minimal.

\definition \label ReachabilityMultiple Let $\Lambda $ be a finite alphabet and let $X\subseteq \Lambda ^{\bf N}$ be a
subshift.  Given a finite set of finite words $B\subseteq \Lang $, and an infinite word $x\in X$, we shall say that $x$
may be \"{collectively reached} from $B$, when there exists a finite word $\gamma $, and a prefix $\alpha $ of $x$, such
that, upon writing $x=\alpha y$, one has that $\beta \gamma y\in X$, for all $\beta \in B$.

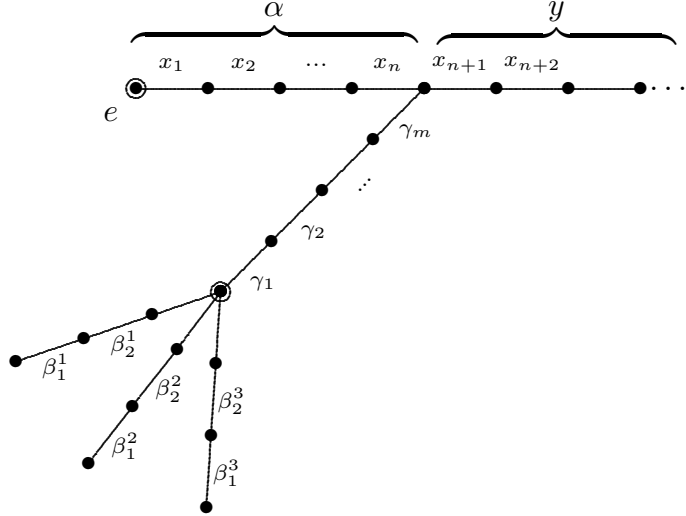
\null \hskip 2cm \beginpicture \setcoordinatesystem units <0.0080truecm, -0.0080truecm> point at 1000 1000 \put {$e$} at
-42 42 \circulararc 360 degrees from 0 16 center at 0 0 \put {$\bullet $} at 0 0 \plot 0 0 120 0 / \put {$\bullet $} at
120 0 \put {$\scriptstyle x_1$} at 60 -40 \plot 120 0 240 0 / \put {$\bullet $} at 240 0 \put {$\scriptstyle x_2$} at
180 -40 \plot 240 0 360 0 / \put {$\bullet $} at 360 0 \put {$\scriptstyle \ldots $} at 300 -40 \plot 360 0 480 0 / \put
{$\bullet $} at 480 0 \put {$\scriptstyle x_n$} at 420 -40 \plot 480 0 600 0 / \put {$\bullet $} at 600 0 \put
{$\scriptstyle x_{n+1}$} at 540 -40 \plot 600 0 720 0 / \put {$\bullet $} at 720 0 \put {$\scriptstyle x_{n+2}$} at 660
-40 \plot 720 0 840 0 / \put {$\bullet $} at 840 0 \put {$\ldots $ } at 895 0 \put {$\overbrace {\hbox to 90pt{\hfill
}}^{\textstyle \alpha }$} at 230 -110 \put {$\overbrace {\hbox to 75pt{\hfill }}^{\textstyle y}$} at 700 -110 \put
{$\bullet $} at 141 339 \plot 141 339 225 255 / \put {$\bullet $} at 225 255 \put {$\scriptstyle \gamma _1$} at 211 325
\plot 225 255 310 170 / \put {$\bullet $} at 310 170 \put {$\scriptstyle \gamma _2$} at 296 240 \plot 310 170 395 85 /
\put {$\bullet $} at 395 85 \put {$\scriptstyle \nedots $} at 381 156 \plot 395 85 480 0 / \put {$\bullet $} at 480 0
\put {$\scriptstyle \gamma _m$} at 466 71 \circulararc 360 degrees from 141 355 center at 141 339 \put {$\bullet $} at
-200 455 \plot -200 455 -87 417 / \put {$\bullet $} at -87 417 \put {$\scriptstyle \beta ^1_1$} at -134 464 \plot -87
417 27 378 / \put {$\bullet $} at 27 378 \put {$\scriptstyle \beta ^1_2$} at -20 426 \plot 27 378 141 339 / \put
{$\bullet $} at 141 339 \put {$\scriptstyle \null $} at 93 387 \put {$\bullet $} at -79 625 \plot -79 625 -6 530 / \put
{$\bullet $} at -6 530 \put {$\scriptstyle \beta ^2_1$} at -18 596 \plot -6 530 68 435 / \put {$\bullet $} at 68 435
\put {$\scriptstyle \beta ^2_2$} at 55 500 \plot 68 435 141 339 / \put {$\bullet $} at 141 339 \put {$\scriptstyle \null
$} at 128 405 \put {$\bullet $} at 117 699 \plot 117 699 125 579 / \put {$\bullet $} at 125 579 \put {$\scriptstyle
\beta ^3_1$} at 151 641 \plot 125 579 133 459 / \put {$\bullet $} at 133 459 \put {$\scriptstyle \beta ^3_2$} at 159 521
\plot 133 459 141 339 / \put {$\bullet $} at 141 339 \put {$\scriptstyle \null $} at 167 401 \endpicture \hfill \null

\bigskip It should be stressed that the \"{bridge} $\gamma $ above is supposed to be the same for all $\beta $.  In fact
an infinite word $x$ may be reached \"{individually} from both $\beta _1$ and $\beta _2$, but not \"{collectively} by
the set $\{\beta _1,\beta _2\}$, as is the case of the even shift with $$ \beta _1=01, \quad \beta _2=011, \and
x=00000\ldots $$

In case a word $x$ is collectively reached from a subset $B\subseteq \Lang $, as in the diagram above, notice that
$\gamma y$ lies in the follower set $\Folow _\beta $, for all $\beta $ in $B$, and in particular the intersection of the
$\Folow _\beta $ is nonempty.  This motivates an extension of the notion of follower sets:

\definition Given a finite subset $B\subseteq \Lang $, we will say that the \"{follower set} of $B$ is the set $\Folow
_B$ defined by $$ \Folow _B =\medcap _{\beta \in B}\Folow _\beta .  $$

\definition We will say that a subshift $X$ is \"{\colcof } if, for every finite set $B\subseteq \Lang $, with $\Folow
_B$ nonempty, and every $x$ in $X$, one has that $x$ may be collectively reached from $B$.

The first relationship between {\colcofty } and minimality is as follows:

\state Theorem \label ColcofOrbX Let $X$ be a subshift.  Then the following are equivalent: \izitem \zitem for every $x$
in $X$, the orbit of $\xi _x$ is dense in $\Spec $, \zitem $X$ is {\colcof }.

\Proof \ipfimply \pfimply (i)(ii) Given any finite subset $B\subseteq \Lang $, such that $\Folow _B$ is nonempty, say
$B=\{\beta _1,\beta _2,\ldots ,\beta _n\}$, we claim that $$ V_{\emptyword ;\,\beta _1,\beta _2,\ldots ,\beta _n}\neq
\emptyset .  $$ To see this, pick $y$ in $\Folow _B$, and notice that $\beta \inv \in \xi _y$, for every $\beta $ in
$B$, by \ref {BinvInXix}, whence $\xi _y$ is an element in the above set, proving it to be nonempty.

Having already chosen $B$, let us also choose any $x$ in $X$, and our task is to show that $x$ may be collectively
reached from $B$.  By hypothesis the orbit of $\xi _x$ is dense in $\Spec $, so there exists $g$ in $\F $ such that $\xi
_x\in \Omega _{g\inv }$, and $$ \vartheta _g(\xi _x)\in V_{\emptyword ;\,\beta _1,\beta _2,\ldots ,\beta _n}.
\equationmark XiXInV $$ It follows that $g\inv \in \xi _x$, by \ref {DescrPAOnOmega.i}, whence we may use \ref
{CharacXix} to write $g\inv =\alpha \gamma \inv $, in reduced form, with $\alpha ,\gamma \in \F _+$, and moreover write
$x=\alpha y$, where $y\in \Folow _\alpha \cap \Folow _\gamma $.  Therefore $$ \vartheta _g(\xi _x) \={Equivariance} \xi
_{\theta _g(x)} = \xi _{\gamma y}.  $$ It then follows from \ref {XiXInV} that $\beta _i\inv \in \xi _{\gamma y}$, for
every $i$, so $\gamma y\in \Folow _{\beta _i}$, by \ref {BinvInXix}, which is the same as saying that $\beta _i\gamma
y\in X$.  This shows that $x$ may be collectively reached from $B$, and hence proves (ii).

\pfimply (ii)(i) Given any $x$ in $X$, we must show that the orbit of $\xi _x$ is dense.  Since $\XO (X)$ is already
known to be dense in $\Spec $, it is enough to prove that $$ \xi _z\in \overline {\Orb (\xi )}, \for z\in X.  $$

Given any $z$ in $X$, let $U$ be an arbitrary neighborhood of $\xi _z$.  By \ref {VerySimpleNBDs} there are finite words
$\alpha ,\ \beta _1,\beta _2,\ldots ,\beta _n$, such that $$ \xi _z\in V_{\alpha ;\,\beta _1,\beta _2,\ldots ,\beta _n}
\subseteq U, \equationmark XisHere $$ and all we need to do is prove that $$ \Orb (\xi _x)\cap V_{\alpha ;\,\beta
_1,\beta _2,\ldots ,\beta _n} \neq \emptyset .  \equationmark Bejetivo $$

Observing that $V_{\alpha ;\,\beta _1,\beta _2,\ldots ,\beta _n}\subseteq \Omega _\alpha $, and that $$ \def \quad {}
\matrix { \vartheta _{\alpha \inv }(&V_{\alpha ;\kern 10pt \beta _1,\beta _2,\ldots ,\beta _n})& \ = \cr \pilar {12pt} &
V_{\emptyword ;\,\alpha ,\beta _1,\beta _2,\ldots ,\beta _n},} $$ notice that any orbit intersecting $V_{\emptyword
;\,\alpha ,\beta _1,\beta _2,\ldots ,\beta _n}$ will also intersect $V_{\alpha ;\,\beta _1,\beta _2,\ldots ,\beta _n}$.
So we may assume without loss of generality that $\alpha =\emptyword $.

Letting $B = \{\beta _1,\beta _2,\ldots ,\beta _n\}$, we claim that $\Folow _B$ is nonempty.  In fact, by \ref {XisHere}
we have that $\beta _i\inv \in \xi _z$, for all $i$, so $z\in \Folow _{\beta _i}$ by \ref {BinvInXix}, whence also $z\in
\Folow _B$, proving our claim.

By hypothesis $x$ may be collectively reached from $B$, so there are finite words $\mu $ and $\gamma $, such that $x=\mu
y$, and $\beta \gamma y\in X$, for all $\beta $ in $B$.  Noticing that for all such $\beta $, $$ \stem _{\mu }(\xi _x) =
y \in \Folow _{\beta \gamma }, $$ we have that $$ \mu (\beta \gamma )\inv \in \xi _x, \equationmark SomeMembers $$ by
\ref {Enfiabilidade}.  The fact that $\beta \gamma y\in X$ implies that $\gamma y\in X$ as well, so a similar argument
gives $\mu \gamma \inv \in \xi _x$, whence $\xi _x\in \Omega _{\mu \gamma \inv }$, and $$ \vartheta _{\gamma \mu \inv
}(\xi _x) = \gamma \mu \inv \xi _x \explain {SomeMembers}\ni \beta \inv .  $$ Consequently $\vartheta _{\gamma \mu \inv
}(\xi _x) \in V_{\emptyword ;\,\beta _1,\beta _2,\ldots ,\beta _n}$, proving \ref {Bejetivo}.  \endProof

The reader is invited to check that the only reason why \ref {ColcofOrbX.i} must be stated just for $\xi _x$, rather
than for a general $\xi $ in $\Spec $, is that the above proof has used the crucial implication ``\implica (b)(a)" of
\ref {Enfiabilidade}, which does not hold in general.

In the presence of {\colcofty }, topological freeness for the spectral partial action may be characterized in a very
simple way:

\state Proposition \label TopFreeWhenColCof Let $X$ be a subshift, and consider the following statements: \iaitem \aitem
The spectral partial action associated to $X$ is topologically free, \aitem $X$ has at least one point which is not
eventually periodic (by eventually periodic we mean a word of the form $\alpha \gamma ^\infty $), \medskip \noindent
Then {\rm (a)} implies {\rm (b)}.  If $X$ is {\colcof }, then also {\rm (b)} implies {\rm (a)}.

\Proof \ipfimply \pfimply (a)(b) It is easy to see that there is only a countable number of eventually periodic infinite
words (even within the full shift $\Lambda ^{\bf N}$).  Assuming by contradiction that every $x$ in $X$ is eventually
periodic, one then has that $X$ is countable, and a standard application of Baire's category Theorem implies that $X$
has at least one isolated point, say $x$.  Since $x$ is eventually periodic let us write $x=\alpha \gamma ^\infty $, for
some circuit $\gamma $.

To say that $x$ is isolated (relative to the topology of $X$), is to say that $\{x\}$ coincides with some cylinder $\Cyl
_\mu $, where $\mu $ is necessarily a prefix of $x=\alpha \gamma ^\infty $.  By increasing the size of $\mu $ up to
$|\alpha |$ plus some multiple of $|\gamma |$, we may suppose that $\mu =\alpha \gamma ^n$, for some integer $n$.

Thus the only infinite word in $X$ that one may produce by starting with $\alpha \gamma ^n$ is $\alpha \gamma ^\infty $,
which may also be expressed by saying that $\Folow _{\alpha \gamma ^n}=\{\gamma ^\infty \}$.  This goes against \ref
{MainTopFree.ii}, and hence also against (a), thus bringing about a contradiction, proving that $X$ must have some non
eventually periodic point.

\pfimply (b)(a) We will prove (a) by verifying \ref {MainTopFree.ii}.  So let $\nu $ be a circuit and let $B$ be a
finite set of finite words such that $\nu ^\infty \in \Folow _B$.  Our task is then to find some element in $\Folow _B$
other than $\nu ^\infty $.

For this, pick any non eventually periodic element $x$ of $X$ and, using that $X$ is {\colcof }, choose a prefix $\alpha
$ of $x$, and a finite word $\gamma $, such that, upon writing $x=\alpha y$, one has that $\beta \gamma y\in X$, for all
$\beta $ in $B$.  This implies that $\gamma y$ lies in $\Folow _B$, and since it is clearly not eventually periodic, one
necessarily has that $\gamma y\neq \nu ^\infty $, thus verifying \ref {MainTopFree.ii}.  \endProof

Having understood the dynamical meaning of {\colcofty } in \ref {ColcofOrbX}, let us now discuss yet another version of
cofinality.  In order to motivate this notion let us refer back to diagram \ref {PombaSimplesDiagrama}.  One might see
the question of reaching $x$ from $\beta $ as an attempt to bringing $x$ into the follower set of $\beta $ by deleting a
few letters from the beginning of $x$, namely the prefix $\alpha $, and then inserting new letters in its place, namely
$\gamma $, after which the resulting word $\gamma y$ is supposed to lie in $\Folow _\beta $.  If we have to pay a price
for each deleted, as well as for each inserted letter, then the cheapest situation is obviously when $x$ itself lies in
the follower set of $\beta $, namely when $\alpha =\gamma =\emptyword $ are enough to do the job.  Beyond this ideal
situation we have:

\definition Let $\beta \in \Lang $, and let $x\in X$.  We shall say that the \"{cost} of reaching $x$ from $\beta $ is
the minimum value of $|\alpha |+|\gamma |$, where $\alpha $ and $\gamma $ are as in \ref {Reachability}.  In symbols, $$
\cost (\beta ,x)=\min \{|\alpha |+|\gamma |: x=\alpha y,\ \beta \gamma y\in X\}.  $$ If no such $\alpha $ and $\gamma $
exist, we shall say that the cost if infinite.

In a cofinal subshift, given $\beta $ in $\Lang $, one may reach any $x$ in $X$, but perhaps at an increasingly high
cost.  This motivates the following:

\definition Let $X$ be a subshift.  We shall say that $X$ is \"{\scof } if $$ \sup _{x\in X}\cost (\beta ,x)<\infty , $$
for every $\beta $ in $\Lang $.

Let us pause for a moment to give an example of a cofinal subshift which is not {\scof }.  Our staple counter-example,
namely the even shift, will not serve us now since it is both cofinal and {\scof }.  In fact, as we already mentioned in
$\ref {BridgeForEvenShift}$, when trying to reach an infinite word $x$ from a finite word $\beta $, it suffices to take
$\alpha =\emptyword $, and either $\gamma =\emptyword $ or $\gamma =\qt 1$, which means that $\cost (\beta ,x)$ is at
most $1$.

However, a subshift based on a similar principle as the even shift will provide the counter-example sought.  Consider
the alphabet $\Lambda =\{0,1\}$, and let us take the following set $\Forbid $ of forbidden words: $$ \Forbid = \{01^n0:
n \hbox { is \underbar {not} a power of } 2\}.  $$

If $X=X_\Forbid $ is the corresponding subshift, then an infinite word $x$ lies in $X$ if and only if, anytime a
contiguous block of 1's occurring in $x$ is delimited by 0's, the amount of said 1's is a power of $2$.

We leave it for the reader to check that $X$ is a duly cofinal subshift, but that reaching the infinite word
$x={1^n0^\infty }$ from the finite word $\beta =\qt {0}$ involves unbounded costs.  Hence $X$ is not {\scof }.

Exploring the relationship between {\scofty } and the dynamical properties of the spectral partial action is our next
goal.

\state Theorem \label BddCost Let $X$ be a subshift.  Then the following are equivalent: \izitem \zitem For every $\beta
$ in $\Lang $, and for every $\xi $ in $\Spec $, one has that $\Orb (\xi )\cap \Omega _\beta $ is nonempty, \zitem X is
{\scof }.

\Proof \ipfimply \pfimply (i)(ii) We will in fact prove that the negation of (ii) implies the negation of (i).  We
therefore suppose that $$ \sup _{x\in X}\cost (\beta ,x)=\infty , $$ for some $\beta $ in $\Lang $.  Then, for every
natural number $n$, there is some $z_n$ in $X$, with $\cost (\beta ,z_n)\geq n$.  By compactness we may then find a
subsequence $x_k=z_{n_k}$, such that $\{\xi _{x_k}\}_k$ converges to some $\xi $ in $\Spec $.  We will accomplish our
task by proving that $$ \Orb (\xi )\cap \Omega _\beta =\emptyset .  \equationmark OrbNoMeet $$

Arguing by contradiction, suppose that there exists some $h$ in $\F $ such that $\xi \in \Omega _{h\inv }$, and
$\vartheta _h(\xi )\in \Omega _\beta $. By \ref {gbOrbit} we deduce that both $h\inv $ and $h\inv \beta $ lie in $\xi $.
In order to simplify our notation, we shall make the change of variables $g=h\inv \beta $, whence $$ g,\ g\beta \inv \in
\xi .  $$

Using \ref {PropForPaint.v} we may write $g=\alpha \gamma \inv $ in reduced form, so that $\alpha \in \xi $, by
{\convexity }.  Employing \ref {BooleanContinuous} we then have that $$ \alpha ,\ g,\ g\beta \inv \in \xi _{x_k}, $$ for
all large enough $k$.  Focusing on the fact that $$ \alpha \gamma \inv = g \in \xi _{x_k}, $$ if follows from \ref
{CharacXix} that $x_k=\alpha y_k$, where $y_k\in \Folow _\alpha \cap \Folow _\gamma $.  Observing that the stem of $\xi
_{x_k}$ at $g$ is $\gamma y_k$, and that $g\beta \inv \in \xi _{x_k}$, we conclude from \ref {Enfiabilidade} that
$\gamma y_k\in \Folow _\beta $, which is to say that $\beta \gamma y_k\in X$.  This implies not only that $x_k$ may be
reached from $\beta $, but also that $$ \cost (\beta ,x_k) \leq |\alpha |+|\gamma |.  $$ Since neither $\alpha $ nor
$\gamma $ depend on $k$, we arrive at a contradiction with the fact that $\lim _{k\to \infty }\cost (\beta ,x_k)=\infty
$.  This proves \ref {OrbNoMeet}, as desired.

\pfimply (ii)(i) Given $\beta $ in $\Lang $, and $\xi $ in $\Spec $, we must prove that $\Orb (\xi )\cap \Omega _\beta
\neq \emptyset $.  Writing $\xi =\lim _n\ \xi _{x_n}$, with $x_n\in X$, by hypothesis we may reach every $x_n$ from
$\beta $ with bounded cost, meaning that for every $n$, there are finite words $\alpha _n$ and $\gamma _n$, such that
$x_n=\alpha _ny_n$, and $\beta \gamma _ny_n\in X$, and moreover $|\alpha _n|+|\gamma _n|\leq M$, where $M$ is a fixed
constant.

As we are working with a finite alphabet, there are finitely many words of any given length.  Therefore the $\alpha _n$
and $\gamma _n$ just obtained must necessarily repeat infinitely many often.  We may then choose a subsequence
$\{x_{n_k}\}_k$, such that $\alpha _{n_k}=\alpha $, and $\gamma _{n_k}=\gamma $, for all $k$.  Therefore $x_{n_k}=\alpha
y_{n_k}$, and $\beta \gamma y_{n_k}\in X$.  We then have that $$ \alpha \in \xi _{x_{n_k}}, \and \stem _\alpha (\xi
_{x_{n_k}}) = y_{n_k} \in \Folow _{\beta \gamma }, $$ whence $\alpha (\beta \gamma )\inv \in \xi _{x_{n_k}}$ by \ref
{Enfiabilidade}.  A similar reasoning, based on the fact that $y_{n_k}$ also lies in $\Folow _\gamma $, gives $\alpha
\gamma \inv \in \xi _{x_{n_k}}$, so by the continuity of Boolean values \ref {BooleanContinuous}, we have that $$ \alpha
\gamma \inv \in \xi , \and \alpha (\beta \gamma )\inv \in \xi .  $$ Setting $g=\alpha (\beta \gamma )\inv $, we have
that $$ g\beta = \alpha \gamma \inv \beta \inv \beta = \alpha \gamma \inv \in \xi , $$ whence $\vartheta _g\inv (\xi
)\in \Omega _\beta $ by \ref {gbOrbit}.  This shows that the orbit of $\xi $ intersects $\Omega _\beta $, concluding the
proof.  \endProof

From \ref {ColcofOrbX} and \ref {BddCost} it is clear that a subshift whose associated spectral partial action is
minimal must necessarily be both {\colcof } and {\scof }.  Our next major goal will be to prove that these two
properties in turn characterize minimality for the spectral partial action.  However there is an important technical
tool we still need to develop before proving this main result.

\state Lemma \label CrazyInclusion Let $X$ be a {\colcof } subshift.  Then for every finite set $B\subseteq \Lang $,
with $\Folow _B$ nonempty, there are $\mu $ and $\nu $ in $\Lang $, such that $$ \Folow _\mu \subseteq \Folow _{B\nu },
$$ where $B\nu $ evidently means the set $\{\beta \nu :\beta \in B\}$.

\Proof Let $B$ be as in the statement.  By hypothesis any $x$ in $X$ may be collectively reached from $B$, so we may
pick $\alpha $, $\gamma $ and $y$, as in \ref {ReachabilityMultiple}, so that $y\in \Folow _\alpha \cap \Folow _{B\gamma
}$, and $x=\theta _\alpha (y)$, whence $ x\in \theta _\alpha (\Folow _\alpha \cap \Folow _{B\gamma }).  $ This implies
that $$ X=\medcup _{\alpha ,\gamma \in \Lambda ^*}\theta _\alpha (\Folow _\alpha \cap \Folow _{B\gamma }).  $$

Observing that $\Folow _\alpha \cap \Folow _{B\gamma }$ is compact, and hence that the sets in the above union are
closed in $X$, we may employ Baire's category Theorem producing $\alpha $ and $\gamma $ in $\Lambda ^*$ such that
$\theta _\alpha (\Folow _\alpha \cap \Folow _{B\gamma })$ has a nonempty interior.  En passant we stress that both
$\alpha $ and $\gamma $ must lie in $\Lang $, or else $\Folow _\alpha \cap \Folow _{B\gamma }$ is the empty set.
Therefore, by the definition of the product topology on $X$, there exists some $\mu $ in $\Lang $ such that $$ \Cyl _\mu
\subseteq \theta _\alpha (\Folow _\alpha \cap \Folow _{B\gamma }).  $$ Assuming without loss of generality that $|\mu
|>|\alpha |$, and noticing that the range of $\theta _\alpha $ is contained in $\Cyl _\alpha $, we have that $\Cyl _\mu
\subseteq \Cyl _\alpha $, so $\alpha $ must be a prefix of $\mu $, and we may then write $\mu =\alpha \delta $, for some
finite word $\delta \in \Lang $.  We will now conclude the proof by showing that the inclusion in the statement holds
once we choose $\nu =\gamma \delta $.

To prove this, let $y\in \Folow _\mu $, so $$ x:=\mu y \in \Cyl _\mu \subseteq \theta _\alpha (\Folow _\alpha \cap
\Folow _{B\gamma }), $$ and then we may write $x=\alpha z$, for some $z\in \Folow _\alpha \cap \Folow _{B\gamma }$.
Since $$ \alpha \delta y = \mu y = x = \alpha z, $$ it follows that $\delta y=z$, so for any $\beta $ in $B$, we have
that $$ \beta \nu y = \beta \gamma \delta y = \beta \gamma z \in X.  $$ This shows that $y\in \Folow _{\beta \nu }$, as
desired, concluding the proof.  \endProof

Let us now give a dynamical interpretation of the conclusion of the above result.

\state Proposition \label CrazyInclusionDynamical Let $B=\{\beta _1,\beta _2,\ldots ,\beta _n\}\subseteq \Lang $, be a
nonempty finite set and let $\mu ,\nu \in \Lang $, be such that $ \Folow _\mu \subseteq \Folow _{B\nu }, $ precisely as
in the conclusion of \ref {CrazyInclusion}.  Then $\Omega _{\mu \inv }$ is contained in the domain of $\vartheta _\nu $
(also known as $\Omega _{\nu \inv }$), and $$ \vartheta _\nu (\Omega _{\mu \inv }) \subseteq V_{\emptyword ;\,\beta
_1,\beta _2,\ldots ,\beta _n}.  $$

\Proof We will first prove that $$ \XO (X)\cap \Omega _{\mu \inv }\subseteq \Omega _{\nu \inv }.  \equationmark
InclusionWithX $$ We thus pick any $x$ in $X$ such that $\xi _x\in \Omega _{\mu \inv }$.  Notice that $$ \xi _x\in
\Omega _{\mu \inv } \explain {DescrPAOnOmega.i}\iff \mu \inv \in \xi _x \explain {BinvInXix}\iff x\in \Folow _\mu .  $$

So $x\in \Folow _\mu $, and by hypothesis we then have that $x\in \Folow _{B\nu }$.  Since $B$ is nonempty, we may pick
any $\beta $ in $B$, and so deduce that $x\in \Folow _{\beta \nu }$, whence $\beta \nu x\in X$.  This implies that $x\in
\Folow _\nu $, and a reasoning similar to the above shows that $\xi _x\in \Omega _{\nu \inv }$, proving \ref
{InclusionWithX}.  A trivial exercise in Topology, using that $\XO (X)$ is dense, $\Omega _{\mu \inv }$ is open, and
$\Omega _{\nu \inv }$ is closed, now shows that $\Omega _{\mu \inv }\subseteq \Omega _{\nu \inv }$, as required.

We will next show that $$ \vartheta _\nu \big (\XO (X)\cap \Omega _{\mu \inv }\big ) \subseteq V_{\emptyword ;\,\beta
_1,\beta _2,\ldots ,\beta _n}.  \equationmark InclusionWithXandTheta $$ We thus again pick any $x$ in $X$ such that $\xi
_x\in \Omega _{\mu \inv }$.  Then $$ \vartheta _\nu (\xi _x) \={Equivariance} \xi _{\theta _\nu (x)} = \xi _{\nu x}.  $$
In addition, as seen above, $x\in \Folow _{B\nu }$, so $\beta \nu x\in X$, for every $\beta $ in $B$, and we see that
$\nu x\in \Folow _\beta $.  We then have by \ref {BinvInXix} that $\beta \inv \in \xi _{\nu x}$, whence $\xi _{\nu x}\in
V_{\emptyword ;\,\beta _1,\beta _2,\ldots ,\beta _n}$, proving \ref {InclusionWithXandTheta}.  The same ``dense plus
open plus closed" argument above, but now also using that $\vartheta _\nu $ is continuous, leads one from \ref
{InclusionWithXandTheta} to the last conclusion in the statement.  \endProof

We have now come to the main result of this section:

\state Theorem \label MainMinimal Let $\Lambda $ be a finite alphabet and let $X\subseteq \Lambda ^{\bf N}$ be a
subshift.  Also let $ \vartheta $ be the spectral partial action of the free group $\F (\Lambda )$ on $\Spec $
introduced in \ref {SpectralAction}.  Then a necessary and sufficient condition for $\vartheta $ to be minimal is that
$X$ be both {\colcof } and {\scof }.

\Proof As already observed, the necessity of the above conditions follows immediately from \ref {ColcofOrbX} and \ref
{BddCost}.

Conversely, suppose that $X$ is both {\colcof } and {\scof }.  In order to prove that $\vartheta $ is minimal, given any
$\xi $ in $\Spec $, we must show that the orbit of $\xi $ is dense.  Arguing as in the proof of ``\implica (ii)(i)" in
\ref {ColcofOrbX}, it suffices to show that if $B=\{\beta _1,\beta _2,\ldots ,\beta _n\}$ is a subset of $\Lang $ such
that $\Folow _B$ is nonempty, then $$ \Orb (\xi )\cap V_{\emptyword ;\,\beta _1,\beta _2,\ldots ,\beta _n}\neq \emptyset
.  $$

Using \ref {CrazyInclusion} and \ref {CrazyInclusionDynamical}, pick $\mu $ and $\nu $ in $\Lang $, such that $$
\vartheta _\nu (\Omega _{\mu \inv }) \subseteq V_{\emptyword ;\,\beta _1,\beta _2,\ldots ,\beta _n}.  $$ Thus, to show
that the orbit of $\xi $ intersects $V_{\emptyword ;\,\beta _1,\beta _2,\ldots ,\beta _n}$, it is enough to show that it
intersects $\Omega _{\mu \inv }$, or even $\Omega _\mu $, since $\vartheta _{\mu \inv }(\Omega _\mu ) = \Omega _{\mu
\inv }$.

In order to prove the latter fact, that is, that $\Orb (\xi )\cap \Omega _\mu $ is nonempty, we just have to note that
it follows immediately from {\scofty } and \ref {BddCost}.  \endProof

Of course one could put together the notions of {\colcofty } and {\scofty }, to form a new high powered notion, as
follows:

\definition Let $X$ be a subshift.  \iaitem \aitem Given a finite subset $B$ of $\Lang $, we shall say that the \"{cost}
of reaching a given $x$ in $X$ from $B$ is the minimum value of $|\alpha |+|\gamma |$, where $\alpha $ and $\gamma $ are
as in \ref {ReachabilityMultiple}.  In symbols, $$ \cost (B,x)=\min \{|\alpha |+|\gamma |: x=\alpha y,\ \beta \gamma
y\in X, \hbox { for all } \beta \in B\}.  $$ If no such $\alpha $ and $\gamma $ exist, we shall say that the cost if
infinite.  \aitem We shall say that $X$ is \"{\hcof } if $$ \sup _{x\in X}\cost (B,x)<\infty , $$ for every finite
subset $B$ of $\Lang $ such that $\Folow _B$ is nonempty.

It is evident that {\hcofty } implies both {\colcofty } and {\scofty }, but the reverse implication is not so obvious to
see.  Nevertheless it is true, as we shall now prove.

\state Proposition \label MainHyper Let $X$ be a subshift.  Then the following are equivalent: \izitem \zitem $X$ is is
both {\colcof } and {\scof }, \zitem $X$ is {\hcof }, \zitem the spectral partial action associated to $X$ is minimal.

\Proof The equivalence of (i) and (iii), repeated here just for emphasis, is precisely the content of \ref
{MainMinimal}.  As we have already seen (ii) implies (i) for obvious reasons, so it suffices to show that (iii) implies
(ii).

So, assuming that $\vartheta $ is minimal, let $B=\{\beta _1,\beta _2,\ldots ,\beta _n\}\subseteq \Lang $ be such that
$\Folow _B$ is nonempty.  Consider the open subset $$ V:= V_{\emptyword ;\,\beta _1,\beta _2,\ldots ,\beta _n}\subseteq
\Spec , $$ which is nonempty by the short argument used in the beginning of the proof of \ref {ColcofOrbX}.  By
minimality, for every $\xi $ in $\Spec $, one has that $\Orb (\xi )$ is dense, so it must intersect $V$, meaning that
there exists some $g$ in $\F $, such that $\xi \in \Omega _{g\inv }$, and $\vartheta _g(\xi )\in V$.  It follows that
$\xi \in \vartheta _{g\inv }(V\cap \Omega _g)$, so we conclude that $$ \Spec = \medcup _{g\in \F }\vartheta _{g\inv
}(V\cap \Omega _g).  $$

This is an open covering of the compact space $\Spec $, so there exists a finite collection $g_1,g_2,\ldots ,g_m$ of
elements in $\F $ such that $$ \Spec = \medcup _{i=1}^m\vartheta _{g_i\inv }(V\cap \Omega _{g_i}).  $$

Given any $x$ in $X$, we may then find some $i\leq m$, such that $\xi _x\in \vartheta _{g_i\inv }(V\cap \Omega _{g_i})$,
so in particular $\xi _x\in \Omega _{g_i\inv }$ and $\vartheta _{g_i}(\xi _x)\in V$.  Writing $g_i=\gamma _i\alpha
_i\inv $, in reduced form, we then have that $\alpha _i$ is a prefix of $x$ and, upon writing $x=\alpha _iy$, one has
that $$ \vartheta _{g_i}(\xi _x)=\xi _{\theta _{g_i}(x)} = \xi _{\gamma _iy}\in V=V_{\emptyword ;\,\beta _1,\beta
_2,\ldots ,\beta _n}.  $$ Consequently $\beta _j\inv \in \xi _{\gamma _iy}$, for all $j\leq n$, whence $\gamma _iy\in
\Folow _{\beta _j}$ by \ref {BinvInXix}, which is to say that $\beta _j\gamma _iy$ lies in $X$.  This shows that $x$ may
be reached from $B$, and moreover that $$ \cost (B,x)\leq |\alpha _j|+|\gamma _j|.  $$ Since the $j$ above may take only
finitely many values, we see that the cost of reaching any $x$ from the given $B$ is bounded, whence $X$ is {\hcof }.
\endProof

Not all subshifts are surjective\fn {If a subshift is built from a two-sided shift, as in many papers on the subject
(e.g.~\cite {MatsuCarl}), then the shift map $\shft $ is surjective, but there are many non-surjective (one-sided)
subshifts as well.}, but when they are, there is another property (to be described in our next result) even stronger
than {\hcofty }, which is still equivalent to minimality.  It says that one may reach infinite words \"{at their source}
with bounded cost.  This property is essentially the same as the condition appearing in \cite [Theorem 4.20]{Thomsen}.
See section \ref {AppSection} for a more thorough discussion of this condition.

\state Proposition \label SuperHiperCond Let $X$ be a subshift and suppose that the shift map $\shft :X\to X$ is
surjective.  Then the spectral partial action associated to $X$ is minimal if and only if, for every finite subset
$B\subseteq \Lambda ^*$, with $F_B$ nonempty, there exists a constant $M>0$, such that, for every $x$ in $X$, one may
find $\gamma \in \Lambda ^*$ with $|\gamma |\leq M$, and such that $\beta \gamma x\in X$, for all $\beta \in B$.  In
other words, $$ \sup _{x\in X}\min \{|\gamma |: \beta \gamma x\in X, \hbox { for all } \beta \in B\}<\infty .  $$

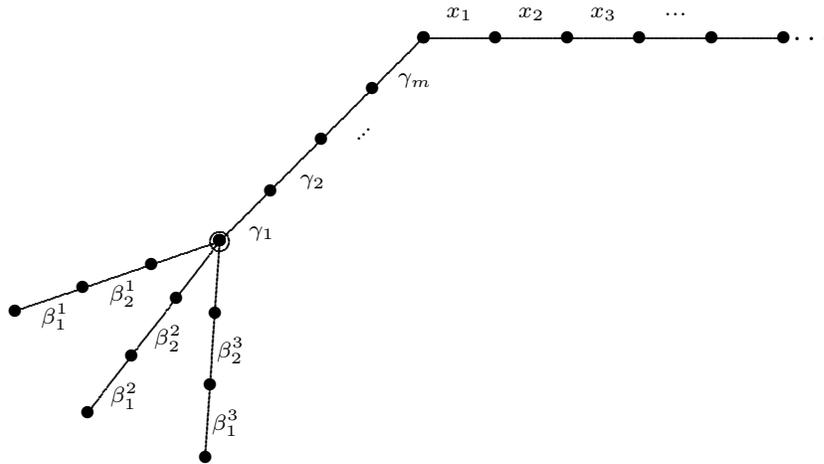
\null \hskip 2cm \beginpicture \setcoordinatesystem units <0.0080truecm, -0.0080truecm> point at 1000 1000 \put
{$\bullet $} at 0 0 \plot 0 0 120 0 / \put {$\bullet $} at 120 0 \put {$\scriptstyle x_1$} at 60 -40 \plot 120 0 240 0 /
\put {$\bullet $} at 240 0 \put {$\scriptstyle x_2$} at 180 -40 \plot 240 0 360 0 / \put {$\bullet $} at 360 0 \put
{$\scriptstyle x_3$} at 300 -40 \plot 360 0 480 0 / \put {$\bullet $} at 480 0 \put {$\scriptstyle \ldots $} at 420 -40
\plot 480 0 600 0 / \put {$\bullet $} at 600 0 \put {$\ldots $ } at 655 0 \put {$\bullet $} at -339 339 \plot -339 339
-255 255 / \put {$\bullet $} at -255 255 \put {$\scriptstyle \gamma _1$} at -269 325 \plot -255 255 -170 170 / \put
{$\bullet $} at -170 170 \put {$\scriptstyle \gamma _2$} at -184 240 \plot -170 170 -85 85 / \put {$\bullet $} at -85 85
\put {$\scriptstyle \nedots $} at -99 156 \plot -85 85 0 0 / \put {$\bullet $} at 0 0 \put {$\scriptstyle \gamma _m$} at
-14 71 \circulararc 360 degrees from -339 355 center at -339 339 \put {$\bullet $} at -680 455 \plot -680 455 -567 417 /
\put {$\bullet $} at -567 417 \put {$\scriptstyle \beta ^1_1$} at -614 464 \plot -567 417 -453 378 / \put {$\bullet $}
at -453 378 \put {$\scriptstyle \beta ^1_2$} at -500 426 \plot -453 378 -339 339 / \put {$\bullet $} at -339 339 \put
{$\scriptstyle \null $} at -387 387 \put {$\bullet $} at -559 625 \plot -559 625 -486 530 / \put {$\bullet $} at -486
530 \put {$\scriptstyle \beta ^2_1$} at -498 596 \plot -486 530 -412 435 / \put {$\bullet $} at -412 435 \put
{$\scriptstyle \beta ^2_2$} at -425 500 \plot -412 435 -339 339 / \put {$\bullet $} at -339 339 \put {$\scriptstyle
\null $} at -352 405 \put {$\bullet $} at -363 699 \plot -363 699 -355 579 / \put {$\bullet $} at -355 579 \put
{$\scriptstyle \beta ^3_1$} at -329 641 \plot -355 579 -347 459 / \put {$\bullet $} at -347 459 \put {$\scriptstyle
\beta ^3_2$} at -321 521 \plot -347 459 -339 339 / \put {$\bullet $} at -339 339 \put {$\scriptstyle \null $} at -313
401 \endpicture \hfill \null

\Proof It is evident that the condition given implies {\hcofty }, and hence minimality by \ref {MainHyper}.  Conversely,
suppose that $X$ is minimal, hence also {\hcof } by \ref {MainHyper}.  Given any finite subset $B$ of $\Lang $, let $$ n
= \sup _{x\in X}\cost (B,x).  $$

Fixing $x$ in $X$, we may use the fact that $\shft $ is surjective to find $z$ in $X$ such that $S^n(z)=x$.  It is then
clear that $z$ is of the form $\mu x$, with $|\mu |=n$.  Spelling out the fact that $z$ may be reached from $B$ with
cost no more than $n$, there are $\alpha $ and $\gamma $ in $\Lambda ^*$, such that $\alpha $ is a prefix of $z$, and,
upon writing $z=\alpha y$, one has that $\beta \gamma y\in X$, for every $\beta $ in $B$, and also $|\alpha |+|\gamma
|\leq n$.

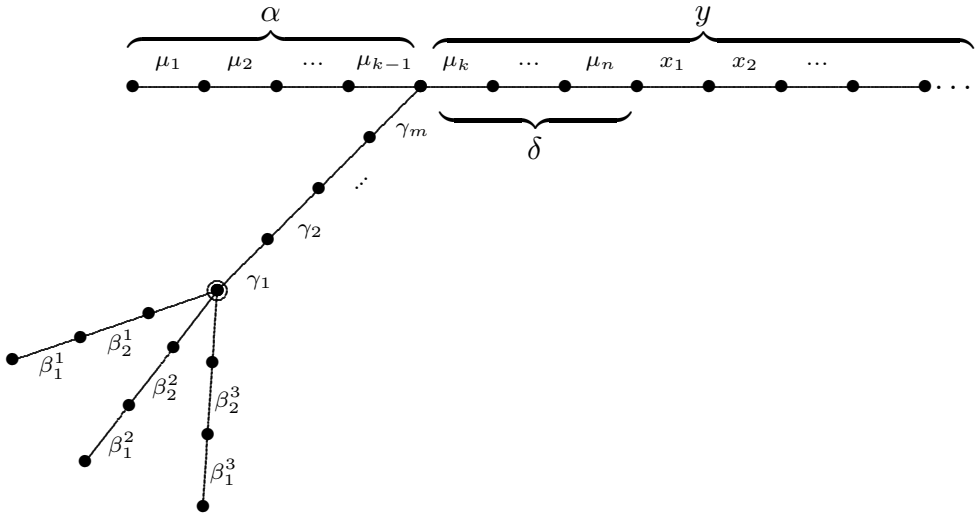
\null \hskip 0.5cm \beginpicture \setcoordinatesystem units <0.0080truecm, -0.0080truecm> point at 1000 1000 \put
{$\bullet $} at 0 0 \plot 0 0 120 0 / \put {$\bullet $} at 120 0 \put {$\scriptstyle \mu _1$} at 60 -40 \plot 120 0 240
0 / \put {$\bullet $} at 240 0 \put {$\scriptstyle \mu _2$} at 180 -40 \plot 240 0 360 0 / \put {$\bullet $} at 360 0
\put {$\scriptstyle \ldots $} at 300 -40 \plot 360 0 480 0 / \put {$\bullet $} at 480 0 \put {$\scriptstyle \mu _{k-1}$}
at 420 -40 \plot 480 0 600 0 / \put {$\bullet $} at 600 0 \put {$\scriptstyle \mu _k$} at 540 -40 \plot 600 0 720 0 /
\put {$\bullet $} at 720 0 \put {$\scriptstyle \ldots $} at 660 -40 \plot 720 0 840 0 / \put {$\bullet $} at 840 0 \put
{$\scriptstyle \mu _n$} at 780 -40 \plot 840 0 960 0 / \put {$\bullet $} at 960 0 \put {$\scriptstyle x_1$} at 900 -40
\plot 960 0 1080 0 / \put {$\bullet $} at 1080 0 \put {$\scriptstyle x_2$} at 1020 -40 \plot 1080 0 1200 0 / \put
{$\bullet $} at 1200 0 \put {$\scriptstyle \ldots $} at 1140 -40 \plot 1200 0 1320 0 / \put {$\bullet $} at 1320 0 \put
{$\ldots $ } at 1375 0 \put {$\overbrace {\hbox to 90pt{\hfill }}^{\textstyle \alpha }$} at 230 -100 \put {$\underbrace
{\hbox to 60pt{\hfill }}_{\textstyle \delta }$} at 670 80 \put {$\overbrace {\hbox to 170pt{\hfill }}^{\textstyle y}$}
at 950 -100 \put {$\bullet $} at 141 339 \plot 141 339 225 255 / \put {$\bullet $} at 225 255 \put {$\scriptstyle \gamma
_1$} at 211 325 \plot 225 255 310 170 / \put {$\bullet $} at 310 170 \put {$\scriptstyle \gamma _2$} at 296 240 \plot
310 170 395 85 / \put {$\bullet $} at 395 85 \put {$\scriptstyle \nedots $} at 381 156 \plot 395 85 480 0 / \put
{$\bullet $} at 480 0 \put {$\scriptstyle \gamma _m$} at 466 71 \circulararc 360 degrees from 141 355 center at 141 339
\put {$\bullet $} at -200 455 \plot -200 455 -87 417 / \put {$\bullet $} at -87 417 \put {$\scriptstyle \beta ^1_1$} at
-134 464 \plot -87 417 27 378 / \put {$\bullet $} at 27 378 \put {$\scriptstyle \beta ^1_2$} at -20 426 \plot 27 378 141
339 / \put {$\bullet $} at 141 339 \put {$\scriptstyle \null $} at 93 387 \put {$\bullet $} at -79 625 \plot -79 625 -6
530 / \put {$\bullet $} at -6 530 \put {$\scriptstyle \beta ^2_1$} at -18 596 \plot -6 530 68 435 / \put {$\bullet $} at
68 435 \put {$\scriptstyle \beta ^2_2$} at 55 500 \plot 68 435 141 339 / \put {$\bullet $} at 141 339 \put
{$\scriptstyle \null $} at 128 405 \put {$\bullet $} at 117 699 \plot 117 699 125 579 / \put {$\bullet $} at 125 579
\put {$\scriptstyle \beta ^3_1$} at 151 641 \plot 125 579 133 459 / \put {$\bullet $} at 133 459 \put {$\scriptstyle
\beta ^3_2$} at 159 521 \plot 133 459 141 339 / \put {$\bullet $} at 141 339 \put {$\scriptstyle \null $} at 167 401
\endpicture \hfill \null

\bigskip Of crucial importance is whether we have reached $z$ before or after the end of the prefix $\mu $.  Observing
that $$ |\alpha | \leq |\alpha |+|\gamma | \leq n = |\mu |, $$ the answer to the above dilemma is before!  Since both
$\alpha $ and $\mu $ are prefixes of $z$, we deduce that $\alpha $ is a prefix of $\mu $, and we may then write $\mu
=\alpha \delta $, for some finite word $\delta $.  It is then clear that $y=\delta x$, and $$ X\ni \beta \gamma y =
\beta \gamma \delta x = \beta \gamma 'x, $$ where $\gamma '=\gamma \delta $.  Noticing that $$ |\gamma '|=|\gamma
|+|\delta | \leq |\alpha |+ |\gamma |+|\mu | \leq 2n, $$ the proof is concluded.  \endProof

\section Applications

\label AppSection In this short section we will draw a few conclusions about the Carlsen-Matsumoto C*-algebra which may
be derived from the work we did so far.  Some of these results are already known, but we may recover them easily given
the many tools available to treat partial crossed product algebras.

\state Theorem {\rm (cf.~\cite [Theorem 17]{CarlsenSilvestrov})} For every subshift $X$, one has that $\CM $ is a
nuclear C*-algebra.

\Proof By \ref {AlgebrasSame} we have that $\CM $ is the reduced crossed product relative to a partial dynamical system
satisfying the approximation property.  The result is then an immediate application of \book {Proposition}{25.10}.
\endProof

The first description of $\CM $ as a groupoid C*-algebra was given in \cite {CarlsenThesis}.  We may recover it here due
to its description as a partial crossed product:

\state Theorem \label AlgAsGpd Let $\G _X$ be the transformation groupoid associated to the spectral partial action for
a given subshift $X$.  Then $\G _X$ is a second countable, Hausdorff, \'etale, amenable groupoid and $$ C^*(\G _X)\simeq
\CM .  $$

\Proof We refer the reader to \cite [Section 2]{AbadieGpd} for the construction of the transformation groupoid relative
to a partial action.  Since the acting group, namely $\F $, is discrete, it is clear that $\G _X$ is \'etale.  It is
also easy to see that $\G _X$ is Hausdorff and, based on the fact that $\Spec $ is metrizable and $\F $ is countable,
$\G _X$ is seen to be second countable.  That $C^*(\G _X)$ is isomorphic to the crossed product $C(\Spec )\rtimes
_\vartheta \F $, and hence also to $\CM $, follows from \cite [Theorem 3.3]{AbadieGpd}.  Finally, since $\CM $ is
nuclear, we deduce from \cite [Theorem 5.6.18]{BrownOzawa} that $\G _X$ is amenable.  \endProof

Our next result involves a weaker alternative to the condition (I) introduced by Matsumoto in \cite
{MatsumotoDimension}.  Under condition (I) the result below has essentially been proved by Matsumoto and Carlsen in
\cite [Lemma 2.3]{MatsuCarl}.  See also \cite [Theorem 16]{CarlsenSilvestrov}.

\state Theorem \label TopFreeFinal Let $X$ be a subshift.  Then the following are equivalent: \izitem \zitem $X$
satisfies \ref {MainTopFree.ii}, that is, for every $\beta _1,\beta _2,\ldots ,\beta _n$ in $\Lambda ^*$, and for every
circuit $\gamma $ such that $\gamma ^\infty \in \bigcap _{i=1}^n\Folow _{\beta _i}$, one has that $\bigcap
_{i=1}^n\Folow _{\beta _i}$ contains some element other than $\gamma ^\infty $, \zitem the spectral partial action
associated to $X$ is topologically free, \zitem $\G _X$ is an essentially principal groupoid (see \cite [Definition
3.1]{RenaultCartan}), \zitem every nontrivial closed two-sided ideal in $\CM $ has a nontrivial intersection with $\Core
$, \zitem every *-homomorphism defined on $\CM $ is injective, provided it is injective on $\Core $.  \medskip \noindent
In case the equivalent conditions above hold, then $\CM $ is naturally isomorphic to $\MatsAlg $.

\Proof The equivalence between (i) and (ii) is precisely the content of \ref {MainTopFree}, while the equivalence
between (ii) and (iii) is evident.  That (ii) implies (iv) was proved in \cite [Theorem 2.6]{ELQ} (see also \cite
[Theorem 4.4]{NHausd}), for the reduced crossed product, which is isomorphic to $\CM $ by \ref {AlgebrasSame}.  A
standard argument equating an ideal with the kernel of a *-homomorphism proves that (iv) and (v) are equivalent.

To close the cycle it suffices to prove that (iv) implies (iii), and this follows from \cite [Proposition 5.5]{BCFS}.

Addressing the last sentence in the statement, it is easy to see that the pair $(i,\pr )$ is a covariant representation
of the spectral partial action in $\MatsAlg $, where $i$ denotes the inclusion of $C(\Spec ) = \Core $ in $\MatsAlg $.
By the universal property of the crossed product \book {Proposition}{13.1} there is an obviously surjective
*-homomorphism $$ \varphi : C(\Spec )\rtimes _\vartheta \F \to \MatsAlg , $$ extending the above covariant pair.  Since
$\varphi $ restricts to the identity map on $C(\Spec )$, it is injective there.  The kernel of $\varphi $ therefore has
trivial intersection with $C(\Spec )=\Core $, whence this kernel itself is trivial by (iv).  Therefore $\varphi $ is
one-to-one and the result follows from the identification between the crossed product and $\CM $ given by \ref
{AlgebrasSame}.  \endProof

It is curious to compare \ref {MainTopFree.ii} with condition (I) from \cite {MatsumotoDimension} already referred to.
Observe that the former is essentially saying that there is no $l$-past equivalence class containing a single
\"{periodic} point, while the latter says that there is no such class containing a single point (periodic or not).  Seen
from this point of view, one realizes that condition (I) is stronger than \ref {MainTopFree.ii}, and therefore \ref
{TopFreeFinal} is a stronger result than \cite [Theorem 16]{CarlsenSilvestrov}, say.  The following example shows that
there is a real difference between the above conditions.

\state Example \rm There exists a subshift satisfying condition \ref {MainTopFree.ii} but not condition (I).

\Proof Let $\Sigma =\{1,2\}$ and let $$ Y=\Sigma ^{\bf N} = \{1,2\}^{\bf N} $$ be the full shift.  Choose any
non-periodic point $z$ in $Y$, and let us construct another subshift on the new alphabet $\Lambda =\{0,1,2\}$.  For this
we let $\Forbid $ be the following set of forbidden words: $$ \Forbid =\{a0:a\in \Lambda \}\cup \{0\alpha : \alpha \in
Y,\ \alpha \hbox { is \underbar {not} a prefix of }z\}, $$ and we claim that the subshift $$ X=Y_\Forbid $$ satisfies
the requirements in the statement.

Since all forbidden words have a $\qt 0$ somewhere, any infinite word not involving $\qt 0$ is allowed, meaning that
$Y\subseteq X$.  Another interesting element of $X$ is the word $0z$, which narrowly escapes being ruled out!  Other
than that, there is nothing else in $X$, meaning that $$ X=Y\cup \{0z\}, $$ as the reader may easily verify.  Denoting
by $$ \Lambda _l(x) = \{\alpha \in \Lambda ^*: |\alpha |=l,\ \alpha x\in X\}, \for x\in X, \for l\in {\bf N}, $$ notice
that $$ \def \quad { } \matrix { \Lambda _1(z) & = & \{0,1,2\},&\hbox {and} \cr \pilar {18pt} \Lambda _1(x) & = &
\{\hfill 1,2\}, & \ \forall x\neq z.}  $$ So the \"{$1$-past equivalence class} of $z$, in the sense of \cite
{MatsumotoDimension}, consists only of $z$, and we thus see that $X$ does not satisfy condition (I).

In order to prove that $X$ satisfies condition \ref {MainTopFree.ii}, let $\gamma $ be a circuit in $X$, and suppose
that $\gamma ^\infty $ lies in $\Folow _B$, for some finite set $B$ of finite words.  Since we have chosen $z$ not
periodic, we necessarily have that $z\neq \gamma ^\infty $.

Given $\beta $ in $B$ we then have that $\beta \gamma ^\infty $ is in $X$, so $\beta $ cannot involve $\qt 0$.  It then
follows that $x\in \Folow _\beta $, for every $x$ in $Y$, whence $Y\subseteq \Folow _B$.  So there are plenty of
elements in $\Folow _B$ other than $\gamma ^\infty $ to choose from, proving that $X$ satisfies condition \ref
{MainTopFree.ii}.  \endProof

Our next Theorem is related to various known simplicity results for C*-algebras of subshifts.

\state Theorem \label Simplicity Let $X$ be a subshift.  Then the following conditions are equivalent: \izitem \zitem
$\CM $ is simple, \zitem $X$ is both {\colcof } and {\scof }, and it satisfies condition \ref {MainTopFree.ii}, \zitem
$X$ is {\hcof } and contains at least one non eventually periodic point.

\Proof \ipfimply \pfimply (i)(ii) By \ref {AlgAsGpd} and \cite [Theorem 5.1]{BCFS} one has that $\G _X$ is an
essentially principal and minimal groupoid.  It then follows that the spectral partial action is topologically free and
minimal, so the conclusion follows from \ref {MainTopFree} and \ref {MainMinimal}.

\pfimply (ii)(i) The spectral partial action is topologically free by \ref {MainTopFree}, and minimal by \ref
{MainMinimal}.  Therefore the reduced crossed product $\Core \rtimes ^{\hbox {\sixrm red}}_\tau \F $ is simple by \cite
[Corollary 2.9]{ELQ}, and hence so is $\CM $ by \ref {AlgebrasSame}.

\BreakImply (ii)$\,\Leftrightarrow \,$(iii) We have seen in \ref {MainHyper} that $X$ is both {\colcof } and {\scof } if
and only if it is {\hcof }.  In this case the existence of a non eventually periodic point is equivalent to topological
freeness of the spectral partial action by \ref {TopFreeWhenColCof}, hence also to condition \ref {MainTopFree.ii}.
\endProof

Of course the last part of \ref {Simplicity.ii} could be interchanged with the last part of \ref {Simplicity.iii}, the
above choice happening to be purely personal.

Similar simplicity results for C*-algebras associated to subshifts are to be found in many works in the literature, such
as \cite [Corollary 6.11]{MatsumotoDimension}, \cite [Proposition 2.6]{MatsuCarl}, \cite [Proposition 4.3]{AdelaR}, and
\cite [Theorem 4.20]{Thomsen}.

We nevertheless feel that the result above is worth the trouble of proving it, mainly because it gives necessary and
sufficient conditions for simplicity.  These conditions also have a distinctively geometrical flavor at the same time
that they are closer to the cofinality hypothesis of classical results about graph algebras, especially if compared with
the intricate condition of ``\"{irreducibility in $l$-past equivalence}" present in \cite {MatsumotoDimension} and \cite
{MatsuCarl}.

Regarding \cite [Proposition 4.3]{AdelaR}, one should be aware that it only holds for subshifts of finite type because
otherwise the shift map is not open, as we have seen in \ref {SFTOpen}.  In particular the statement that $\CM $ is
simple when $X$ is the even shift, given in the second paragraph of page 222, is not correct.  As we have seen in \ref
{EvenNotMinimal}, the spectral partial action for the even shift is not minimal, so it fails to be {\hcof } by \ref
{MainHyper}, whence $\CM $ is not simple by \ref {Simplicity}.

Among the references cited above, the only one providing for necessary and sufficient conditions is \cite [Theorem
4.20]{Thomsen}, which however only applies to surjective subshifts.  The condition given there is essentially the same
as the one described in \ref {SuperHiperCond}, hence equivalent to minimality in the surjective case, but it is much too
strong in general, even for Markov subshifts.  For example, the Markov subshift on the alphabet $\Lambda =\{1,2,3\}$,
corresponding to the graph

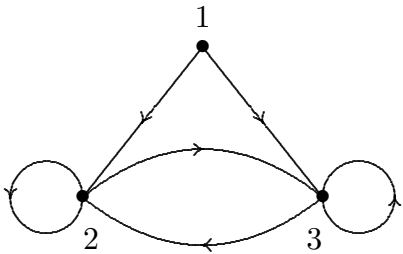
\null \hfill \beginpicture \setcoordinatesystem units <0.0080truecm, -0.0080truecm> point at 1000 1000 \put {$\bullet $}
at 200 -250 \put {$1$} at 200 -300 \put {$\bullet $} at 0 0 \put {$2$} at 14 69 \put {$\bullet $} at 400 0 \put {$3$} at
386 69 \setquadratic \plot 200 -250 100 -125 0 0 / \setlinear \arrow <0.11cm> [0.3,1.2] from 101 -127 to 99 -123
\setquadratic \plot 200 -250 300 -125 400 0 / \setlinear \arrow <0.11cm> [0.3,1.2] from 299 -127 to 301 -123
\setquadratic \plot 0 0 200 -80 400 0 / \setlinear \arrow <0.11cm> [0.3,1.2] from 198 -80 to 202 -80 \setquadratic \plot
400 0 200 80 0 0 / \setlinear \arrow <0.11cm> [0.3,1.2] from 202 80 to 198 80 \circulararc 360 degrees from 0 0 center
at -60 0 \arrow <0.11cm> [0.3,1.2] from -120 -2 to -120 2 \circulararc 360 degrees from 400 0 center at 460 0 \arrow
<0.11cm> [0.3,1.2] from 520 2 to 520 -2 \endpicture \hfill \null

\bigskip \noindent (recall that our subshifts consist of infinite paths formed by vertices, as opposed to edges) is
non-surjective: since $\qt 1$ is a \"{source}, any infinite word beginning with $\qt 1$ is not in the range of $\shft $.
The associated Carlsen-Matsumoto algebra (which in this case is well known to be the Cuntz-Krieger algebra for the
adjacency matrix of the above graph) is simple by \ref {Simplicity}, as the reader may easily verify.  However its
associated subshift does not satisfy Thomsen's condition, namely the condition given in \ref {SuperHiperCond}, e.g.~for
$$ B = \{\qt 2\}, \and x={12222\ldots }, $$ simply because $x$ is not in any follower set whatsoever.

\references

\Article AbadieGpd F. Abadie; On partial actions and groupoids; Proc. Amer. Math. Soc., 132, (2003), 1037-1047

\Article AdelaR C. Anantharaman-Delaroche; Purely infinite C*-algebras arising form dynamical systems;
Bull. Soc. Math. France, 125 (1997), 199-225

\Article ArchSpiel R. J. Archbold and J. Spielberg; Topologically free actions and ideals in discrete dynamical systems;
Proc.  Edinburgh Math. Soc., 37 (1993), 119-124

\Article BCFS J. Brown, L. O. Clark, C. Farthing and A. Sims; Simplicity of algebras associated to \'etale groupoids;
Semigroup Forum, 88 (2014), 433-452

\Bibitem BrownOzawa N.P. Brown and N. Ozawa; C*-algebras and finite-dimensional approximations; Graduate Studies in
Mathematics, 88, American Mathematical Society, Providence, RI, 2008

\Bibitem CarlsenThesis T. M. Carlsen; Operator algebraic applications in symbolic dynamics; Ph.D. thesis, University of
Copenhagen (2004)

\Article CarlsenCuntzPim T. M. Carlsen; Cuntz-Pimsner C*-algebras associated with subshifts; Internat. J. Math., 19
(2008), 47-70

\Bibitem CarlsenNotes T. M. Carlsen; C*-algebras associated to shift spaces; Notes for the Summer School Symbolic
Dynamics and Homeomorphisms of the Cantor set, University of Copenhagen, 23--27 June 2008

\Article MatsuCarl T. M. Carlsen and K. Matsumoto; Some remarks on the C*-algebras associated with subshifts;
Math. Scand., 95 (2004), 145-160

\Article CarlsenSilvestrov T. M. Carlsen and S. Silvestrov; C*-crossed products and shift spaces; Expo. Math., 25
(2007), 275-307

\Article CuntzKrieger J. Cuntz and W. Krieger; A class of C*-algebras and topological Markov chains; Invent. Math., 63
(1981), 25-40

\Article amena R. Exel; Amenability for Fell bundles; J. reine angew. Math., 492 (1997), 41-73

\Article ortho R. Exel; Partial representations and amenable Fell bundles over free groups; Pacific J. Math., 192
(2000), 39-63

\Article NewLook R. Exel; A new look at the crossed-product of a C*-algebra by an endomorphism; Ergodic Theory
Dynam. Systems, 23 (2003), 1733-1750

\Article NHausd R. Exel; Non-Hausdorff \'etale groupoids; Proc. Amer. Math. Soc., 139 (2011), 897-907

\Bibitem PDSFB R. Exel; Partial Dynamical Systems, Fell Bundles and Applications; to be published in a forthcoming NYJM
book series.  Available from {\tt http://mtm.ufsc.br/$\sim $exel/papers/pdynsysfellbun.pdf}

\Article infinoa R. Exel and M. Laca; Cuntz-Krieger algebras for infinite matrices; J. reine angew. Math., 512 (1999),
119-172

\Article ELQ R. Exel, M. Laca and J. Quigg; Partial dynamical systems and C*-algebras generated by partial isometries;
J. Operator Theory, 47 (2002), 169-186

\Article ItoTaka S. Ito and Y. Takahashi; Markov subshifts and realization of $\beta $-expansions; J. Math. Soc. Japan,
26 (1974), 33-55

\Article KPR A. Kumjian, D. Pask and I. Raeburn; Cuntz-Krieger algebras of directed graphs; Pacific J. Math, 184 (1998),
161-174

\Article KPRR A. Kumjian, D. Pask, I. Raeburn, and J. Renault; Graphs, groupoids, and Cuntz-Krieger algebras;
J. Funct. Anal., 144 (1997), 505-541

\Bibitem Kurka P. K${\buildrel {\scriptscriptstyle \circ } \over {\hbox {u}}}$rka; Topological and symbolic dynamics;
Cours Sp\'ecialis\'es, 11. Soci\'et\'e Math\'e\-matique de France, Paris, 2003. xii+315 pp

\Bibitem LindMarcus D. Lind and B. Marcus; An introduction to symbolic dynamics and coding; Cambridge Univ. Press, 1999

\Article MatsuOri K. Matsumoto; On C*-algebras associated with subshifts; Internat. J. Math., 8 (1997), 357-374

\Article MatsumotoDimension K. Matsumoto; Dimension groups for subshifts and simplicity of the associated C*-algebras;
J. Math. Soc. Japan, 51 (1999), 679-698

\Article MatsuAutomorph K. Matsumoto; On automorphisms of C*-algebras associated with subshifts; J. Operator Theory, 44
(2000), 91-112

\Article Matsumoto K. Matsumoto; Stabilized C*-algebras constructed from symbolic dynamical systems; Ergodic Theory
Dynam. Systems, 20 (2000), 821-841

\Article Parry W. Parry; Symbolic dynamics and transformations of the unit interval; Trans. Amer. Math. Soc., 122
(1966), 368-378

\Article RenaultCartan J. Renault; Cartan subalgebras in C*-algebras; Irish Math. Soc. Bull., 61 (2008), 29-63

\Bibitem Royer D. Royer; Representa\c c\~oes parciais de grupos; Master Thesis, Universidade Federal de Santa Catarina, 2001

\Article Thomsen K. Thomsen; Semi-\'etale groupoids and applications; Ann. Inst. Fourier (Grenoble), 60 (2010), 759-800

\endgroup \bye